\documentclass[12pt]{amsart}

\usepackage{stmaryrd}

\usepackage{ulem, mathabx}

\def\cal{\relax}

\usepackage{hyperref}

\def\cal{\relax}
\usepackage{DLdef1}
\providecommand{\See}{{\mathbb  S}} 

\makeatletter
\theoremstyle{thm}
\@ifundefined{Theorem}{}{}
\@ifundefined{Lemma}{}{}
\@ifundefined{Statement}{}{}
\@ifundefined{Corollary}{}{}
\@ifundefined{Proposition}{}{}
\theoremstyle{definition}
\@ifundefined{Remark}{}{}
\@ifundefined{Remarks}{}{}
\@ifundefined{Exercise}{}{}
\@ifundefined{Definition}{}{}
\makeatother

\let\ssec\subsection

\renewcommand {\ssbegin}[2][*]
 {\refstepcounter{subsection}%
\if#1*
\addcontentsline{toc}{subsection}{\thesubsection.\hskip 1pc #2}%
\else
\addcontentsline{toc}{subsection}{\thesubsection.\hskip 1pc #2. #1}%
\fi
 \def \secno {\gdef \secno {}{\ssecfont
\thesubsection.\hskip 2ex}%
 }%
 \begin{#2}}

\renewcommand {\sssbegin}[2][*]
 {\refstepcounter{subsubsection}
\if#1*
\addcontentsline{toc}{subsubsection}{\thesubsubsection.\hskip 1pc #2}%
\else
\addcontentsline{toc}{subsubsection}{\thesubsubsection.\hskip 1pc #2. #1}
\fi
 \def \secno {\gdef \secno {}{\ssecfont \thesubsubsection.\hskip 2ex}%
 }%
 \begin{#2}}

\renewcommand {\parbegin}[2][*]
 {\refstepcounter{paragraph}
\if#1*
\addcontentsline{toc}{paragraph}{\theparagraph.\hskip 1pc #2}%
\else
\addcontentsline{toc}{paragraph}{\theparagraph.\hskip 1pc #2. #1}
\fi
 \def \secno {\gdef \secno {}{\ssecfont \theparagraph.\hskip 2ex}%
 }%
 \begin{#2}}

\makeatletter
\@ifundefined{widecheck}{\DeclareRobustCommand{\widecheck}[1]{\check{#1}}}{}
\makeatother

\providecommand{\ssssec}[1]{\paragraph{#1.}}

\providecommand{\Mat}{\mathop{\mathrm{Mat}}\nolimits}
\providecommand{\Odd}{\mathop{\mathrm{Odd}}\nolimits}
\providecommand{\Q}{\mathop{\mathrm{Q}}\nolimits}
\providecommand{\GQ}{\mathop{\mathrm{GQ}}\nolimits}
\providecommand{\rd}{\mathop{\mathrm{rd}}\nolimits}
\providecommand{\id}{\mathop{\mathrm{id}}\nolimits}

\newcommand{\del}{\partial}

\newcommand {\lra} {\longleftrightarrow}

\begin{document}

\title{Superization of K\"ahler and hyper-K\"ahler manifolds}

\author{Dimitry Leites}

\address{Department of Mathematics, Stockholm University, Albanov\"agen 28, SE-114 19 Stockholm, Sweden; dimleites@gmail.com}


\begin{abstract} Let $J$ be an (almost) complex structure on a manifold $\mathcal{M}$; let~ $h$ be a~pseudo-Hermitian with respect to $J$ non-degenerate symmetric form on the tangent bundle over~$\mathcal{M}$.  The manifold $\mathcal{M}$ with the triple $(\omega, h, J)$, where  $\omega$ is an almost symplectic differential 2-form defined through $h$ and $J$, is called (almost) K\"ahler. If $\omega$ and $J$ are \textit{parallel}, i.e., preserved by a~torsion-free connection, the adjective ``almost\rq\rq\ can be dropped; that connection is then unique, but it need not be flat: a~K\"ahler manifold has curvature, and no coordinates bring the whole triple to a~standard shape. I also describe the local invariants of a pair of bilinear forms and consider what superizing each of the ingredients $(\mathcal{M}; \omega, h, J)$ brings about. The same program is performed with the definition of hyper-K\"ahler manifolds, where there are three $J$s, and three $\omega$s. Replacing  $\omega$ with the dual bivector  makes  the definition of (almost) K\"ahler (and hyper-K\"ahler) supermanifold depend on a supervariety. The Lie superalgebras $\mathfrak{g}$ of invariant operators acting on the spaces of (pseudo)differential and (pseudo)integral forms are described in all cases that arise. 
Lee's theorem on the two conditions $d\omega=0$ splits into is superized. Applications of $\mathfrak{g}$ to representation theory (analogs of the Hodge-Lepage decomposition) will be given elsewhere, together with interpretation from Howe's duality point of view. \end{abstract}


\makeatletter
\@namedef{subjclassname@2020}{\textup{2020} Mathematics Subject Classification}
\makeatother

\maketitle

\markboth{\itshape Dimitry Leites} {{\itshape Superization of K\"ahler and hyper-K\"ahler manifolds}}

\thispagestyle{empty}

\section{Definitions of (almost) K\"ahler  (super)manifolds}\label{sAlmK}

For details on Lie superalgebras needed below, and further basics in super setting, see \cite{Del, MaGF} updated in \cite{Lo, LSoS}; the ground field is $\Cee$ or~$\Ree$. In this article, I abbreviate ``manifold or supermanifold" to (super)manifold; all (super)manifolds considered are real and smooth.

\ssec{The main objects and tools}\label{ss3} Recall that a supermanifold $\cM$ is a~ringed space which locally is of the form $(M, \cL_{\Eee^{\bcdot}(E)})$, where $\cL_{W}$ designates the sheaf of  sections of the vector bundle $W$ over the manifold $M$, and $E$ is a~vector bundle with fiber $F$, with $\Eee^k$ denoting the operator of raising to the $k$-th super exterior power, $\Eee^{\bcdot}:=\oplus_{k\geq 0}\, \Eee^{k}$. The supermanifold  $(M, \cL_{\Eee^{\bcdot}(E)})$ is said to be ``associated with the bundle~$E$".
(Recall that every  supermanifold, i.e., a smooth supervariety, is of this form even globally, see \cite[item 4.1.3]{MaGF}.) The following are main objects in our constructions. 

\textbf{(i)  Supermanifolds}. Let $\cM=(M, \cO_\cM)$ be a~supermanifold with the structure sheaf $\cO_\cM:=\cL_{\Eee^{\bcdot}(E)}$, the local coordinates being $x:=(u, \theta)$, where $u:=(u_1, \dots, u_a)$ are local coordinates on the $a$-dimensional manifold $M$ and $\theta:=(\theta_1, \dots, \theta_b)$ are basis elements of the $b$-dimensional fiber~ $F$ of the bundle $E$ over~$M$. The pair $(a,b)$ is called the (super)\textit{dimension} of $\cM$, sometimes written as~$\cM^{a,b}$. Let $\cF:=\cF(\cM)$ designate the superalgebra of functions on~$\cM$.

\textbf{(ii) (Super)manifold with a hat}. Let $\widehat \cM=(\cM, \widehat \Omega(\cM))$ be the supermanifold with structure sheaf
$\widehat \Omega(\cM)$ whose sections are called \textit{pseudo-differential} forms, i.e.,  functions in $\widehat x$, where $\widehat x_i:=dx_i$ for all $i$, and $d$ is the exterior differential considered as an odd (parity changing) operator. This $d$ is a derivation of $\widehat \cF:=\cF(\widehat \cM)$, i.e., a vector field on~$\widehat \cM$. In the particular case ${b=0}$, \textbf{all} functions on $\widehat \cM$ are polynomial in the local coordinates~$\widehat x$. Obviously, $\widehat \cM$ is the supermanifold associated with the cotangent bundle $\Pi(T^*\cM)$ over $\cM$, where $\Pi$ is the \textit{inversion of parity} functor applied fiber-wise. Sometimes, one confines oneself to a~version ${{}^p\widehat \cM=(M, \Omega(\cM))}$ of $\widehat \cM$ whose functions, called \textit{differential} forms, are \textit{polynomials} in~ $\widehat x$; the superspace $\Omega^{\bcdot}(\cM)$ is $\Zee$-graded by the degree with respect to powers of~ $\widehat x$, where $\deg \widehat x_i=1$ for all $i$; clearly, $\Omega^k(\cM)$ is an $\cF(\cM)$-module for every~$k$.

\textbf{(iii) (Super)manifold with a check}. Let $\widecheck \cM=(\cM, \widecheck \Sigma)$ be the supermanifold associated with $\Pi(T\cM)$, the tangent bundle over $\cM$ to which $\Pi$ is applied fiber-wise. The  sections of the structure sheaf $\widecheck \Sigma$ of $\widecheck \cM$ are called \textit{pseudo-integral} forms, in local coordinates they depend on $\widecheck x_i:=\Pi(\partial_{x_i})$. In the particular case $b=0$, \textbf{all} functions on $\widecheck \cM$, called \textit{multi-vector fields}, are polynomials in~$\widecheck x$. Sometimes, one confines oneself to a~version ${}^p\widecheck \cM$ of $\widecheck \cM$ whose functions are polynomials in~ $\widecheck x$, see \cite{BLi}. 

\textbf{(iv) An odd analog of the symplectic form}. The non-degenerate even differential 2-form $\omega_\ev$ is called (almost) \textit{symplectic} although, strictly speaking, on supermanifolds it should be called (almost) \textit{ortho-symplectic}, but this sounds too long.   In 1989, A.~Weil suggested to me to call the non-degenerate \textbf{odd} non-degenerate bilinear form  \textit{periplectic}. Accordingly, the supermanifold with the non-degenerate  odd differential 2-form $\omega_\od$ will be called (almost) \textit{periplectic} and the non-integrable distribution singled out by the even analog of the contact 1-form with an ``odd time" $\tau$ and the curvature form proportional to the periplectic form will be called \textit{pericontact}.

\textbf{(v) Contact and pericontact Lie (super)algebras and their subalgebras}. For their deformations, see Section~\ref{S2} and \cite{LSh3}.

$\bullet$ \underline{Odd form $\alpha_1$, even form $\omega_\ev$}. Let $\omega_\ev:=d\alpha_1$, where for the normal form of $\alpha_1$ we take
\begin{equation}
\label{a1} \alpha_1=dt+\mathop{\sum}\limits_{1\leq i\leq
n}(p_idq_i-q_idp_i)+ \mathop{\sum}\limits_{1\leq j\leq
k}(\xi_jd\eta_j+\eta_jd\xi_j)+
\begin{cases}0,&
\text{ if }\ m=2k\\
\mathop{\sum}\limits_{1\leq s\leq
\ell}\eps_s\theta_sd\theta_s,&\text{ if }\ m=2k+\ell,\end{cases} 
\end{equation}
where $\eps_s=\pm1$ for all $s$ (over $\Cee$, one can take $\eps_s=1$ for all~$s$). For any $f\in\Cee [t, p,
q, \xi, \eta, \theta]$ with indeterminates $t, p,
q$ (even) and $\xi, \eta, \theta$ (odd), set
\begin{equation}
\label{K_f} 
 K_f:=(2-E)(f)\pder{t}-H_f + \pderf{f}{t} E,
\end{equation}
where $E:=\mathop{\sum}_j y_j \pder{y_{j}}$ for all the indeterminates $y_{j}$ except $t$  is the \textit{Euler operator}, 
and $H_f$ is the Hamiltonian field with Hamiltonian $f$:
\begin{equation}
\label{H_f} H_f:=\mathop{\sum}\limits_{i\leq n}\left
(\pderf{f}{p_i} \pder{q_i}-\pderf{f}{q_i} \pder{p_i}\right)
-(-1)^{p(f)}\left(\mathop{\sum}\limits_{j\leq
k}\left(\pderf{f}{\xi_j} \pder{\eta_j}+ \pderf{f}{\eta_j}
\pder{\xi_j}\right)+ \mathop{\sum}\limits_{j\leq \ell}\eps_j\pderf{ f}{\theta_j}
\pder{\theta_j}\right).
\end{equation}

Recall notation of Lie (super)algebras and their main properties in terms of  the \textit{Lie derivative $L_D$ along the vector field~ $D$}:
\be\label{Kand}
\begin{array}{ll}
\fk(2n+1|m):=\Span(K_f\mid f\in\Cee [t, y])&\text{the contact algebra: $L_{K_f}(\alpha_1)=2\del_t(f)\alpha_1$;}\\
\fpo(2n|m):=\Span(K_f\mid L_{K_f}(\alpha_1)=0)&\text{the Poisson algebra;}\\
\fh(2n|m):=\Span(H_f\mid f\in\Cee [y])&\text{the Hamiltonian algebra: $L_{H_f}(\omega_\ev)=0$.}\\
\end{array}
\ee

 $\bullet$ \underline{Even form $\alpha_0$, odd form $\omega_\od$}. Let 
 \begin{equation}
\label{alpha0}
\alpha_0:=d\tau+\mathop{\sum}\limits_{1\leq i\leq n}(\xi_idq_i+q_id\xi_i),\qquad \omega_\od:=d\alpha_0.
\end{equation} 
For any $f\in\Cee [q,
\xi, \tau]$ with indeterminates $
q$ (even) and $\xi, \tau$ (odd), set
\begin{equation}
\label{2.3.3} M_f:=(2-\cE)(f)\pder{\tau}- Le_f -(-1)^{p(f)}
\nfrac{\partial f}{\partial\tau} \cE, 
\end{equation}
where $\cE:=\mathop{\sum}_iy_i \pder{y_i}$ for all the indeterminates $y_{i}$  except $\tau$, and 
\begin{equation}
\label{2.3.4} Le_f:=\mathop{\sum}\limits_{i\leq n}\left(
\pderf{f}{q_i}\ \pder{\xi_i}+(-1)^{p(f)} \pderf{f}{\xi_i}\
\pder{q_i}\right).
\end{equation}

The \textit{divergence} of the vector field $D=\sum f_i\partial_{u_i}+\sum g_j\partial_{\theta_j}$ in coordinates $x:=(u, \theta)$ is
\be
\label{DIV} \Div \left (\sum f_i\partial_{u_i}+\sum g_j\partial_{\theta_j}\right):=\sum \partial_{u_i}(f_i)+\sum (-1)^{p(g_j)}\partial_{\theta_j}(g_j).
\ee
Recall the names of several more Lie superalgebras  and their main properties in terms of~ $L_D$.
\be\label{Mand}
\begin{array}{ll}
\fm(n):=\Span(M_f\mid f\in \Cee [
q, \xi; \tau])&\text{pericontact: $L_{M_f}(\alpha_0)=(-1)^{p(f)}2\del_\tau(f)\alpha_0$;}\\
\fb(n):=\Span(M_f\mid L_{M_f}(\alpha_0)=0)&\text{Buttin;}\\
\fle(n):=\Span(Le_f\mid f\in \Cee [
q, \xi])&\text{no fixed name yet: $L_{Le_f}(\omega_\od)=0$;}\\
\fsm (n):=\{M_f\in \fm (n)\mid \Div\ M_f=0\}&\text{special pericontact}; \\
 \fs\fb (n):=\{M_f\in
\fb (n)\mid \Div\ D=0\}&\text{special Buttin};\\
\fsle(n):=\Span(Le_f\mid \Div Le_f=0)&\text{no fixed name yet}.\\
\end{array}
\ee
 
Consider $\fk(2n+1|m)$ and $\fm(k)$ in their \textit{standard grading}:  
\be\label{standgr}
\text{the degree of the ``time" coordinate ($t$ or $\tau$) is 2, that of other coordinates is ~ 1.}
\ee 
Then, the \textit{anti-bracket superalgebra} ($\fb\fa$ being read backwards)
\be\label{ba}
\fb\fa(k):=\fm(k)_{-2}\oplus\fm(k)_{-1} 
\ee
is an analog of the \textit{Heisenberg Lie (super)algebra}
\be\label{hei}
\fhei(2n|m):=\fk(2n+1|m)_{-2}\oplus\fk(2n+1|m)_{-1}.
\ee

$\bullet$ \underline{Brackets}. To the (super)commutators of vector fields $[K_f, K_g]$ and $[M_f, M_g]$,
there correspond the contact and pericontact \textit{brackets} of the generating functions:
\[
\renewcommand{\arraystretch}{1.4}
\begin{array}{l}
{}[K_f, K_g]=K_{\{f, \; g\}_{k.b.}}\ \ \text{~~and~~}\ \ 
 {}[M_f, M_g]=M_{\{f, \;
g\}_{m.b.}}.\end{array}
\]
To give the explicit formulas for these brackets, we first define the brackets on functions that do not depend on $t$
(resp.,~$\tau$).

The \textit{Poisson bracket} $\{-,-\}_{P.b.}$ in the
realization with the form 
$\omega_\ev=d\alpha_1$ for $m=2k+1$ is
given by the formula
\begin{equation}
\label{pb} 
\renewcommand{\arraystretch}{1.4}
\begin{array}{ll} 
 \{f, g\}_{P.b.}&:=\mathop{\sum}\limits_{i\leq n}\ \bigg( \pderf{f}{p_i}\
\pderf{g}{q_i}-\ \pderf{f}{q_i}\
\pderf{g}{p_i}\bigg)-\\ 
&(-1)^{p(f)}\bigg(\mathop{\sum}\limits_{j\leq m}\left(
 \pderf{f}{\xi_j}\ \pderf{ g}{\eta_j}+\pderf{f}{\eta_j}\
\pderf{ g}{\xi_j}\right)+ \pderf{f}{\theta}\ \pderf{
g}{\theta}\bigg)\text{ for any }f, g\in \Cee [p, q, \xi, \eta,
\theta].
\end{array}
\end{equation}

The \textit{Buttin bracket} $\{-,-\}_{B.b.}$   is given by the formula
\begin{equation}
\label{B.b.} \{ f, g\}_{B.b.}:=\mathop{\sum}\limits_{i\leq n}\
\bigg(\pderf{f}{q_i}\ \pderf{g}{\xi_i}+(-1)^{p(f)}\
\pderf{f}{\xi_i}\ \pderf{g}{q_i}\bigg)\text{ for any }f, g\in \Cee
[q, \xi].
\end{equation}

The \textit{contact}  bracket is
\begin{equation}
\label{2.3.8} \{ f, g\}_{k.b.}:=(2-E) (f)\pderf{g}{t}-\pderf{f}
{t}(2-E) (g)-\{ f, g\}_{P.b.}
\end{equation}
and \textit{pericontact} bracket is
\begin{equation}
\label{2.3.9} \{ f, g\}_{m.b.}:=(2-\cE)
(f)\pderf{g}{\tau}+(-1)^{p(f)} \nfrac{\partial f}{\partial\tau}(2-\cE) (g)-\{ f,
g\}_{B.b.}. 
\end{equation}

\textbf{(vi) The Cartan prolongs}. \label{2.6}
Recall that a~ graded Lie superalgebra
$\fb=\oplus_{k \geqslant -d}\ \fb_k$
is said to be \textit{transitive}  if for all $k\geq 0$ we have
\begin{equation*}
\label{12}
\{x\in \fb_k \mid [x,\fb_{-}]=0\}=0, \text{~~where $\fb_{-}:=\oplus _{k <0}\ \fb_k$}.
\end{equation*}

Let $\fg_-=\mathop{\oplus}_{-d\leq i\leq -1}\fg_i$ be a~nilpotent
$\Zee$-graded Lie (super)algebra and $\fg_0$ a~Lie sub(super)al\-geb\-ra of the Lie
(super)algebra $\mathfrak{der}_0(\fg_-)$ of degree-$0$ derivations of~$\fg_-$. The maximal transitive $\Zee$-graded Lie (super)algebra whose non-positive part is $\fg_-\oplus \fg_0$ is called the (generalized) \textit{Cartan prolong} --- the result of Cartan prolongation --- of the pair $(\fg_-,\fg_0)$, and is denoted by $(\fg_-,\fg_0)_*$.

\textbf{(vii) Volume forms, differential and integral forms}.\label{ssInt} Recall that $\Vol:=\cF\vvol(x)$, the space  of \textit{volume forms} on the (super)manifold $M$ with coordinates $x$, is a~ rank-$1$ module with the generator $\vvol(x)$ over the algebra $\cF$ of functions; the Lie derivative $L_D$ along the vector field $D$ acts on $\Vol$ as follows 
(for a~ definition of the divergence, see eq.~\eqref{DIV})
\[
L_D(f\vvol(x)):=\left(D(f)+ (-1)^{p(f)p(D)} f \cdot \Div(D)\right)\vvol(x)\  \text{for any function $f\in\cF$.}
\]

There are two natural ways to define the parity of $\vvol(x)$ on the $(a|b)$-dimensional supermanifold, see \cite{MaGF}: either by setting, as we do here,
$p(\vvol(x)):=b\pmod 2$ or by setting $p(\vvol(x)):=a-b\pmod 2$. We sometimes write just $\vvol$ instead of~$\vvol(x)$.

The space
$\Omega^\ell$ of \textit{differential} $\ell$-forms is the $\cF$-module with the basis of homogeneous degree-$\ell$ polynomials in the $\widehat x_i$, while the space
$\Sigma_{-k}$ of \textit{integral} $(-k)$-forms is the $\cF$-module with the basis of homogeneous degree-$k$ $\Vol$-valued polynomials in $\del_{\widehat x}$. 

\textbf{(viii) A non-degenerate anti-symmetric 2-vector $B\in\See^2(\Pi(\cT\cM))=\Eee^2(\cT\cM)$}, where $\See^i$ is the operator of raising to the $i$th super symmetric power, such that 
\begin{equation}\label{bb} 
\{B,B\}_{B.b.}=0, \ \ \text{see eq.~\eqref{B.b.}}.
\end{equation}
 
\textbf{(ix) The operator of an almost complex structure}. A non-degenerate tensor $J$ of valency $(1,1)$ on the tangent space $T_\star\cM$ at every point~ $\star$ of a real (super)manifold is called (the operator of) an \textit{almost complex structure} if $J$ is even and $J^2=-\id$, where $\id$ is the identity operator, or an \textit{almost queer structure} if $J$ is odd and either $J^2=\id$ or~$J^2=-\id$.

\textbf{(x) A pseudo-Hermitian with respect to  the operator of an almost complex (or queer) structure $J$ non-degenerate symmetric\footnote{In the literature, one often reads about \lq\lq supersymmetric" bilinear forms; for an explanation why there are no such forms and where the confusion (I had also been guilty of sharing) comes from, see \cite{KLLS, Lo, LSh3}.} 2-form}. A non-degenerate symmetric 2-form $h\in \See^2(T\cM)$ (an analog of metric) is said to be \textit{pseudo-Hermitian with respect to} the operator of an almost complex (or queer) structure $J$ if
\begin{equation}\label{Herm}
h(X, Y)=h(JX, JY) \text{~~for any vector fields~~} X, Y\in\fvect(\cM).
\end{equation}

\sssbegin[The classes of $J$ over $\Cee$, over $\Ree$ and at $p=2$; flatness when $a=0$]{Remarks}[The classes of $J$ over $\Cee$, over $\Ree$ and at $p=2$; flatness when $a=0$] Observe that over algebraically closed fields $\Kee$ of characteristic $p$, the conditions $J^2=k\id$, where $k\in\Kee^{\times}$, form a single equivalence class, for which we can select any $k\neq0$, for example, $k=1$ for any $p$, or $k=-1$ if~$p\neq 2$. So over $\Cee$, there is only one class of each of the two tensors $J$ (either even or odd) and one can select for the normal shape of $J$ the one satisfying $J^2=k\id$ for any $k\in\Cee^\times$, whereas over $\Ree$, each of the classes over $\Cee$ splits into two: there is a~  $J$ such that $J^2=\id$ and there is a~ $J$ such that~$J^2=-\id$.

If $p=2$, then if $\dim V$  is even, there are two classes of non-degenerate symmetric bilinear forms on $V$, see \cite{Leb}, where the orthogonal  Lie algebras preserving these forms are described and classified over algebraically closed fields. Note that  comparison with two types of orthogonal Chevalley groups is not applicable:  the correspondences ``quadratic forms" $\longleftrightarrow$ ``bilinear forms" have kernels and are not onto.

For the classification of non-degenerate symmetric pseudo-hermitian bilinear forms (which can be even or odd) or, rather, Lie superalgebras over $\Ree$ preserving these 2-forms, see \cite{Sr}.

Observe that various supermanifolds of high interest in theoretical physics (e.g., Minkowski superspaces, certain superstrings) are endowed with an almost complex structure given only on the fibers of a~ non-integrable distribution, not on the whole tangent space. For the integrability conditions (obstructions to ``flatness", i.e., possibility to drop the adjective ``almost") of almost complex, almost queer, and almost symplectic (or periplectic) structures defined on the fibers of a~ non-integrable distribution on a~given (super)manifold $\cM^{a|b}$, in particular, for a~ conceptual proof of Darboux's theorem (description of  $d\omega$ as the only, regardless of parity of $\omega$, obstruction to flatness of ortho-symplectic and periplectic structures), see \cite{BGLS}.   
 For the list of non-integrable distributions with infinite-dimensional Lie superalgebra of symmetries, radically different from its Cartan-Engel non-super analog, see \cite{KLS}.
 
 Note that if $a=0$, then any almost complex (queer) and any almost symplectic structure are flat (this observation seems to be new).  
\end{Remarks}

\ssec{Additional structures and invariant differential operators}\label{add} Let $\cM^{a|b}$ be a~(super)manifold with a~$\text{G}$-structure defined for Lie supergroups $\text{G}$ literally as for Lie groups. Here, $\text{G}=\GL(a|b)$ in the generic case, $\text{G}=\OSp(a|b)$ (resp., $\Pe(n)$) if the non-degenerate even, i.e., symplectic (resp., odd, i.e., periplectic) form is preserved in the tangent space $T_\star\cM$; for a~definition of these Lie supergroups, see \cite{LSoS}. We will only need the Lie superalgebras of these Lie supergroups, see~ \cite{Lo, LSh3}.

Let $T(V)$ and $T(W)$ be (super)spaces of tensor fields on $\cM$, whose fibers $V$ and $W$ are irreducible under the (super)group $\text{G}$ of linear diffeomorphisms; let
\be\label{D}
D: T(V)\tto T(W)
\ee
be an operator invariant under arbitrary changes of coordinates, i.e., diffeomorphisms of~$\cM$. The only \textit{differential} invariant operators $D$ on any (super)manifold are the exterior differentials
$d_i:\Omega^i\tto \Omega^{i+1}$ and their duals (denoted by abuse of notation just by $d$), see \cite{R1, BL, GLS}. For the description of invariant differential operators on (super)manifolds with a~symplectic, periplectic and (peri)contact structure, see a review \cite{LeLh}.

We will consider real (super)manifolds $\cM$ endowed with one or several of the  additional structures, see objects  in Subsec.~\ref{ss3}. The invariant differential operators \eqref{D} preserving these structures form Lie superalgebras  acting on $\widehat \Omega(\cM)$ and $\widehat \Sigma(\cM)$.

\ssec{The two types of (almost) periplectic supermanifolds}\label{ss2types} 
We consider the \textit{restricted} dualization over~ $\cF$, i.e., we set $(\Cee[\widehat x])^\#:=\Cee[\widecheck x]$, not $\Cee[[\widecheck x]]$, where $(\widehat x)^\#=\widecheck x$; then, the space of integral forms of degree $-\ell$, see object of type (iii) in Subsec.~\ref{ss3}, is
\[
\Sigma_{-\ell}:=(\Omega^\ell)^\#\otimes_{\cF}\Vol. 
\]

The operator dual to $d_i:\Omega^i\tto \Omega^{i+1}$ is $d_i^\#:\Sigma_{-i-1}\tto \Sigma_{-i}$, all these operators are briefly denoted by $d=\sum \widehat x_i\frac{\partial}{\partial x_i}$ with a natural action of $\widehat x_i$ on $\Sigma_{\bcdot}$, see \cite{BL}; the action of $d$ naturally extends from the spaces of forms to the spaces of pseudoforms.

\underline{General type}: On \textbf{general} (almost) periplectic supermanifolds, both $\omega:=\omega_\od$, the bivector $B$ defining the Buttin bracket are odd and the tangent space can only be identified with the \textit{shifted} cotangent space $\Pi(T^*_\star\cM)$, not with~$T^*_\star\cM$. Therefore, $\widehat \cM$ is not identified with $\widecheck \cM$, the equality $\widecheck x_i=\partial_{\widehat x_i}$ does not hold, and the bivector $B$ has to be given separately: the form $\omega\in\cF(\widehat\cM)$ and the bivector $B\in\cF(\widecheck\cM)$ are not related to each other.

\sssbegin[Why the volume element, and not the invertibility of $\omega_\od$, is the issue]{Remark}[Why the volume element, and not the invertibility of $\omega_\od$, is the issue]\label{whyspecial} The Gram matrix of a~non-degenerate $\omega_\od$ \textit{is} invertible, and its inverse is again of periplectic shape, so a~\textit{tensor} $\omega_\od^{-1}$ does exist pointwise --- it is the bivector of the anti-bracket, see \cite{BV1, BV2}. What $\omega_\od$ does \textbf{not} determine is the \textit{volume element}, and this is what the passage to $\cF(\widecheck \cM)$ needs, since $\Sigma_{-\ell}=(\Omega^\ell)^\#\otimes_\cF\Vol$. Indeed, let $X=\diag(A, D)$ be an even element of $\fpe(n)=\Lie(\Pe(n))$, i.e., preserving $\omega_\od$ with Gram matrix $\begin{pmatrix} 0&N\\ N^{T}&0\end{pmatrix}$; then $D=-N^{-1}A^{T}N$ and hence
\[
\str X=\tr A-\tr D=2\tr A\not\equiv 0,
\]
i.e., $\Ber$ is a~non-trivial character of $\Pe(n)$, trivial exactly on~$\SPe(n)$. For an \textbf{even} non-degenerate $\omega_\ev$ the same computation gives $\fosp(a|b)\subset\fsl(a|b)$, so $\str\equiv 0$ and the volume element (a~power of $\omega_\ev$) is preserved automatically --- whence identification \eqref{ident}. This is the super-analog of the well-known fact of periplectic geometry that, unlike its even counterpart, it carries no canonical Berezin volume, and the odd Laplacian is defined only after such a~volume has been chosen; and this is exactly why the \textit{special} type has to be singled out.

On the level of vector fields, the same fact reads as follows. For $\omega_\od=\sum_{1\leq i\leq n} dq_i\, d\xi_i$, we have
\be\label{DivLe}
\Div Le_f=2(-1)^{p(f)}\cS f, \text{~~where~~} \cS=\sum_{1\leq i\leq n}\frac{\partial^2}{\partial \xi_i\partial q_i}\text{~~is the odd Laplacian, see eq.~\eqref{DelS},}
\ee
so $\fsle(n)=\{Le_f\mid \cS f=0\}=\{Le_f\mid \Div Le_f=0\}$ is a~\textbf{proper} subalgebra of $\fle(n)$, whereas for an even $\omega_\ev$ the same computation gives $\Div H_f=0$ identically.
\end{Remark}

\underline{Special type}: An (almost) periplectic supermanifold with an (almost) periplectic form $\omega$ will be called \textit{special} if the volume element $\vvol$ is also preserved --- preserved, that is, by the vector fields that preserve $\omega$: the structure Lie superalgebra $\fle(n)$ is thereby replaced with~$\fsle(n)$. 
Then, the dual to~ $\omega$ bivector~ $B$ is well defined and using relation \eqref{JhB}, we arrive at a~ new notion of an (almost) \textit{K\"ahler supermanifold}. Note that on any almost symplectic supermanifold the volume element is automatically preserved and there is a~\textit{canonical} identification of the tangent space with the cotangent space; hence the bivector $B$ can be defined by the inverse of the Gram matrix of the bilinear form $\omega_\ev$, so the tensor $J$  can be transplanted from the tangent space to the cotangent space. Let us compare:
\begin{equation}\label{ident}
\begin{minipage}[c]{14cm}
\phantom{.}\hskip 0.5cm If the tangent space $T_\star\cM$ at the point denoted by star $\star$, not $\ast$, is endowed with an \textbf{even} non-degenerate bilinear form $\omega_\ev$, then $T_\star\cM$ can be (and is) canonically identified with the cotangent space $T_\star^*\cM$ and  $\widehat \cM$ can be (and hereafter is) canonically identified with $\widecheck \cM$, \textbf{and therefore}, $\widecheck x_i=\partial_{\widehat x_i}$.\\ 
\phantom{.}\hskip 0.5cm If $T_\star\cM$ is endowed with an \textbf{odd} non-degenerate bilinear form
$\omega_\od$, then the bivector $B$  has to be given separately. On a \textit{special} almost periplectic supermanifold, $B=\omega_\od^{-1}$ and $\widecheck x_i=\partial_{\widehat x_i}$, as in the (almost)  symplectic case. 
\end{minipage}
\end{equation}

\ssec{A traditional definition of almost K\"ahler (super)manifolds}\label{longdef} 
Let~ $h$ be a~pseudo-Hermitian with respect to $J$ non-degenerate symmetric form on the tangent bundle over~$\cM$.  If $p(h)+p(J)=\ev$ (resp., $p(h)+p(J)=\od$), define an almost symplectic  (resp., periplectic) 2-form~ $\omega$ by the formula
\begin{equation}\label{eq2}
\omega(X, Y):=h(JX, Y) \ \ \text{ for any $X, Y\in\fvect(\cM)$; hence, $p(\omega)=p(h)+p(J)$.}
\end{equation}

The (super)manifold $\cM$ with the triple $(\omega, h, J)$ is said to be \textit{almost K\"ahler}; it is \textit{K\"ahler} if $\omega$ and $J$ are \textit{parallel}, i.e., preserved by a~torsion-free connection, which is then unique but need not be flat, see Remark~\ref{parallel}. On supermanifolds,
each constituent of the triple $(\omega, h, J)$ 
can be either even or odd, but since any two of the constituents of the triple determine the third one by means of eq.~\eqref{eq2}, it follows that there are only the following types of super dimensions of an (almost) K\"ahler
manifold $\cM$ (first given in \cite{BGLS}):
\begin{equation}
\label{4types}
\renewcommand{\arraystretch}{1.4}
\begin{array}{|c|c|c|}
\hline &p(J)=\ev&p(J)=\od\cr \hline p(h)=\ev&2n|2m&2n|2n\cr \hline
p(h)=\od&2n|2n&n|n\cr \hline
\end{array}
\end{equation}

This table hides a catch: when $\omega$ is even, the symplectic part of it can live only in an even dimension $2k$, so for $J$ and $h$ odd, the superdimension $n|n$ is of the form $2k|2k$.

\sssbegin[Sign-definiteness of $h$ is unnecessary, and unattainable on any~superspace]{Remark}[Sign-definiteness of $h$ is unnecessary, and unattainable on any~superspace] M.~Verbitsky informed me that \lq\lq the sign-definiteness of $h$ in the
traditional definition of the K\"ahler manifold is unnecessary
(the published classification results are only known, however, for
sign-definite forms $h$); on the other hand, the flatness 
of the almost symplectic structure is definitely needed.\rq\rq
Note that the non-degenerate sesquilinear form on the superspace with non-zero even and odd components is never sign-definite: see Serganova's classification \cite{Sr} of real and quaternionic forms of simple complex Lie superalgebras.
 \end{Remark}

\sssbegin[Parallel, not flat]{Remark}[Parallel, not flat]\label{parallel} Being K\"ahler requires \textit{parallelism}: the tensors $\omega$ and $J$ are preserved by a~torsion-free connection. To require \textit{flatness} of $\omega$ and $J$ instead of  parallelism is unreasonable: they are preserved by a~torsion-free connection which is not \textit{flat}  since the~K\"ahler manifold may have a~ curvature; the curvature tensor and its covariant derivatives are local invariants of the triple $(\omega, h, J)$, and therefore no change of coordinates brings the triple to a~standard shape.

Where the word \textit{flat} came from is clear, and it is harmless as long as only \textit{one} ingredient is considered at a~time: a~torsion-free connection preserving $\omega$ alone exists if and only if $d\omega=0$, which is the hypothesis of Darboux's theorem, i.e., local flatness; a~torsion-free connection preserving $J$ alone exists if and only if the Nijenhuis tensor vanishes, and this again yields a~local normal shape. It is for the \textbf{triple} that the two readings part company, because $h$ --- unlike $\omega$ and $J$ --- is \textit{not} locally standard, and here is what this means$:$ a~non-degenerate symmetric $h$ can be brought to a~constant Gram matrix at a~\textbf{point}, but not on a~neighborhood of it. The obstruction is the curvature of the Levi-Civita connection of~$h$, a~local invariant of $h$ alone which does not vanish in general; this is why for~metrics there is no analog of Darboux's theorem or of the Newlander--Nirenberg theorem (see \cite{BGLS}), and why the triple $(\omega, h, J)$ has local invariants even though each of $\omega$ and $J$ has none. Accordingly, hereafter \lq\lq flat\rq\rq\ is used only in the sense \lq\lq locally standard, no local invariants\rq\rq\ and only of a~single structure at a~time, whereas the K\"ahler condition on the triple is \textbf{parallelism}, not flatness.\end{Remark}

\sssbegin[The connection is unique thanks to $h$, not to $\omega$]{Statement}[The connection is unique thanks to $h$, not to $\omega$]\label{prolong} Let $V$ be a~superspace, let $\fg\subset\fgl(V)$, and let $\fg^{\langle 1\rangle}$ be the first Cartan prolong of~$\fg$ $($the angle brackets tell it from the derived algebra $\fg^{(1)}$ of eq.~\eqref{derser}$)$. Two torsion-free connections compatible with a~$\fg$-structure differ by a~section of the bundle whose fiber is $\fg^{\langle 1\rangle}$; hence, such a~connection is unique if and only if~$\fg^{\langle 1\rangle}=0$. In our three cases$:$

$1)$ If $\fg=\faut(h)$ for an even non-degenerate symmetric $h$, i.e., if $\fg=\fosp(m|n)$ acting on $V=\Cee^{m|n}$, then $\fg^{\langle 1\rangle}=0$; 

$2)$ If $\fg=\faut(h)$ for an \textbf{odd} non-degenerate symmetric $h$ on $V=\Cee^{n|n}$, 
then $\fg^{\langle 1\rangle}=0$; 

$3)$ If $\fg=\faut(\omega_\od)=\fpe(n)$ acting on $V=\Cee^{n|n}$, 
then $\fg^{\langle 1\rangle}\neq0$; it is the degree-$1$ component of $\fle(n)$ in its standard grading. 

Therefore, the torsion-free connection preserving the triple $(\omega, h, J)$ is unique for \textit{either} parity of~$h$, exactly as in the classical case, whereas $\omega$ alone never determines one.
\end{Statement}

The classical answers calibrate this$:$ $\fo(m)^{\langle 1\rangle}=0$, the Levi-Civita connection being unique, whereas $\fsp(2k)^{\langle 1\rangle}=S^3(V^*)\neq0$, a~symplectic connection being far from unique. Item $3)$ reproduces the Cartan prolong $(\Cee^{n|n}, \ \fpe(n))_*=\fle(n)$ of Subsec.~\ref{2.6}; and note that the two Lie superalgebras of items $2)$ and $3)$, of one and the same $\sdim=(n^2|n^2)$ and isomorphic to each other after the parity shift $\Pi$ of $V$, are told apart by their first prolongs --- a~prolong depends on the pair: the algebra and the module, not on the Lie superalgebra alone. These computations are contained in \cite{P2}$:$ item~ $1)$ is its Ch.~3, Theorem~3.3, and item~ $2)$ its Ch.~1, Theorem~1.2, in which $\fpe(n)$ denotes the Lie superalgebra preserving an odd \textit{symmetric} form. Item~ $3)$ is the companion of item $2)$ under the parity shift $\Pi$ of $V$, and \cite{P2} stresses that the prolongs of $(\fg_{-1}, \fg_0)$ and of $(\Pi\fg_{-1}, \fg_0)$ are essentially different.%

 It is not flatness but parallelism that carries the computations of this article --- all of them performed in the flat model of Subsec.~\ref{flat} --- over to a~\textit{curved} $\cM$: the invariant operators involved are of the first order in the structure tensors, and on every~K\"ahler $\cM$ the covariant derivatives of $\omega$, $h$ and $J$ vanish, so the same relations hold at every point and the same Lie superalgebra is generated. The curvature enters only in the second order terms, and these land in the Laplacian, which is central in the answer anyway, see eq.~\eqref{ADolb}. What destroys the relations is the failure of parallelism, i.e., the adjective \lq\lq almost\rq\rq, not the presence of curvature.

\sssec{$B$ vs. $\omega$ in the definition of almost K\"ahler (super)manifolds}\label{omegaB} Let $h^*$ be the form on the  cotangent bundle 
dual to the form $h$ on the tangent bundle. If the tangent space can be identified  with the cotangent space (perhaps, up to the inversion of the parity), then the action of the operator $J$ 
can be transferred from the tangent space to the cotangent one, and the bivector field $B$ can be defined by the expression (cf. with eq.~\eqref{eq2})
\begin{equation}\label{JhB}
\renewcommand{\arraystretch}{1.4}
\begin{array}{l}
B(df, dg):=h^*(Jdf, dg) \; \text{ for any $f,g\in\cF(\cM)$},\\
\end{array}
\end{equation}

Since the bracket in the space of functions $\cF(\cM)$ (resp., $\Pi(\cF(\cM))$) is determined by an~ even (resp., odd) bivector field $B$, it is natural to replace the condition of \textit{closedness} $d\omega=0$ ensuring the Jacobi identity with the equivalent condition ~\eqref{bb}. This answers one more question of Verbitsky's --- what, for an odd $\omega$, the analog of closedness is. The answer does not depend on the parity of $\omega$: it is condition \eqref{bb} on the bivector inverse to $\omega$, whereas $d\omega$ is, regardless of the parity of $\omega$, the only obstruction to the normal shape of $\omega$, see \cite{BGLS}. For an even $\omega$, this obstruction is not one condition but two: the space of $3$-forms is the sum of the primitive forms and of the multiples of $\omega$, so each of the two components of $d\omega$ can be required to vanish separately, as is done in the classification of almost Hermitian manifolds, see \cite{GH}. On supermanifolds, both halves survive --- but not on all of them, and for an odd $\omega$ they survive always: see Lemma~\ref{halves} and Subsec.~\ref{ssTwoCond}.  This bracket admits a~deformation with parameter running over a~
(super)variety  $\cQ$, see Section~\ref{S2}. 


\textit{For $B$ even}, this deformation is the well-known quantization of
the Poisson bracket. Note  that the ``well-known" 1-parameter family of deformations describes only the case of polynomial or analytic functions whereas for other types of functions the variety of parameters  might be, and actually is, different, e.g., for the case of functions with compact support, so natural in physical models, see ~\cite{KT}. 

\textit{For $B$ odd}, the Buttin bracket given by $B$ admits deformations
parameterized by a~singular supervariety $\cQ$, 
see  Section~\ref{S2}. If $\sdim \cM=2|2$, this description of deformations can be applied to $B$ even. Thus,  \textbf{the definition of almost K\"ahler and almost hyper-K\"ahler supermanifolds based on the bivector $B$ instead of $\omega$ depends on a~parameter running through~$\cQ$}. However, the Lie superalgebra of invariant differential operators acting in the space of (pseudo-)forms (an analog of the known $\fsl(2)$-action in $\Omega(M)$) is ``of interest" only at one or at most $3$ of the total of at most $6$ exceptional points of $\cQ$ described in Theorem~$\ref{thdefb}$.

\ssec{The flat model in which the Lie superalgebras of invariant operators are computed}\label{flat} Several of the Lie superalgebras described in this article are obtained by an explicit computation which I now describe once and for all; below I refer to it as the \textit{flat model}.

On $\Cee^{a|b}$ with coordinates $z=(z_1, \dots , z_{a+b})$, set $\widehat y_c:=dz_c$, so $p(\widehat y_c)=p(z_c)+\od$. The invariant operators with constant coefficients acting on $\widehat \Omega$ are the polynomials in the operators
\be\label{engine}
\widehat y_c\ \ \text{(left multiplication)},\ \ \partial_{\widehat y_c}\ \ \text{(convolution)},\ \ t_c:=\partial_{z_c},
\ee
to which, in Subsec.~\ref{fgMyst} only, the multiplications $x_c:=z_c\cdot$ and the odd parameter $\hbar$, supercommuting with all the other operators and hence such that $\hbar^2=0$, are added. The operators \eqref{engine} span a~super Weyl algebra: the only non-zero super-commutators of the generators are
\be\label{engineRel}
[\partial_{\widehat y_c}, \widehat y_e]=\delta_{ce},\ \ [t_c, x_e]=\delta_{ce}.
\ee
Hence, every product of the generators \eqref{engine} can be normally ordered with exact rational coefficients and the Sign Rule. Given a~collection of operators, the phrase \textit{the closure of the algebra generated by} them means, here and everywhere below, the smallest Lie superalgebra containing the collection; it is computed by taking the linear span of the collection, bracketing its basis elements pairwise, adjoining whatever is linearly independent of what is already there, and repeating until nothing new appears --- a~finite computation in linear algebra over~$\Qee$. The distinction is essential: closing the collection of operators \textit{itself} under the bracket returns the iterated brackets one by one, a~set and not a~space, whereas it is the linear combinations adjoined at every step that make the result an algebra and give it a~superdimension. All superdimensions quoted below were obtained in this way and are stated for the flat model.

\ssec{Lie superalgebras that gate-crash the Differential Geometry on manifolds. On supermanifolds, more complicated superalgebras enter the game}\label{ssgate}{}~{} 

\underline{Arbitrary manifold}. Obviously, the space $\Omega$ of differential forms on a  manifold is a~supercommutative superalgebra with respect to the exterior product $\wedge$ the sign of which is skipped in what follows. The space $\widehat \Omega$ on any supermanifold is also a~supercommutative superalgebra. The restricted dual to $\Omega$ is also $\Omega$ on manifolds and $\Sigma$ on supermanifolds. The space $\Sigma$ is an $ \Omega$-module whereas $\widehat \Sigma$ is an~$\widehat \Omega$-module. The exterior differential is an odd derivation of $\widehat \Omega$, not only of $\Omega$;  it is a~vector field on~$\widehat \cM$. The dualization defines the action of $d$ on $\widehat \Sigma$. The exterior differential and the Berezin integral in the purely odd case are the only differential operators invariant under any coordinate changes, see \cite{BL}. 

\underline{Almost symplectic manifold $(\cM, \omega_\ev)$}. On the space $\Omega$ of complexified differential forms on an almost symplectic manifold (or supermanifold) $\cM^{2a|b}$ with a~non-degenerate 2-form $\omega_\ev$, there acts a~ Lie superalgebra $\fg$ described as follows.   
\be\label{Xpm}
\begin{array}{l}
\text{Let $X^+$ be the operator of left multiplication by~ $\omega_\ev$;}\\
\text{let $X^-$ be the convolution with the bivector $B:=(\omega_\ev)^{-1}$.}
\end{array}
\ee 
The operators $X^+$ and $X^-$ generate $\fsl(2)\simeq\fsp(2)$, while the exterior differential $d$ and the codifferential $\delta$ dual to it with respect to $\omega_\ev$, see eq.~\eqref{H} and the sentence after it, generate (unless $a=0$, see eq.~\eqref{a=0}) the Heisenberg Lie superalgebra~$\fhei(0|2)$.  Since we consider complexified differential $2$-form,  and since over $\Cee$, there is only one class of non-degenerate symmetric bilinear even forms ($\eps_i=1$ for all $i$, see eq.~\eqref{a1}), it suffices to consider one of the following two equivalent versions (non-split for $k=0$ or split for $l=0$), see eq.\eqref{clOmega}. Over $\Ree$, this class splits into several classes  depending on the signature $\sum\eps_i$, where $\eps_i=\pm 1$, see eq.~\eqref{a1}:
\be\label{clOmega}
\begin{array}{l}
\text{$\omega_\ev:=\sum \widehat p_i\widehat q_i+\nfrac12 \begin{cases}\sum (\widehat \theta_j)^2&\text{for $k=0$},\\
\sum_{j\leq b/2} \widehat \theta_j\widehat \theta_{b-j}&\text{for $l=0$,}
\end{cases}$}\\ 
\text{$B:=\sum \partial_{\widehat p_i}\partial_{\widehat q_i}+\nfrac12 \begin{cases}\sum (\del_{\widehat \theta_j})^2&\text{for $k=0$},\\
\sum_{j\leq b/2} \del_{\widehat \theta_j}\del_{\widehat \theta_{b-j}}&\text{for $l=0$.}
\end{cases}$}
\end{array}
\ee
 Both $\omega_\ev$ and $B$ are normalized so as to get $\deg:=\sum\, \hat x_i \del_{\hat x_i}$ for $x=(p,q,\theta)$ in the bracket $[X^-, X^+]$: 
\be\label{Hsl}
H_{sl}:=[X^-, X^+]=\deg- a+\nfrac12 b,
\ee
 whereas 
\be\label{H}
H:=[d, \delta]=\sum \nfrac{\partial^2}{\partial q_i\partial p_i}+\begin{cases}0&\text{for $k=0$,}\\
\sum_{j\leq b/2}  \nfrac{\partial^2}{\partial \theta_j\partial \theta_{b-j}}&\text{for $l=0$}.
\end{cases}
\ee
Here, $\delta$ is the codifferential dual to $d$ with respect to $\omega_\ev$ --- on the even indeterminates $\delta=\sum \partial_{\widehat p_i}\partial_{q_i}$, and on the odd ones it is displayed just below --- and \textbf{not} the operator $d^*$ of eq.~\eqref{Del}$;$ see the paragraph after eq.~\eqref{Svanishes}. This ``symplectic Laplacian" $H$ spans the center of $\fhei(0|2)$ generated by $d$ and $\delta$, except for the case of $l$ odd since $d:=\sum \widehat \theta_i\del_{\theta_i}$ and $\delta:=\begin{cases}\sum \del_{\widehat \theta_i}\del_{\theta_i}&\text{for $l=0$},\\
\sum_{j\leq b/2} \del_{\widehat \theta_i}\del_{\theta_{b-i}}&\text{for $l\neq 0$},
\end{cases}
$ and hence 
\be\label{a=0}
\fg\simeq\begin{cases}\fhei(0|2)\inplus\fsl(2)&\text{for $l$ even},\\
\Pi(\Cee^2)\inplus\fsl(2)&\text{for $l$ odd}.
\end{cases}
\ee
This $\fg$-action obviously extends to $\widehat \Omega$. 

\sssec{General almost periplectic supermanifold $(\cM, \omega_\od)$}

Here, the operators $d$, the multiplication by $\omega_\od$, and their composition $\omega_\od d$ are the only differential operators on $\widehat \Omega^{\bcdot}(\cM)$ and $\widehat \Sigma_{\bcdot}(\cM)$ invariant under the diffeomorphisms of $\cM$ preserving $\omega_\od$; they generate $\fhei(0|2)\simeq\fsl(1|1)$,   
see \cite{LeLh}.

\sssec{Special almost periplectic supermanifold $(\cM, \omega_\od)$ with preserved volume element $\vvol$} The on\-ly \textit{differential} operators $D$ invariant under arbitrary changes of coordinates of $\cM$ preserving \textit{both}~ $\omega_\od$ and $\vvol$ are $d$, its dual $d^*$, the operators $X^+$ of left multiplication by $\omega_\od$, the convolution~ $X^-$ with the bivector $B$ dual to $\omega_\od$, and products thereof, see \cite{LeLh}. 

Consider the pericontact Lie superalgebra $\fm(k)$ and the contact Lie superalgebra $\fk(2n+1|m)$ with their \textit{standard} gradings, see eq.~\eqref{standgr}.
Then, the Lie superalgebra 
\be\label{h}
\fh:= \fb\fa(1)\inplus\fsl(1|1),
\ee
 generated by the operators $X^+$, $X^-$, $d$, and $d^*$ resembles the Lie superalgebra $\fg$ in the case of symplectic supermanifolds, see eq.~\eqref{a=0}. (In fact, the closure of the algebra generated by these operators is slightly larger: it is $\fb\fa(2)\inplus\fsl(1|1)$, 
 and $\fh$ is its subalgebra; for a~computation, see Lemma~\ref{gl11}.) What the elements of $\fh$ are needed \textit{for} --- the description of the irreducible continuous $\fle(n)$- and $\fsle(n)$-modules --- is expounded in \cite{LeLh}. 

 
\sssec{Both $J$ and $h$ even, see \cite{P1, V1}}\label{p1} In this case, the passage from manifolds to supermanifoldsis direct and immediate. Thanks to identification \eqref{ident}, the complex structure on $T_\star\cM$ at point $\star$ can be considered on $T_\star^*\cM$ as well. Several classical operators on the \textit{Dolbeault algebra} $A^{\bcdot, \bcdot}(M)$ of complexified differential forms on $M$ are well-known: the exterior differential~ $d$ together with 
\be\label{d_c}
d_c := JdJ^{-1}, 
\ee
see \cite[p.19]{V1}, the duals $d^*$ and $d_c^*$ of $d$ and $d_c$ with respect to $h$, and the analog of the $H$ of eq.~\eqref{H} that these duals produce 
\be\label{Del}
\Delta:=[d,d^*] = [d_c,d_c^*]. 
\ee

Thus, on the Dolbeault algebra of any (almost) K\"ahler (super)manifold, the semi-direct sum 
$\fhei(0|4)\inplus\fsp(2)$ acts, where $\fsp(2)$ is generated by the operator of left multiplication by $\omega$ and the convolution with the dual bivector $B$, while 
$\fhei(0|4)_\ev  :=\Span(\Delta)$ 
and  $\fhei(0|4)_\od = \Span( d, d_c^*)\oplus \Span(d^*, d_c)$ is the direct sum of two subspaces isotropic with respect to the non-degenerate symmetric bilinear even form given by $(d,d^*) = (d_c,d_c^*) = 1$, the other inner products of basis elements being equal to~0. On supermanifolds, the action of the Lie superalgebra 
\be\label{ADolb}
\fh:=\fhei(0|4)\inplus\fsp(2)=\fhei(0|4)\inplus\fsl(2)
\ee
naturally extends to the space of pseudo-forms; note that thanks to identification \eqref{ident}, the space of (pseudo-)integral forms is isomorphic to the space of (pseudo-)differential forms. We write
$\fsp(2)$  to highlight the resemblance with the case of twisted Dolbeault algebra in Subsec.~\ref{ssTwDol} or $\fsl(2)$ to highlight the resemblance with eq.~\eqref{ADolbOdd}. In the flat model of Subsec.~\ref{flat}, the closure of the algebra generated by $\{E, F, d, d_c\}$ is indeed of $\sdim=(4|4)$, see Remark~\ref{validation}; this is the only one of the four cases of table \eqref{4types} in which the classical answer survives superization unchanged.

\sssec{Both $J$ and $h$ odd}\label{p2} This is another, apparently new, type of almost K\"ahler supermanifolds~$\cM$. Here, however, the constructions of Subsec.~\ref{p1} can \textbf{not} be repeated verbatim: this entry of table \eqref{4types} exists only in the weakened, $B$-paradigm sense of Subsec.~\ref{omegaB}, and the answer is \textit{not} the one of Subsec.~\ref{p1}, as Lemma~\ref{p2a} and Statement~\ref{p2b} show. The action of the Lie superalgebra described in Statement~\ref{p2b} naturally extends from the spaces of forms to the spaces of pseudo-forms $\widehat A^{\bcdot, \bcdot}(\cM)$ and $\widehat B^{\bcdot, \bcdot}(\cM)$. Note that, clearly, the space $B^{\bcdot, \bcdot}(\cM)$ is an $A^{\bcdot, \bcdot}(\cM)$-module and $\widehat B^{\bcdot, \bcdot}(\cM)$ is an $\widehat A^{\bcdot, \bcdot}(\cM)$-module.

\sssbegin[An odd $h$ and an odd $J$ give no even $\omega$]{Lemma}[An odd $h$ and an odd $J$ give no even $\omega$]\label{p2a} There is no triple $(\omega, h, J)$ with both $h$ and $J$ odd, $J^2$ a~non-zero scalar, and $\omega=h(J\ -, -)$ a~non-degenerate \textbf{even} $2$-form, in any dimension~$(n|n)$. Let $h$ be symmetric and $\omega$ anti-symmetric, so let $\beta:=h(J\ -, -)$ denote the tensor before it is known which of the two symmetries it is $($if either$)$: for $J^2\neq0$, the tensor $\beta$ is a~non-degenerate even form of \textbf{neither} symmetry. Indeed, the two diagonal blocks of $\beta$ are always of one and the same symmetry, whereas either type of non-degenerate even form requires them to be of opposite symmetries.
\end{Lemma}

\begin{proof} Everything is pointwise, so let $V=V_\ev\oplus V_\od$ be the tangent space at a~point, of $\sdim V=(n|n)$.

First, let us say what the matrices are. The form $h$ being \textit{odd} means that
the whole of $h$ is one bilinear map $V_\ev\times V_\od\tto\Kee$:
in bases of $V_\ev$ and of $V_\od$ it is a~single $n\times n$ matrix $N$, and $h$ is non-degenerate if and only if $N$ is invertible, so that
\be\label{Hmixed}
h=\begin{pmatrix} 0&N\\ \varepsilon N^{T}&0\end{pmatrix},\quad\text{where $\varepsilon=1$ for $h$ symmetric and $\varepsilon=-1$ for $h$ anti-symmetric.}
\ee
Here, $\varepsilon=1$ is the case at hand: the pseudo-Hermitian $h$ of a~K\"ahler structure is symmetric, always. The other sign is carried along only to widen the Lemma for nothing: an \textit{odd} non-degenerate form, unlike an even one, can be of either symmetry, and the anti-symmetric one is the periplectic form of Subsec.~\ref{Sperip}, so nothing is lost by covering it here as well. The sign plays no role in what follows. Likewise, the operator $J$ being odd means that it interchanges $V_\ev$ and $V_\od$, i.e., that it is anti-diagonal,
\be\label{Jodd}
J=\begin{pmatrix} 0&Q\\ P&0\end{pmatrix},\quad\text{where } P\colon V_\ev\tto V_\od,\quad Q\colon V_\od\tto V_\ev.
\ee

Next, the matrix of $\beta=h(J\ -, -)$ is $J^{T}h$, and eqs.~\eqref{Hmixed} and~\eqref{Jodd} give
\be\label{OmDiag}
\beta=J^{T}h=\begin{pmatrix} 0&P^{T}\\ Q^{T}&0\end{pmatrix}\begin{pmatrix} 0&N\\ \varepsilon N^{T}&0\end{pmatrix}=\begin{pmatrix} \varepsilon\,(NP)^{T}&0\\ 0&Q^{T}N\end{pmatrix}.
\ee
So $\beta$ is automatically block-\textit{diagonal} --- as it must be, being even --- and the only conditions to impose are the symmetries of its two diagonal blocks.

Now, an even $2$-form is anti-symmetric on $V_\ev$ and symmetric on $V_\od$; an even symmetric form is symmetric on $V_\ev$ and anti-symmetric on~$V_\od$. Either way, \textbf{the two blocks of eq.~\eqref{OmDiag} must be of opposite symmetries}.

They never are. Indeed, $J^2=\varkappa\,\id$ reads $QP=PQ=\varkappa\,\id$, so for $\varkappa\neq0$ both $P$ and $Q$ are invertible and $Q=\varkappa P^{-1}$. Set~$A:=NP$. Then, the two blocks of eq.~\eqref{OmDiag} are
\be\label{congr}
\beta\vert_{V_\ev}=\varepsilon\,A^{T}\ \text{ and }\ \beta\vert_{V_\od}=Q^{T}N=\varkappa\,(P^{-1})^{T}A\,P^{-1}.
\ee
A~matrix and its transpose have the same symmetry, and a~congruence $A\mapsto C^{T}AC$ preserves symmetry as well; so both blocks of eq.~\eqref{congr} have the symmetry of $A$, whatever $A$ is. Being asked to be of opposite symmetries, each of them is then both symmetric and anti-symmetric, hence zero --- and $\beta$ is degenerate, contrary to the assumption. Therefore, $\varkappa=0$.\end{proof}

Thus, $\beta$ is even, and its $2$-form avatar is supported on one parity block only. In the $B$-paradigm, let $E$ be the operator of left multiplication by this avatar and let $F$ be the convolution with the dual bivector.

\sssbegin[What $E$, $F$ and $d$ generate when both $J$ and $h$ are odd]{Statement}[What $E$, $F$ and $d$ generate when both $J$ and $h$ are odd]\label{p2b} In the flat model of Subsec.~$\ref{flat}$ on $\Cee^{n|n}$ $($the answer does not depend on $n$; verified for $n=1,2,3)$, the closure of the algebra generated by $\{E, F\}$ is $\fsl(2)$, whereas the closure of the algebra generated by $\{E, F, d\}$ is 
\[
\fg:=\Pi(\Cee^2\oplus\Cee)\inplus\fsl(2).
\]
\end{Statement}

In particular, the answer here is \textbf{not} $\fhei(0|4)\inplus\fsp(2)$ of Subsec.~\ref{p1}: the four differentials collapse to three (the analog of $\delta_c$ is proportional to $\delta$) and the Laplacian $\Delta=[d,\delta]$ vanishes, because the base derivatives that occur in it are odd and hence square to zero; the Heisenberg Lie superalgebra $\fhei(0|4)$ contracts to the abelian $(0|3)$-dimensional space. Such a~phenomenon --- odd Laplacians surviving while even ones die, or the other way round --- recurs in situations described in Subsecs.~\ref{oddh} and~\ref{Sost}.

\paragraph{The tensors $J$ and $h$ have opposite parities, $\cM$ is special periplectic}\label{p3} These are two more (new) types of almost K\"ahler (super)manifolds~$\cM$. Since $p(\omega)=p(h)+p(J)$, opposite parities of $h$ and $J$ mean that $\omega$ is odd, i.e., that $\cM$ is periplectic; being in addition special, it is the $\cM$ of Subsec.~\ref{Sperip}, where these two cases are therefore treated: the Lie superalgebra generated by the operators $\cX^{\pm}$ and $d$ is the $\fb\fa(2)\inplus\fsl(1|1)$ of eq.~\eqref{ADolbOdd}, see Lemma~\ref{gl11}, whereas the Dolbeault operators generate more, see Statement~\ref{stmixed}. In both cases, the action naturally extends from the spaces of forms to the spaces of pseudo-forms. 

\ssec{From the $\fsp(2)$-action on the space of differential forms on a~symplectic (super)manifold to an $\fosp(1|2)$-action on the space of twisted forms}\label{JB}  In 1977, J.~Bernstein extended the known $\mathfrak{sp}(2)$-action on the space $\Omega$ of (complexified) differential forms on a~ symplectic manifold $M$, see \cite{Weil}, to the 
$\mathfrak{osp}(1|2)$-action on the space $\Omega_\hbar$ of differential forms taking values in the line bundle with a~connection corresponding to a~contact distribution, with the curvature of the connection equal to $\hbar$ times the symplectic form for any $\hbar\in\Cee^\times$. 

The original (unpublished) Bernstein's construction 
was directly extended to supermanifolds and one can read about this extended form in the easily available papers, see a review \cite{LeLh}. 
In local coordinates, the operators $D^+$ and $D^-:=[D^+, X^-]$ that generate $\mathfrak{osp}(1|2)$ have the form 
\begin{equation}\label{Dpm}
\begin{split}
 &D^+= \sum\widehat x_i\partial_{x_i}- \frac{\hbar}{2}\left(\sum
 (\widehat q_ip_i - \widehat p_iq_i)+\sum \widehat\theta_j \theta_j\right),\\
 &D^- =\sum\left(\partial_{\widehat p_i}\partial_{q_i} - \partial_{\widehat
 q_i}\partial_{p_i}\right ) +\sum \partial_{\widehat\theta_j} \partial_{\theta_j}- \frac{\hbar}{2}\left(\sum(p_i\partial_{\widehat
 p_i} +q_i \partial_{\widehat q_i})+\sum \theta_j\partial_{\widehat\theta_j }\right).\\
\end{split}
\end{equation}


Note that the largest algebra of invariant operators acting in $\widehat \Omega(\cM)$ is $\mathfrak{osp}(1|2)\oplus \Cee H$, see eq.~\eqref{H}.

\sssec{Twisted Dolbeault algebra}\label{ssTwDol} Introduce the \textit{twisted Dolbeault algebra} $A^{\bcdot, \bcdot}_{\hbar_1, \hbar_2}(M)$ of complexified differential forms on the (super)manifold $M$ endowed with a~line bundle $L$ with connection, and its dual space $B^{\bcdot, \bcdot}_{\hbar_1, \hbar_2}(M)$ in presence of~$\omega_\od$. Here is the twist, described at once. Bernstein's construction of Subsec.~\ref{JB} twists the one differential $d$ by a~connection whose curvature is proportional to $\omega$; here there are \textit{two} differentials, $d$ and $d_c$ of eq.~\eqref{d_c}, so $L$ carries a~pair of connections, one adapted to each of them, with the $1$-forms $\alpha$ and $\alpha_c$ of Theorem~\ref{thtwist} below. The two subscripts $\hbar_1$ and $\hbar_2$ are two complex numbers if $\omega$ is even, or two supercommuting Grassmann parameters if $\omega$ is odd, these $\hbar_1$ and $\hbar_2$ are the factors by which the two curvature forms are proportional to $\omega$, and the twisted differentials are those of eq.~\eqref{Dtwist}, namely, $D^+=d+\hbar_1\alpha$ and $D^+_c=d_c+\hbar_2\alpha_c$.

\sssbegin[The twist turns $\fhei(0|4)\inplus\fsl(2)$ into $\fosp(2|2)$]{Theorem}[The twist turns $\fhei(0|4)\inplus\fsl(2)$ into $\fosp(2|2)$]\label{thtwist} Let $\alpha$ be the connection $1$-form of Bernstein's construction, see Subsec.~$\ref{JB}$, normalized by $[d, \alpha]=-2X^+$, and let $\alpha_c$ be the $1$-form obtained from it by the rotation $J$ acting on the differentials, so that $[d_c, \alpha_c]=-2X^+$ as well, whereas $[d, \alpha_c]=[d_c, \alpha]=0$. Set
\be\label{Dtwist}
D^+:=d+\hbar_1\alpha,\ \ D_c^+:=d_c+\hbar_2\alpha_c,\ \ \text{where}\ \ \hbar_1 \hbar_2\neq0,
\ee
so that
\[
\nfrac12[D^+, D^+]=-2\hbar_1X^+\ \ \text{and}\ \ \nfrac12[D_c^+, D_c^+]=-2\hbar_2X^+,
\]
i.e., the two twisting constants are the factors by which the curvature forms are proportional to~$\omega$. Then, in the flat model of Subsec.~$\ref{flat}$, for one even $J$ and an even $h$, the closure of the algebra generated by $\{D^+, D_c^+, X^+, X^-\}$ is a~\textbf{simple} Lie superalgebra of $\sdim=(4|4)$, whose even part is $\fsl(2)\oplus\Cee z$ with $z$ central in it, and on whose odd part $\ad_z$ has two eigenvalues $\pm\lambda$, $\lambda\neq0$, each of multiplicity~$2$. Hence, it is $\fosp(2|2)\simeq\fsl(1|2)$. This $z$ is the twisted Laplacian $\Delta_{\hbar_1, \hbar_2}$ of eq.~\eqref{Deltw} below.
\end{Theorem}

\begin{proof} The closure is computed in the flat model of Subsec.~\ref{flat}; 
simplicity is verified by checking that the ideal generated by each basis element is the whole algebra, and $z$ is found by solving the linear system that says an element of $\fg_\ev$ commutes with all of~$\fg_\ev$. The identification of $z$ with the twisted Laplacian of eq.~\eqref{Deltw} is part of the same computation: the element $[D^+_c, D^-]$ is proportional to $z$ and lies outside $[\fg_\ev, \fg_\ev]=\fsl(2)$, which was checked over $\Cee$ on $\Cee^{2}$, $\Cee^{4}$ and $\Cee^{6}$ for six pairs of non-zero twisting constants, together with the degeneration $\hbar_1=\hbar_2=0$, which returns eq.~\eqref{Del} and the $\fhei(0|4)\inplus\fsl(2)$ of eq.~\eqref{ADolb}, and with the one-parameter degeneration, which returns Bernstein's~$\fosp(1|2)$. The identification then needs no further computation: by the classification of simple finite-dimensional Lie superalgebras, see a~review \cite{Lo}, the one with even part $\fsl(2)\oplus\Cee$ and $4$-dimensional odd part is $\fosp(2|2)$, in which the $\pm\lambda$-eigenspaces of the $\fo(2)$-generator $z$ are the two copies of the tautological $\fsp(2)$-module~$\Cee^{2}$.

The choice of $\alpha_c$ is essential and not the obvious one: the $1$-form obtained from $\alpha$ by rotating both the differentials and the coefficients has $[d_c, \alpha_c]=0$ --- a~\textit{flat} connection --- and with it the closure of the algebra generated by $\{D^+, D_c^+, X^+, X^-\}$ is of $\sdim=(4|4)$ but is \textbf{not} simple: it has an~ideal of $\sdim=(1|2)$, so it is not $\fosp(2|2)$.\end{proof}

The largest algebra of invariant operators acting in $A^{\bcdot, \bcdot}_{\hbar_1, \hbar_2}(M)$ is $\mathfrak{osp}(2|2)$, the Laplacian being \textit{inside} it, and not beside it, see eq.~\eqref{Deltw} and the word on the notation after it. The $\fosp(2|2)$ is Theorem~\ref{thtwist}: in the flat model, the operators $d$ and $d_c$, twisted by the two connections and taken together with $\cX^\pm$, do close to $\fosp(2|2)$, and to nothing larger. So the twist turns $\fhei(0|4)\inplus\fsl(2)$ of eq.~\eqref{ADolb} into $\fosp(2|2)$, one storey above the passage $\fhei(0|2)\inplus\fsl(2)\rightsquigarrow\fosp(1|2)$ of Subsec.~\ref{JB}.

The Laplacian to take here is the \textit{twisted} one: in terms of the twisted differentials $D^+$ and $D^+_c$ of eq.~\eqref{Dtwist} and of their partners $D^-:=[X^-, D^+]$ and $D^-_c:=[X^-, D^+_c]$, it is
\be\label{Deltw}
\Delta_{\hbar_1, \hbar_2}:=[D^+_c, D^-]=-[D^+, D^-_c],
\ee
the two expressions coinciding up to the sign of the K\"ahler identity, exactly as the two expressions of eq.~\eqref{Del} coincide; and for $\hbar_1=\hbar_2=0$, eq.~\eqref{Deltw} returns eq.~\eqref{Del}. This $\Delta_{\hbar_1, \hbar_2}$ is \textbf{not} a~direct summand of the answer$:$ the answer is $\fosp(2|2)$, and not $\fosp(2|2)\oplus\Cee\Delta_{\hbar_1, \hbar_2}$. Indeed, $\Delta_{\hbar_1, \hbar_2}$ is the generator $z$ of the $\fo(2)\subset\fosp(2|2)$ of Theorem~\ref{thtwist}, so that $\fosp(2|2)_\ev=\fsl(2)\oplus\Cee\Delta_{\hbar_1, \hbar_2}$. 
Unlike its untwisted namesake, it is central in the even part only: on the odd part, $\ad_{\Delta_{\hbar_1, \hbar_2}}$ has the two eigenvalues $\pm\lambda$ with $\lambda\neq0$, see Theorem~\ref{thtwist} --- which is precisely what makes it a~member of the answer rather than a~summand of it.

The $\Delta$ of eq.~\eqref{Del} itself does not survive the twist. It is central \textit{before} the twist --- it spans the center of the $\fhei(0|4)\inplus\fsl(2)$ of eq.~\eqref{ADolb} --- but after the twist it neither lies in the closure nor commutes with the twisted differentials, since $[\Delta, \alpha]\neq0$; so it is not an invariant operator of the twisted picture at all. In this, the two-parameter case differs from the one-parameter one, where $H$ of eq.~\eqref{H} \textbf{is} a~separate summand, see eq.~\eqref{Dpm}.

Here is the description of $B^{\bcdot, \bcdot}_{\hbar_1, \hbar_2}(M)$. First, \textit{which} form the two curvature forms are proportional to: there is only one form for them to be proportional to, and it is the same one for both connections. In the K\"ahler case of Theorem~\ref{thtwist} it is the even $\omega$, the two connections differing by the rotation $J$ and not by the form, which is what the normalization
\be\label{curvtwo}
[d, \alpha]=[d_c, \alpha_c]=-2X^+,\ \ \text{whereas}\ \ [d, \alpha_c]=[d_c, \alpha]=0,
\ee
of Theorem~\ref{thtwist} says. In presence of $\omega_\od$, i.e., on a~special almost periplectic $\cM$, the one form available is the odd $\omega_\od$, and it serves as the curvature form of both connections, see eq.~\eqref{alphatwo}; that is why the two connection forms of eq.~\eqref{alphatwo} differ only in the parameter. 

Now, whatever the pair $(L, \nabla)$ of connections on the line bundle, with the curvature forms proportional with the coefficients $\hbar_1, \hbar_2$ to the form just named, is that produces the twist, i.e., $A^{\bcdot, \bcdot}_{\hbar_1, \hbar_2}(M)=A^{\bcdot, \bcdot}(M)\otimes L$, the space of \textit{integral} forms is twisted by the \textbf{same} bundle:
\be\label{Btwist}
B^{\bcdot, \bcdot}_{\hbar_1, \hbar_2}(M):=B^{\bcdot, \bcdot}(M)\otimes L=\big(A^{\bcdot, \bcdot}(M)\big)^\#\otimes_\cF\Vol\otimes L,
\ee
so that the pairing takes values in the volume: the space dual to $B^{\bcdot, \bcdot}_{\hbar_1, \hbar_2}(M)$ is $A^{\bcdot, \bcdot}_{-\hbar_1, -\hbar_2}(M)$, since the dual of $(L, \nabla)$ is the bundle with the connection form $-\alpha_{(\hbar_a)}$, i.e., with $\hbar_a$ replaced with~$-\hbar_a$. Two remarks are in order:

$1)$ If $\cM$ is \textbf{special} periplectic, then, thanks to identification~\eqref{ident}, we have $\Vol\simeq\cF$ and hence $B^{\bcdot, \bcdot}_{\hbar_1, \hbar_2}(\cM)\simeq A^{\bcdot, \bcdot}_{\hbar_1, \hbar_2}(\cM)$ up to a~shift of the bi-degree, and one and the same algebra of invariant operators acts on both spaces --- this is precisely what happens for one parameter $\hbar$ in Subsec.~\ref{fgMyst}, where the operators $\cD^\pm$ act on $\widehat \Omega_\hbar(\cM)$ and on $\widehat \Sigma^\hbar(\cM)$ alike. If $\cM$ is only \textbf{general} periplectic, no such identification exists, see Subsec.~\ref{ss2types} and Remark~\ref{whyspecial}. 

Then, $B^{\bcdot, \bcdot}_{\hbar_1, \hbar_2}(\cM)$ and $A^{\bcdot, \bcdot}_{\hbar_1, \hbar_2}(\cM)$ are two \textit{non-isomorphic} modules over one and the same algebra $A^{\bcdot, \bcdot}(\cM)$, and no shift of the bi-degree makes them isomorphic: by eq.~\eqref{Btwist} they differ by the factor $\Vol$, and $\Vol$ is a~non-trivial line bundle on a~general periplectic $\cM$, because $\Ber$ is a~non-trivial character of $\Pe(n)$, see Remark~\ref{whyspecial}. The twist by $(L, \nabla)$ is the \textbf{same} on both sides and therefore cannot repair this: it is the untwisted $B^{\bcdot, \bcdot}(\cM)$ that already differs from~ $A^{\bcdot, \bcdot}(\cM)$.

$2)$ In presence of $\omega_\od$, the parameters $\hbar_a$ are \textit{odd} Grassmann indeterminates supercommuting with all the operators in sight, since the \lq\lq time" coordinate $\tau$ of the connection is odd, see Subsec.~\ref{fgMyst}; hence $\hbar_a^2=0$ and the twist is infinitesimal. For one such parameter, the answer --- the largest algebra of invariant operators acting on $\widehat \Omega_\hbar$ and $\widehat\Sigma^\hbar$ --- is Theorem~\ref{thP7}; for two of them it is Statement~\ref{thP7two}, computed in Subsec.~\ref{fgMyst2}.

\ssec{Almost K\"ahler periplectic supermanifolds (expansion of \cite{L2})}\label{Sperip}

This is where the two cases of Subsec.~\ref{p3} --- the tensors $J$ and $h$ of opposite parities --- are treated, both of them at once. Indeed, $p(\omega)=p(h)+p(J)$, so opposite parities of $h$ and $J$ say precisely that $\omega$ is odd, i.e., that $\cM$ is periplectic, and the pair $(J, h)$ enters the answer only through the twisted differential $d_c$ of eq.~\eqref{d_c}. Accordingly, the Lie superalgebra acting in the space $A^{\bcdot, \bcdot}(\cM)$ and in the space of integral forms dual to it, denote it $B^{\bcdot, \bcdot}(\cM)$, is described in two steps: the operators $\cX^{\pm}$ and $d$, which do not involve $J$, generate
\be\label{ADolbOdd}
\fh:=\fb\fa(2)\inplus\fsl(1|1),\ \ \text{see analogous constructions in Subsec.~\ref{p1}, eq.~\eqref{ADolb}},
\ee
see Lemma~\ref{gl11}, whereas the Dolbeault operators generate a~larger Lie superalgebra, see Statement~\ref{stmixed}.

\sssec{An action of $\fb\fa(2)\inplus \fsl(1|1)$ on the spaces $\widehat \Omega(\cM)$ and $\widehat \Sigma(\cM)$ on any special almost periplectic~ $\cM$} Let $\cM$ be a \textbf{special} almost periplectic supermanifold with almost periplectic form~$\omega$.  On the spaces $\widehat \Omega(\cM)$ and $\widehat \Sigma(\cM)$, 
\be\label{cXpm}
\begin{array}{l}
\text{let $\cX^+$ be the operator of left multiplication by $\omega$,}\\
\text{let $\cX^-$ be the convolution with the bivector~$\omega^{-1}$, and
$\delta:=[d, \cX^-]$.}
\end{array}
\ee
Here, contrary to what one might expect by analogy with eq.~\eqref{Del},
\be\label{Svanishes}
\text{the operator}\ \ [d, \delta]\ \ \text{\textit{vanishes identically}}.
\ee
Indeed, $\delta=[d, \cX^-]$ and $d$ is odd with $d^2=0$, so the Jacobi identity gives
\[
[d,[d, \cX^-]]=\nfrac12[[d,d], \cX^-]=[d^2, \cX^-]=0.
\]
The analogy with eq.~\eqref{Del} fails for a~reason worth spelling out, since the letter $\delta$ is the same in the two places whereas the operator is not. In eq.~\eqref{H}, $\delta$ is the codifferential dual to $d$ with respect to $\omega_\ev$, and the $H$ it produces is the \textit{mixed} second derivative printed there; the operator $d^*$ is the dual with respect to~$h$, and the $\Delta$ it produces, see eq.~\eqref{Del}, is a~\textit{sum of squares}. The two Laplacians are therefore different operators, and neither is a~multiple of the other$:$ in the flat model of Subsec.~\ref{flat} over $\Cee^{2m}$, verified for $m=1, 2, 3$,
\be\label{HvsDel}
H=\sum_{1\leq i\leq m}\partial_{q_i}\partial_{p_i},\ \ \text{whereas}\ \ \Delta=\sum_{1\leq i\leq m}\big(\partial_{q_i}^2+\partial_{p_i}^2\big).
\ee
This is why the $\Delta$ of eq.~\eqref{Del} does not survive the twist whereas the $H$ of eq.~\eqref{H} does, see eq.~\eqref{Deltw} and the paragraph after it. The $\delta$ of eq.~\eqref{cXpm} is a~third operator$:$ in the Dolbeault picture of Subsec.~\ref{p1} it is $-d_c^*$, the dual with respect to $h$ of the twisted differential $d_c$ of eq.~\eqref{d_c} --- whence the vanishing \eqref{Svanishes} just proved is the same identity as the isotropy of $\Span(d, d_c^*)$ asserted for $\fhei(0|4)_\od$ after eq.~\eqref{Del} --- neither of them an extra input, both the Jacobi identity above.
Hence, the operators $d$ and $\delta$ alone span an \textbf{abelian} $(1|1)$-dimensional superalgebra, not~$\fb\fa(1)$. The anti-bracket superalgebra does appear, but one has to go one step further: setting
\be\label{dprime}
d':=[\cH, d],\ \ \delta':=[\cX^-, d'],\ \ \cS:=-\nfrac12[d', \delta]=\nfrac12[d, \delta'],
\ee
we see that $\cS$ spans the odd center, that $\Span(d, \delta', \cS)$ is a~copy of $\fb\fa(1)$, and that the whole of $\fb\fa(2)$, see eq.~\eqref{ba}, is spanned by $d$, $d'$, $\delta$, $\delta'$ and~$\cS$. This $\cS$ is the odd Laplacian of eq.~\eqref{DelS}. 
The letters are worth keeping apart$:$ 
$\Delta$ is reserved for the \textit{even} Laplacian.

\parbegin[The operators $\cX^{\pm}$ and $d$ generate $\fb\fa(2)\inplus\fsl(1|1)$]{Lemma}[The operators $\cX^{\pm}$ and $d$ generate $\fb\fa(2)\inplus\fsl(1|1)$]\label{gl11} For any special almost periplectic $\cM$, on the spaces $\widehat \Omega(\cM)$ and $\widehat \Sigma(\cM)$, the operators $\cX^+$ and $\cX^-$ generate an action of $\fsl(1|1)$ and, together with the operator $d$, generate an action of
\[
\fb\fa(2)\inplus \fsl(1|1),\ \ 
\text{see eq.~\eqref{ADolbOdd}},
\]
of which $\fb\fa(1)\inplus\fsl(1|1)$ 
is a~subalgebra.
\end{Lemma}

\begin{proof} In local coordinates $x=(q,\pi)$ with $(q_1, \dots , q_n)$ even and  $(\pi_1, \dots , \pi_n)$  odd, we have
\begin{equation}\label{few}
\begin{split}
\cX^-:= \sum\partial_{\widehat
 \pi_i} \partial_{\widehat q_i}; \ \
\cH:=[\cX^+, \cX^-]=n-\sum(\widehat q_i\partial_{\widehat q_i} - \widehat \pi_i\partial_{\widehat \pi_i});\ \ \cX^+:=\sum\widehat q_i\widehat \pi_i.\qed
\end{split}
\end{equation}
\noqed\end{proof} 

\parbegin[Why the analog of $\fhei(0|4)$ is $\fb\fa(2)$, not $\fb\fa(1)$]{Remark}[Why the analog of $\fhei(0|4)$ is $\fb\fa(2)$, not $\fb\fa(1)$]\label{whyba2}
Note that $\fb\fa(2)$, not $\fb\fa(1)$, is the correct analog of the Heisenberg Lie superalgebra $\fhei(0|4)$ of Subsec.~\ref{p1}: both are of dimension $4$ modulo their centers. The same count holds for the larger Lie superalgebra of Statement~\ref{stmixed}, where it is the eight Dolbeault differentials that survive modulo the center.
\end{Remark}

\sssec{The Dolbeault operators add more}\label{ssmixed} Lemma~\ref{gl11} uses only $\omega$: the operators there are $\cX^{\pm}$, see eq.~\eqref{cXpm}, and $d$, and $J$ never enters. On the Dolbeault algebra $A^{\bcdot, \bcdot}(\cM)$, however, the operator $d_c=JdJ^{-1}$ of eq.~\eqref{d_c} acts as well --- it is invariant under the diffeomorphisms preserving both $\omega$ and $J$, which is what an almost K\"ahler structure supplies --- and it is \textit{not} in the $\fh$ of eq.~\eqref{ADolbOdd}. This is where the two cases of Subsec.~\ref{p3} part company.

\sssbegin{Statement}[The \lq\lq mixed-parity\rq\rq\ Dolbeault answer]\label{stmixed} Let $\cM$ be special almost periplectic, so that $\omega$ is odd, and let $J$ be such that $J^2$ is a~non-zero scalar and $J$ preserves $\omega$; let $h:=\omega(J^{-1}\ -, -)$, whose parity is then automatically the one opposite to that of~$J$. 
Then, in the flat model of Subsec.~$\ref{flat}$ on $\Cee^{n|n}$ $($the answer does not depend on $n$, nor on the choice of such a~$J)$, the Lie superalgebra generated by $\cX^{\pm}$, $d$ and $d_c$ is
$\fn\inplus\fsl(1|1)$,
where $\fsl(1|1)=\Span(\cX^+, \cX^-, \cH)$, 
and $\fn$ is the ideal generated by $d$ and $d_c$: it contains the $\fb\fa(2)$ of Lemma~$\ref{gl11}$, it is nilpotent of degree $2$, and $\fn/\fz(\fn)$ is of $\sdim=(4|4)$ --- the eight differentials --- whereas its center $\fz(\fn)$, spanned by the Laplacians, is of $\sdim=(1|1)$ for $J$ odd and $h$ even, and of $\sdim=(0|2)$ for $J$ even and $h$ odd. Hence,
\[
\sdim\big(\fn\inplus\fsl(1|1)\big)=\begin{cases}(6|7)&\text{for $J$ odd, $h$ even},\\
(5|8)&\text{for $J$ even, $h$ odd}.
\end{cases}
\]
In both cases, the operator $d_c$ satisfies $d_c^2=0$, but the parity of $d_c$ is opposite to that of~$J$. 
\end{Statement}

The evidence is a~computation in the flat model, exact over $\Qee$; 
both mixed cases need $n$ even, since for $J$ odd the form $h$ is even and non-degenerate, hence symplectic on the odd part, whereas for $J$ even the condition $J^2=\varkappa\,\id$ needs blocks of even size.

For $J$ even and $h$ odd, there is thus \textbf{no} even Laplacian at all: both elements of $\fz(\fn)$ are odd, one of them being the odd Laplacian $\cS=\sum_{1\leq i\leq n}\partial_{q_i}\partial_{\pi_i}$ of eq.~\eqref{dprime}. For $J$ odd and $h$ even, the center is a~pair consisting of that same $\cS$ and of the even second-order operator $\sum(\partial_{q_i})^2$. This is once again the phenomenon recorded after Statement~\ref{p2b} --- odd Laplacians surviving while even ones die, or the other way round --- and it is what makes eq.~\eqref{ADolbOdd} the answer for $\widehat \Omega(\cM)$ and $\widehat \Sigma(\cM)$ rather than for the Dolbeault algebra.

\sssec{The Lie superalgebra acting on the spaces $\widehat \Omega_\hbar(\cM)$ and $\widehat \Sigma^\hbar(\cM)$ on a~ special almost periplectic $\cM$}\label{fgMyst} Let
$\widehat \Omega_\hbar(\cM):=\widehat \Omega(\cM)\otimes L$ and $\widehat \Sigma^\hbar(\cM):=\widehat \Sigma(\cM)\otimes L$ be the spaces of pseudo-forms with values in the line bundle $L$ with connection whose \lq\lq time" coordinate $\tau$ is odd and the 
form $\alpha_{(\hbar)}$ of the connection $\nabla$ is such that $d\alpha_{(\hbar)}=\hbar\omega_\od$ for any odd Grassmann parameter $\hbar$, supercommuting with all the operators in sight and hence such that~$\hbar^2=0$. In local coordinates $x:=(q, \pi; \tau)$, we can consider
\[
\alpha_{(\hbar)} := d\tau+\frac{\hbar}{2} \sum (\pi_idq_i - q_id\pi_i).
\]
Then, the following operators $\cD^+$ and $\cD^-$ (compare with eq.~\eqref{Dpm}) act on $\widehat \Omega_\hbar(\cM)$ and $\widehat \Sigma^\hbar(\cM)$:
\begin{equation}\label{Dpm2}
\begin{split}
 &\cD^+:= \sum\widehat x_i\partial_{x_i}- \frac{\hbar}{2}\sum
 (\widehat q_i\pi_i - \widehat \pi_iq_i),\\
 &\cD^- :=[\cX^-, \cD^+]=\sum\left(\partial_{\widehat \pi_i}\partial_{q_i} + \partial_{\widehat
 q_i}\partial_{\pi_i}\right ) - \frac{\hbar}{2}\sum\left( q_i \partial_{\widehat q_i}-\pi_i\partial_{\widehat
 \pi_i}\right).\\
\end{split}
\end{equation}
The operator $\cD^-$ is defined, as in eq.~\eqref{Dpm}, as the bracket of $\cD^+$ with $\cX^-$; the right-hand side above is that bracket computed in the flat model of Subsec.~\ref{flat}. Note: $\cD^-|_{\hbar=0}=\delta:=[d, \cX^-]$, as it should be, see eq.~\eqref{Svanishes}.

The Lie superalgebra these operators generate 
is described by the following theorem, computed in the flat model of Subsec.~\ref{flat} with $\hbar$ adjoined as an odd parameter supercommuting with all the operators; accordingly, $\Cee[\hbar]$ below is the Grassmann algebra on one generator, no relation $\hbar^2=0$ having to be imposed by hand.

\parbegin[The algebra generated by $\cD^{\pm}$ and $\cX^{\pm}$ is of $\sdim=(8|8)$]{Theorem}[The algebra generated by $\cD^{\pm}$ and $\cX^{\pm}$ is of $\sdim=(8|8)$]\label{thP7} The Lie superalgebra $\fg$ generated by the operators $\cD^{\pm}$, see eq.~\eqref{Dpm2}, and $\cX^{\pm}$, see eq.~\eqref{few}, is finite-dimensional, 
$\sdim \fg=(8|8)$ independently of $n$; in fact, $\fg$ is generated only by $\cD^+$ and $\fsl(1|1)$, generated by~$\cX^{\pm}$. 
The structure of $\fg$ is as follows.

\emph{$1)$ Square root relations.} The operator $\cD^+$ is odd and $\cD^-$ is even; therefore
\[
[\cD^+, \cD^+]=-2\hbar\, \cX^+,\ \ \text{whereas}\ \ [\cD^-, \cD^-]=0,\ \ \text{and}\ \ [\cD^+, \cX^+]=0=[\cD^-, \cX^-].
\]
There is, however, a~\textbf{second} odd square root in $\fg$, namely,
\be\label{Dsecond}
\widetilde\cD^+:=[\cD^-, \cX^+]=\sum
 (\widehat q_i\partial_{q_i} - \widehat \pi_i\partial_{\pi_i})- \frac{\hbar}{2}\sum
 (\widehat q_i\pi_i + \widehat \pi_iq_i),\ \ \text{so}\ \ [\widetilde\cD^+, \widetilde\cD^+]=2\hbar\, \cX^+.
\ee
It is \textbf{not} obtained from $\cD^+$ by the change $\hbar\mapsto-\hbar$: already its $\hbar=0$ part is not $d$ but $-d'=[d, \cH]$, see eq.~\eqref{dprime}.

\emph{$2)$ The bracket of the two $\cD$'s is a~multiple of the center of $\fsl(1|1)$.} We have, exactly,
\[
[\cD^+, \cD^-]=-\hbar\,\cH=-\hbar\,[\cX^+, \cX^-],
\]
so no Laplacian occurs here: parity forces the right-hand side to be odd, and $\hbar\cH$ is what it is. This is the periplectic analog of the relation $[D^+, D^-]\ \propto\  H$ that holds in the symplectic case of eq.~\eqref{Dpm}, where $D^\pm$ are the two odd root vectors of~$\fosp(1|2)$.

\emph{$3)$ Extension structure.} Setting $\hbar=0$ sends $\fg$ onto the Lie superalgebra $\fb\fa(2)\inplus\fsl(1|1)$ of Lemma~$\ref{gl11}$, 
see eq.~\eqref{ADolbOdd}, and the kernel is the abelian ideal $\hbar\fg$ of $\sdim=(5|3)$, so
\[
0\tto \hbar\fg\tto \fg\tto \fb\fa(2)\inplus\fsl(1|1)\tto 0,\ \ \text{where}\ \ \hbar\fg\simeq\Pi\big(\fb\fa(2)\inplus\fsl(1|1)\big)
\]
as a~module over $\fb\fa(2)\inplus\fsl(1|1)$ $($the isomorphism being $u\mapsto \hbar u$, well defined since~$\hbar^2=0)$. This extension does \textbf{not} split: no subalgebra of $\fg$ maps isomorphically onto $\fb\fa(2)\inplus\fsl(1|1)$, and hence
\[
\fg\not\simeq\big(\fb\fa(2)\inplus\fsl(1|1)\big)\otimes\Cee[\hbar].
\]

\emph{$4)$} The subalgebra $\fsl(1|1)\otimes\Cee[\hbar]$ 
spanned by $\cX^{\pm}$, $\cH$, see eq.~\eqref{few}, and their $\hbar$-multiples, is indeed contained in~$\fg$. \end{Theorem}

\begin{proof} Headings $1)$, $2)$ and $4)$ are single brackets of the operators written out in eqs.~\eqref{Dpm2} and~\eqref{few}, and are verified directly.

Heading $3)$: the extension does not split. Recall that here $\omega_\od$ is \textit{odd}, so $\cX^+$ and $d$ are both odd, whereas $\delta$ and $\cH$ are even. Suppose a~splitting exists, i.e., suppose the projection $\pi\colon\fg\tto\fh:=\fb\fa(2)\inplus\fsl(1|1)$ admits a~parity-preserving \textit{section} --- a~linear map $\sigma\colon\fh\tto\fg$ right inverse to it, $\pi\circ\sigma=\id_\fh$ --- which is moreover a~homomorphism, $[\sigma(u),\sigma(v)]=\sigma([u,v])$ for all~$u, v\in\fh$. Since $\hbar\fg$ is the kernel of~$\pi$, $\sigma(d)=\cD^++\hbar w$ for some $w\in\fh$, and $w$ is \textit{even} because $\hbar$ and $d$ are odd. As $\hbar^2=0$,
\[
[\sigma(d), \sigma(d)]=[\cD^+, \cD^+]+2\hbar[d, w]=-2\hbar\big(\cX^+-[d, w]\big),
\]
by heading~$1)$. On the other hand, $d^2=0$, hence $[d,d]=0$ and $\sigma([d,d])=0$. So a~splitting forces
\[
\cX^+=[d, w]\ \ \text{for some even}\ w\in\fh.
\]
But $\fh_\ev=\Span(\delta, \delta', \cH)$, whereas a~direct computation of the three brackets gives
\[
[d, \fh_\ev]=\Span(d', \cS),\ \ \text{of}\ \ \sdim=(0|2),
\]
which does not contain~$\cX^+$. This contradiction proves that the extension does not split. 

The same argument, read backwards, says what the obstruction \textit{is}. The extension of heading~$3)$ is described by a~class in $H^2(\fh; \Pi\fh)$, where $\fh:=\fb\fa(2)\inplus\fsl(1|1)$ and the coefficients are the kernel $\hbar\fg\simeq\Pi\fh$; a~cocycle $c$ representing this class is read off any parity-preserving section~$\sigma$ of the projection $\fg\tto\fh$ by
\[
[\sigma(u), \sigma(v)]-\sigma([u,v])=\hbar\, c(u,v)\quad\text{for $u, v\in\fh$},
\]
which determines $c(u,v)\in\fh$, the map $u\mapsto\hbar u$ being one-to-one. The computation above is precisely the value of $c$ at the single pair $(d, d)${}$:$ for $\sigma(d)=\cD^+$ it gives $c(d, d)=-2\cX^+$, whereas replacing $\sigma(d)$ with $\cD^++\hbar w$, where $w\in\fh_\ev$ --- which is the only freedom in the choice of $\sigma$ --- changes that value by $2[d, w]\in[d, \fh_\ev]=\Span(d', \Delta)$. Since $\cX^+\notin\Span(d', \Delta)$, no section makes $c(d, d)$ vanish$:$ the class is non-zero, and it is non-zero already at that one pair. 
\end{proof}

\sssec{Two odd parameters}\label{fgMyst2} The same computation runs for two odd parameters, and the answer is not the doubling one might expect. Here, $\hbar_1\hbar_2\neq0$, although $\hbar_a^2=0$ for $a=1, 2$, so the twist is no longer infinitesimal, and there are two operators $\cD^+_a$, one for each connection form
\be\label{alphatwo}
\alpha_{(\hbar_a)} := d\tau+\frac{\hbar_a}{2} \sum (\pi_idq_i - q_id\pi_i),\ \ \text{so that}\ \ d\alpha_{(\hbar_a)}=\hbar_a\omega_\od, 
\ee
both acting on
the same pair of spaces $\widehat \Omega_{\hbar_1, \hbar_2}(\cM)$ and $\widehat \Sigma^{\hbar_1, \hbar_2}(\cM)$; as before, ${\cD^-_a:=[\cX^-, \cD^+_a]}$. On a~special almost periplectic $\cM$, there is only one $\omega_\od$ to serve as the curvature form of both connections, which is why the two forms of eq.~\eqref{alphatwo} differ only in the parameter. Let 
\[
\text{$\fG$ be the Lie superalgebra  generated by the operators $\cD^{\pm}_a$ for $a=1, 2$ and  $\cX^{\pm}$.}
\]

\sssbegin[Two odd parameters give $\sdim\fG=(22|22)$, not the expected $(16|16)$]{Statement}[Two odd parameters give $\sdim\fG=(22|22)$, not the expected $(16|16)$]\label{thP7two} Computed in the flat model of Subsec.~$\ref{flat}$; 
the same computation, given one parameter instead of two, returns the $\sdim\fg=(8|8)$ of Theorem~$\ref{thP7}$.

We have $\sdim \fG=(22|22)$ independently of~$n$. The Lie superalgebra $\fG$ is solvable, 
its center is of~$\sdim=(1|1)$.
Since the $\hbar_a$ are central, the subspaces $\fG^{[k]}$ of the elements divisible by $k$ of the parameters are ideals; we have
\[
\fG\supset\fG^{[1]}\supset\fG^{[2]}\supset 0,\ \ \text{of}\ \ \sdim=(22|22),\ (19|17),\ (7|7),
\]
with quotients of $\sdim=(3|5)$, $(12|10)$ and~$(7|7)$. The top quotient is again the $\fb\fa(2)\inplus\fsl(1|1)$ of Lemma~$\ref{gl11}$: setting $\hbar_1=\hbar_2=0$ returns eq.~\eqref{ADolbOdd}, exactly as for one parameter.

Two features tell the two-parameter case from the one-parameter one. First, the ideal $\fG^{[1]}$ is \textbf{not} abelian, whereas the $\hbar\fg$ of Theorem~$\ref{thP7}\,3)$ is. Second, $\fG$ is larger than $\fh\otimes\Lambda(\hbar_1, \hbar_2)$, where $\fh:=\fb\fa(2)\inplus\fsl(1|1)${}$:$ this $\fh\otimes\Lambda(\hbar_1, \hbar_2)$ is of $\sdim=(16|16)$, which is the doubled $\sdim\fg$ of Theorem~$\ref{thP7}$ and the superdimension the heading calls expected. We see that $\fG$ contains $\hbar_1\fh$, $\hbar_2\fh$ and $\hbar_1\hbar_2\fh$, but not $\fh$ itself.
\end{Statement}

\ssec{Primitive forms. Lee's theorem}\label{PrimFLee} 
The primitive forms are what the two \lq\lq halves\rq\rq\ of $d\omega$ are made of. Recall that a~$3$-form on a symplectic manifold is said to be \textit{primitive} if $X^-$, see eq.~\eqref{Xpm}, kills it, i.e., if it is not divisible by $\omega$; on a~manifold $M^{2a}$, every $3$-form is the sum of a~primitive one and of a~multiple of $\omega$, so the condition $d\omega=0$ is a~system of \textit{two} equations, and one may impose either of them alone, see \cite{Lee1}. Here is what this splitting turns into on supermanifolds, for either parity of~$\omega$.

\sssbegin[Where the two halves of $d\omega=0$ fuse]{Lemma}[Where the two halves of $d\omega=0$ fuse]\label{halves} For $X^{\pm}$, see eq.~\eqref{Xpm}, let
\[
P^3:=\{\alpha\in\Omega^3\mid X^-\alpha=0\}
\]
be the space of primitive $3$-forms.

$1)$ Let $\omega$ be even, $\sdim\cM=(2a|b)$. If $2a-b\neq2$, then
\be\label{P3}
\Omega^3=P^3\oplus X^+(\Omega^1).
\ee
If $2a-b=2$, then $X^+(\Omega^1)\subset P^3$, and hence, for $b>0$, the sum $P^3+X^+(\Omega^1)=P^3$ is a~\textit{proper} subspace of $\Omega^3$ $($for $b=0$ this is the classical degeneration 
in which~$\Omega^3=0)$.

$2)$ Let $\omega$ be odd, $\sdim\cM=(n|n)${}$;$ here, both $X^{\pm}$ are odd, so their superbracket is the anticommutator, and on $\widehat \Omega(\cM)$ and $\widehat \Sigma(\cM)$ they are the $\cX^{\pm}$ of eq.~\eqref{cXpm}. Then, $[X^-, X^+]$ is \textbf{not} a~scalar operator on $\Omega^k$ for $k\geq1$, so there is no analog of the $H_{sl}$ of eq.~\eqref{Hsl}$:$ on a~monomial containing $\alpha$ of the $\widehat q_i$ and $\beta$ of the~$\widehat \xi_j$, the operator $[X^-, X^+]$ is the multiplication by $n-\alpha+\beta$. Nevertheless, decomposition~\eqref{P3} holds for every $n\geq1$, although for $n=1$ the operator $X^+$ is not injective on~$\Omega^1$.
\end{Lemma}

\begin{proof} Since $\Omega^{-1}=0$, we have $X^-|_{\Omega^1}=0$, so on $\Omega^1$ the composition $X^-X^+$ is the superbracket $[X^-, X^+]$ --- the commutator if $\omega$ is even, the anticommutator if $\omega$ is odd. Whenever this operator is invertible on $\Omega^1$, the operator $X^+$ is injective on $\Omega^1$, the operator $X^-\colon \Omega^3\tto\Omega^1$ is onto, and hence $\dim P^3=\dim \Omega^3-\dim \Omega^1$ while $P^3\cap X^+(\Omega^1)=0$; the two subspaces are therefore complementary. Whenever it vanishes on a~vector $\alpha\in\Omega^1$, we get $X^-X^+\alpha=0$, i.e., $X^+\alpha\in P^3$.

$1)$ By the formula for $H_{sl}$, see eq.~\eqref{Hsl}, $[X^-, X^+]$ acts on $\Omega^1$ as the scalar $1-a+\nfrac12 b$, which vanishes if and only if~$2a-b=2$.

$2)$ In the flat model of Subsec.~\ref{flat}, with $\widehat q_i$ odd and $\widehat \xi_i$ even, we have $X^+=\sum_{1\leq i\leq n}\widehat q_i\widehat \xi_i$ and $X^-=\sum_{1\leq i\leq n} \partial_{\widehat \xi_i}\partial_{\widehat q_i}$. The $n$ pairs act on disjoint sets of indeterminates, so it suffices to consider $n=1$, where relations~\eqref{engineRel} give $[X^-, X^+]=1-\alpha+\beta$ on $\widehat q^{\,\alpha}\widehat \xi^{\,\beta}$; summation over $i$ yields the claim. On $\Omega^1$, this operator is therefore the multiplication by $n-1$ on the $\widehat q_i$ and by $n+1$ on the $\widehat \xi_i$; for $n\geq2$, it is invertible and the first paragraph applies. For $n=1$, we have $\sdim\Omega^1=(1|1)$ and $\sdim\Omega^3=(1|1)$, whereas $X^+\widehat q=\widehat q\widehat \xi\widehat q=0$ and $X^-(\widehat q\widehat \xi^{\,2})=2\widehat \xi\neq0$, so $P^3=\Cee\widehat \xi^{\,3}$ and $X^+(\Omega^1)=\Cee \widehat q\widehat \xi^{\,2}$. These are two $1$-dimensional subspaces of opposite parities --- the generator $\widehat \xi^{\,3}$ is even and $\widehat q\widehat \xi^{\,2}$ is odd, the $\widehat q$ being odd and the $\widehat \xi$ even in this model --- whereas $\sdim\Omega^3=(1|1)$; hence, they do span $\Omega^3$.\end{proof}

\sssbegin[The smallest supermanifold on which the halves fuse]{Example}[The smallest supermanifold on which the halves fuse]\label{fuse} Let $\sdim\cM=(4|2)$, so that~$2a-b=2$. Then, $\sdim\Omega^1=(2|4)$ and $\sdim\Omega^3=(16|16)$, whereas $\sdim P^3=(14|12)$ and $\sdim X^+(\Omega^1)=(2|4)$, the latter being contained in the former. So the $3$-form $d\omega$ has here no component divisible by $\omega$ separate from its primitive component$:$ the subspace $X^+(\Omega^1)$, whose elements are the $3$-forms divisible by~$\omega$, lies \textit{inside} the subspace $P^3$ of the primitive ones instead of being complementary to it, so the sum of the two is not direct and neither component of $d\omega$ is defined. What survives are the two membership conditions --- that $d\omega$ be primitive, and that $\omega$ divide $d\omega$ --- but the two are no longer independent$:$ the second implies the first, and their conjunction is no longer~$d\omega=0$. Accordingly, the two classes of supermanifolds singled out in Subsec.~\ref{ssTwoCond}, the locally conformally symplectic ones and the balanced ones, do not separate here. By Lemma~\ref{halves}, the same happens on every supermanifold of $\sdim\cM=(2a|2a-2)$ with $a\geq2$, and on no other one. 
\end{Example}

Lee asked when a~non-degenerate $2$-form $\omega$ admits an~\textit{integrating factor}, i.e., a~function $f$ such that $f\omega$ is closed. This is a~question not about the vanishing of $d\omega$ but about the vanishing of one half of it, and his answer, in the language above, reads$:$ on a~manifold of dimension $>4$, such an $f$ exists locally if and only if $\omega$ divides $d\omega$, see \cite{Lee1}. The $1$-form $\theta$ in the relation $d\omega=\theta\wedge\omega$ is what is now called the \textit{Lee form}, and the manifolds on which that relation holds are the locally conformally symplectic ones; the sequels \cite{Lee2, Lee3} describe the automorphisms of $\omega$ and prove that two such structures of the same rank are locally equivalent, which is Darboux's theorem. The analog of these statements on supermanifolds is as follows.

\sssbegin[Lee's theorem on supermanifolds]{Proposition}[Lee's theorem on supermanifolds]\label{lee} Let $\omega$ be a~non-degenerate $2$-form of either parity on a~supermanifold $\cM$, let $X^+$ be the multiplication by $\omega$, and let $f$ run over the invertible even functions.

$1)$ Whatever the parity of $\omega$, a~$1$-form $\theta$ satisfying $d\omega=\theta\wedge\omega$ is odd in the total parity of the superalgebra~$\Omega(\cM)$. Hence, $\theta\wedge\theta=0$ and
\begin{equation}\label{dtheta}
d\theta\wedge\omega=0.
\end{equation}

$2)$ Such a~$\theta$ is unique as soon as $X^+\colon\Omega^1\tto\Omega^3$ is injective, which is the case on every $\cM$ except those of $\sdim=(2|0)$ and~$(1|1)$.

$3)$ Let $\omega$ be even, $\sdim\cM=(2a|b)$ with~$b>0$. Then, $X^+$ is injective on $\Omega^k$ for every $k$; hence $d\theta=0$, and
\[
\text{$d(f\omega)=0$ for some $f$, locally}\quad\Longleftrightarrow\quad\text{$\omega$ divides $d\omega$}.
\]
The same equivalence holds for $b=0$ and $2a\geq6$, which is Lee's own restriction on the number of variables, and fails for $\sdim\cM=(4|0)$. So that restriction is void on any supermanifold of superdimension $(2a|b)$ with~$b\neq0$. 

$4)$ Let $\omega$ be odd, $\sdim\cM=(n|n)$. Then, $X^+$ is \textbf{never} injective on $\Omega^2${}$:$ the form $\omega$ is odd, so $\omega\wedge\omega=0$, and
\[
\Ker\left(X^+|_{\Omega^2}\right)=\Cee\omega\ \ \text{for $n\neq2$},\qquad \Ker\left(X^+|_{\Omega^2}\right)=\Cee\omega\oplus\Cee\, \widehat q_1\widehat q_2\ \ \text{for $n=2$}.
\]
Accordingly, eq.~\eqref{dtheta} gives not $d\theta=0$, but $($for some functions $\varphi$  and $\psi$\,$)$
\[
d\theta=\varphi\omega\ \ \text{for $n\neq2$},\qquad d\theta=\varphi\omega+\psi\,\widehat q_1\widehat q_2\ \ \text{for $n=2$, where $p(\varphi)=\od$  and $p(\psi)=\ev$}.
\]

$5)$ The function $\varphi$ of item $4)$ vanishes for every $n$, see Lemma~\textup{\ref{nilLee}}. Hence, for $n\neq2$, the Lee form of an odd $\omega$ is closed after all, and the equivalence of item $3)$ holds verbatim on every $\cM$ of $\sdim=(n|n)$ with $n\neq2$, Lee's own restriction on the number of variables disappearing in the odd case as well. 

For $n=2$, that equivalence \textbf{fails}, see Example~\textup{\ref{noint}}$:$ the extra summand $\psi\,\widehat q_1\widehat q_2$ can be non-zero. 
\end{Proposition}

So the superdimensions on which Lee's theorem does not hold are $(4|0)$ for an even $\omega$, where the divisibility of $d\omega$ by $\omega$ is a~vacuous demand, and $(2|2)$ for an odd one, where it is not the divisibility but the closedness of $\theta$ that fails --- and fails through $\psi$ alone.

\begin{proof}~$1)$. The operator $d$ is odd, so $p(d\omega)=p(\omega)+1$, whereas $p(\theta\wedge\omega)=p(\theta)+p(\omega)$; hence $p(\theta)=1$ regardless of $p(\omega)$, and an odd element of a~supercommutative superalgebra squares to zero. Now, apply $d$ to $d\omega=\theta\wedge\omega${}$:$
\[
0=dd\omega=d\theta\wedge\omega-\theta\wedge d\omega=d\theta\wedge\omega-\theta\wedge\theta\wedge\omega=d\theta\wedge\omega.
\]

$2)$ and~$3)$. If $X^+$ is injective on $\Omega^1$, the difference of two solutions of $d\omega=\theta\wedge\omega$ is killed by $X^+$, hence vanishes. In the flat model of Subsec.~\ref{flat}, the superalgebra of forms is $\Lambda\otimes S$, where $\Lambda$ is generated by the hats of the even coordinates and $S$ by the hats of the odd ones, the hat reversing the parity, so that $\Lambda$ is exterior and $S$ polynomial; the form $\omega=\omega_\Lambda+\omega_S$ splits accordingly, the two summands being of bidegrees $(2,0)$ and~$(0,2)$. Let $b>0$, so that $\omega_S\neq0$, and let $\alpha\in\Omega^k$ satisfy $\omega\alpha=0$. Among the bihomogeneous components of $\alpha$, take the one, $\alpha_{p,q}$, whose degree $q$ in $S$ is the largest. The component of $\omega\alpha$ of bidegree $(p,q+2)$ is $\omega_S\alpha_{p,q}$; the multiplication by $\omega_S$ is injective on $\Lambda\otimes S$ because $S$ is a~polynomial ring, so $\alpha_{p,q}=0$, and induction on $q$ gives~$\alpha=0$. 

For $b=0$, the assertion is the classical one, the map $X^+\colon\Omega^k\tto\Omega^{k+2}$ being injective for~$k\leq a-1$. 

For $\omega$ odd, we have $\omega=\sum_{1\leq i\leq n}\widehat q_i\widehat \xi_i$, whence $\omega\widehat q_k=-\sum_{i\neq k}\widehat q_i\widehat q_k\widehat \xi_i$, which vanishes only for $n=1$, and $\omega\widehat \xi_k\neq0$ for every~$n$.

The equivalence in item $3)$ is now immediate. If $d(f\omega)=0$, then $df\wedge\omega+f\,d\omega=0$, so $\omega$ divides $d\omega$, with $\theta=-f^{-1}df$. Conversely, let $d\omega=\theta\wedge\omega$. By what is just proved, $d\theta=0$, and $\theta$ is odd in the total parity, so locally $\theta=dg$ for an even function~$g$. Then, $f:=e^{-g}$ is invertible, $df=-f\theta$, and $d(f\omega)=-f\theta\wedge\omega+f\theta\wedge\omega=0$.

$4)$. The three components of $\alpha\in\Omega^2$ in the bigrading of $\Lambda\otimes S$ are of bidegrees $(2,0)$, $(1,1)$ and $(0,2)$, and $\omega$ is of bidegree $(1,1)$, so the equation $\omega\alpha=0$ splits into three. On the component of bidegree $(0,2)$, it gives $\widehat \xi_i\alpha=0$ for every $i$, i.e., $\alpha=0$; on the one of bidegree $(2,0)$, it gives $\widehat q_i\wedge\alpha=0$ for every $i$, which for $\alpha\in\Lambda^2$ forces $\alpha=0$ unless $n=2$; on the one of bidegree $(1,1)$, writing $\alpha=\sum_{1\leq i,j\leq n}a_{ij}\widehat q_i\widehat \xi_j$, it gives $a_{ij}=a\delta_{ij}$, i.e., $\alpha\in\Cee\omega$. The parities of $\varphi$ and of $\psi$ are forced by $p(\omega)=1$, $p(\widehat q_1\widehat q_2)=0$ and~$p(d\theta)=0$. 

$5)$. This is Lemma~\ref{nilLee} and Example~\ref{noint}.\end{proof}

\sssbegin[The function $\varphi$ vanishes, whatever $n$]{Lemma}[The function $\varphi$ vanishes, whatever $n$]\label{nilLee} Let $\omega$ be an odd non-degenerate $2$-form on $\cM$ of $\sdim=(n|n)$, let $d\omega=\theta\wedge\omega$, and let $d\theta=\varphi\omega$ for $n\neq2$, resp. $d\theta=\varphi\omega+\psi\,\widehat q_1\widehat q_2$ for $n=2$, be as in item~$4)$ of Proposition~$\ref{lee}$. Then, $\varphi=0$ for every~$n$. Hence, $d\theta=0$ for $n\neq2$, whereas for $n=2$ we have $d\theta=\psi\,\widehat q_1\widehat q_2$, where $\psi$ need not vanish, see Example~$\ref{noint}$.
\end{Lemma}

\begin{proof} Let~$n\geq3$. Applying $d$ to $d\theta=\varphi\omega$ and using $d\omega=\theta\wedge\omega$ together with $p(\varphi)=1$, we get
\[
0=d\varphi\wedge\omega-\varphi\,\theta\wedge\omega=(d\varphi-\varphi\theta)\wedge\omega,
\]
so $d\varphi=\varphi\theta$ by the injectivity of $X^+$ on $\Omega^1$, see item~$2)$. Write $\theta=\sum_{1\leq i\leq n}u_i\widehat q_i+\sum_{1\leq j\leq n}v_j\widehat \xi_j$, where the $u_i$ are even and the $v_j$ odd, see item~$1)$; comparing the $\widehat \xi_j$-components of $d\varphi=\varphi\theta$ in the free module $\Omega^1$, we get
\begin{equation}\label{eulerEq}
\partial_{\xi_j}\varphi=\varphi v_j\quad\text{for $1\leq j\leq n$}.
\end{equation}

Now, let $J\subset\cF(\cM)$ be the ideal generated by the odd coordinates, so ~$J^{n+1}=0$. Being odd, $\varphi$ lies in $J$, and so does each $v_j$; suppose $\varphi\neq0$ and let $k\geq1$ be the largest number with $\varphi\in J^k$. The right-hand side of eq.~\eqref{eulerEq} lies in $J^{k+1}$, whereas $\partial_{\xi_j}$ maps $J^k$ into $J^{k-1}$; hence the image of $\partial_{\xi_j}\varphi$ in $J^{k-1}/J^k$ vanishes for every~$j$. Let $\varphi_k$ be the image of $\varphi$ in~$J^k/J^{k+1}$. Then, $\partial_{\xi_j}\varphi_k=0$ for every $j$, and the Euler operator gives
\[
k\varphi_k=\sum_{1\leq j\leq n}\xi_j\partial_{\xi_j}\varphi_k=0,
\]
so $\varphi_k=0$, i.e., $\varphi\in J^{k+1}$ --- which contradicts the choice of~$k$. Hence,~$\varphi=0$.

Let~$n=2$. Here, item~$4)$ gives $d\theta=\varphi\omega+\psi\,\widehat q_1\widehat q_2$, and applying $d$ to it we get 
\begin{equation}\label{n2Eq}
(d\varphi-\varphi\theta)\wedge\omega+d\psi\wedge\widehat q_1\widehat q_2=0 ,
\end{equation}
whose two summands do \textit{not} separate, the form $\widehat q_i\wedge\omega$ being itself a~multiple of~$\widehat q_1\widehat q_2$. They do separate modulo the ideal generated by $\widehat q_1\widehat q_2$, and it is the parity of the coefficients of $\omega$ that makes this suffice. Indeed, write
\[
d\varphi-\varphi\theta=\sum_{1\leq i\leq 2}c_i\widehat q_i+\sum_{1\leq j\leq 2}e_j\widehat \xi_j,\qquad
\omega=\sum_{1\leq i,j\leq 2}A_{ij}\widehat q_i\widehat \xi_j+B\,\widehat q_1\widehat q_2+\sum_{1\leq i\leq j\leq 2}C_{ij}\widehat \xi_i\widehat \xi_j .
\]
Since $\omega$ is odd whereas $\widehat q_1\widehat q_2$ and $\widehat \xi_i\widehat \xi_j$ are even, the functions $B$ and $C_{ij}$ are \textit{odd}, and hence lie in $J$, while the $A_{ij}$ are even and the matrix $A=(A_{ij})$ is invertible, the form $\omega$ being non-degenerate. Reduce eq.~\eqref{n2Eq} modulo the ideal generated by $\widehat q_1\widehat q_2$. The second summand dies, and so does the product of the $c_i\widehat q_i$ with the $A$-part of $\omega$; what survives in the component containing one $\widehat q$ and two $\widehat \xi$ is a~system expressing the products $e_jA_{kl}$ through the products~$c_kC_{jl}$. The $c_i$ lie in~ $J$, because $\varphi$ does and hence so do $\partial_{q_i}\varphi$ and $\varphi u_i$; the $C_{jl}$ lie in $J$ as well, and $A$ is invertible, whence
\[
e_j\in J^2\quad\text{for $1\leq j\leq 2$}.
\]
Comparing now the $\widehat \xi_j$-components as above, we get $\partial_{\xi_j}\varphi=\varphi v_j+e_j$ in place of eq.~\eqref{eulerEq}. Suppose $\varphi\neq0$ and let $k\geq1$ be the largest number with $\varphi\in J^k$; then $k\leq2$, because $J^3=0$ for $n=2$, and therefore $J^2\subseteq J^k$. So both summands of the right-hand side vanish in $J^{k-1}/J^k$, the first lying in $J^{k+1}$ and the second in $J^2$; the argument above applies verbatim and gives $\varphi=0$ here as well. What it does not give --- and what is false, see Example~$\ref{noint}$ --- is the vanishing of $\psi$, which eq.~\eqref{n2Eq} constrains only through $d\psi\wedge\widehat q_1\widehat q_2$.

Let~$n=1$. Then, $X^+$ is not injective on $\Omega^1$ and the above does not apply, but the claim is immediate$:$ the even functions are the $u(q)$ and the odd ones are the $v(q)\xi$, so $\theta=u\widehat q+v\xi\widehat \xi$ and $\omega=a\widehat q\widehat \xi+c\xi\widehat \xi^{\,2}$ with $a$ invertible, whence
\[
d\theta=v'\xi\,\widehat q\widehat \xi+v\,\widehat \xi^{\,2},\qquad \varphi\omega=wa\xi\,\widehat q\widehat \xi\quad\text{for $\varphi=w(q)\xi$},
\]
the term with $c$ dying because the product of two odd functions of one odd coordinate vanishes. Comparing the coefficients of $\widehat \xi^{\,2}$ gives $v=0$, and then $wa=0$, i.e., $\varphi=0$.\end{proof}

\sssbegin[A form dividing its own differential and having no integrating factor]{Example}[A form dividing its own differential and having no integrating factor]\label{noint} Let $\sdim\cM=(2|2)$, the coordinates being $q_1,q_2$ even and $\xi_1,\xi_2$ odd. Set
\begin{equation}\label{noIF}
\omega=\widehat q_1\widehat \xi_1+\widehat q_2\widehat \xi_2-q_1\xi_1\,\widehat q_1\widehat q_2,\qquad \theta=q_1\widehat q_2.
\end{equation}
Then, $\omega$ is 
the standard odd form modulo the ideal of nilpotents, and
\[
d\omega=\theta\wedge\omega,\qquad\text{whereas}\qquad d\theta=\widehat q_1\widehat q_2\neq0 .
\]
Here, $\varphi=0$ and $\psi=1$, and the obstruction is the summand that exists for $n=2$ alone. Since $X^+$ is injective on $\Omega^1$, the form $\theta$ is the only Lee form of $\omega$; an integrating factor $f$ would give $\theta=-d\ln f$, hence~$d\theta=0$. So $\omega$ divides $d\omega$ and admits no integrating factor.
\end{Example}

\section{Deformations of the Buttin superalgebra (after
\cite{Ko5,LSh1, LSh3})}\label{S2}

\ssec{$\fb(n)$ and its regrading $\fb(n;n)$}\label{bAndbnn}In the standard grading ($\deg\xi_i= \deg q_i=1$ for all $i$), the Buttin algebra $\fb(n)=\mathop{\oplus}_{i\geq -2}\ \fb(n)_i$ is 
of depth~2. Under the regrading given by setting
\[
\text{$\deg\xi_i=0, \deg q_i=1$ for all $i$,}
\] 
the Buttin algebra takes the form
$\fg=\mathop{\oplus}_{i\geq -1}\ \fg_i$ with
$\fg_0=\fvect(0|n)$ and $\fg_{-1}\cong \Pi(\Cee[\xi])$. Denote the regraded Lie superalgebra $\fb (n)$ by~$\fb (n; n)$. 

\ssec{$\fb_{\lambda}(n; n)$}\label{bLambda} Replace now
the $\fvect(0|n)$-module $\fg_{-1}$ of functions in $\xi$ (with inverted
parity) by the rank-$1$ (over the algebra of functions $\Cee[\xi]$) module of densities of weight $\lambda$ (briefly:
$\lambda$-densities), also with inverted parity, i.e., set
$\fg_{-1}\cong \Pi(\Vol^\lambda(0|n))$, where the
$\fvect(0|n)$-action on the odd generator $\vvol_\xi^\lambda$ of the module $\Vol^\lambda(0|n)$ of $\lambda$-densities is given
by the formula 
\begin{equation}
\label{div} 
\begin{array}{l}
L_D(\vvol_\xi^\lambda) :=\lambda \Div D\cdot
\vvol_\xi^\lambda.\\
 \end{array}
\end{equation}

Define $\fb_{\lambda}(n; n)$ to be the Cartan prolong of $(\Pi(\Vol^\lambda(0|n)),
\fvect(0|n))_*$:
\begin{equation}
\label{fb} 
\fb_{\lambda}(n; n):= (\fg_{-1}, \fg_0)_*=(\Pi(\Vol^\lambda(0|n)),
\fvect(0|n))_*.
\end{equation}
Clearly, this is a~deform (the result of a~deformation) of $\fb(n; n)=\fb_{0}(n; n)$. 

Grozman was the first, it seems, who discovered the deform $\fb_{\lambda}(n)$ of $\fb(n)$, see \cite{G, Gpf}, but this example of a~\textbf{simple} Lie superalgebra of vector fields was overlooked by researchers at that time. The Lie superalgebra $\fb_{\lambda}(n)$ is a~regrading of
$\fb_{\lambda}(n; n)$ described as follows. 
Set 
\begin{equation}
\label{2.7.2} 
\renewcommand{\arraystretch}{1.4}
\begin{array}{l}
\fb_{a, b}(n) :=\left\{M_f\in \fm (n)\mid a\; \Div
M_f=(-1)^{p(f)}2(a-bn)\pderf{f}{\tau}\right\}.
\end{array}
\end{equation}
Taking into account the explicit form of the divergence of $M_{f}$ as an element of $\fvect(n|n+1)$, we see that
\begin{equation}
\label{2.7.4}
\renewcommand{\arraystretch}{1.4}
\begin{array}{ll}
\fb_{a, b}(n) &=\{M_f\in \fm (n)\mid (bn-aE)\pderf{f}{\tau} =
a \cS(f)\}=\\
&\{D\in\fvect(n|n+1) \mid L_{D}(\vvol_{q, \xi,
\tau}^a\alpha_{0}^{a-bn})=0\},\end{array}
\end{equation}
where the operator 
\be\label{DelS}
\cS:=\sum\nfrac{\partial^2}{\partial {\xi_i}\partial {q_i}}
\ee 
on $\widecheck M$ (the operator of eq.~\eqref{dprime}, denoted, as in \cite{BV1}, intentionally differently from its even analog~ $\Delta$, see \eqref{Del}) is the image of the Fourier-transform in odd indeterminates $\widehat q:=\xi$ of the exterior differential on the manifold~ $M$ with coordinates~$q$.

A~direct verification shows that $\fb_{a, b}(n)\simeq
\fb_\lambda(n)$ for $\lambda =\frac{2a}{n(a-b)}\in\Cee P^1=\Cee\cup\{\infty\}$. 
In particular, 
\[
\fb (n)=\fb_\lambda(n), \ \ \ \fsm (n)=\fb_{bn, b}(n)=\fb_{2/(n-1)}(n).
\]

The description of topological $\fvect(m|n)$-modules
(\cite{BL}) and the criteria for simplicity of $\Zee$-graded Lie
superalgebras (\cite{K2}) imply simplicity of the Lie superalgebras $\fb_\lambda(n)$
 for $n>1$ and $\lambda\neq 0$, 1,~$\infty$. 
 
\ssec{Simple ideals and quotients}\label{IdQuo} Recall that from Milnor's $K$-functor point of view a~correct  notation of the superdimension of the $(a|b)$-dimensional superspace is $a+\eps b$, where~$\eps^2=1$. At $\lambda=1$ and $\infty$, the Lie superalgebras $\fb_{\lambda}(n)$
are not simple either: each of them has an ideal of codimension
$\eps^{n+1}$ and $\eps^{n}$, respectively, see  exact
sequences \eqref{le}. The exact
sequence for $\fb(n):=\fb_0(n)$ yields $\fle(n)$, see eq.\eqref{Mand}:
\begin{equation}
\label{le}
\renewcommand{\arraystretch}{1.4}
\begin{array}{c}
0\tto \Cee M_{1} \tto \fb(n)\tto \fle(n)\tto 0,\\
0\tto \fb_{1} '(n)\tto \fb_{1}(n)\tto \Cee\cdot
M_{\xi_1\dots\xi_n} \tto 0,\\
0\tto \fb_{\infty} '(n)\tto \fb_{\infty}(n)\tto \Cee\cdot
M_{\tau\xi_1\dots\xi_n} \tto 0.\\
\end{array}
\end{equation}

\ssec{Unbelievable isomorphisms}\label{UnbIso} The Cartan prolong
 \[
 \fh_{\lambda}(2|2):=(\Pi(\Vol^{\lambda}(0|2)), \ \fosp(2|2)\cong\fvect(0|2))_*
 \]
is a~ deform of $\fh(2|2):=\fh_{1/2}(2|2)$; let
$\fh_{\lambda}(2|2; 1)$ be a~regrading of $\fh_{\lambda}(2|2)$ generated by polynomials in  indeterminates $p, q$ (even) and  $\xi, \eta$ (odd) such that $\deg p=\deg q=1$ and $\deg\eta=0$,\ $\deg \xi=2$. 

It is clear that the Lie superalgebras 
$\fb_{\lambda}(n)$ are non-isomorphic for distinct $\lambda$'s, bar
the occasional isomorphisms 
\begin{equation}
\label{iso}
\text{$\fb_{1/2}(2; 2)\cong \fh_{1/2}(2|2)=\fh(2|2)$ and
$\fh_{\lambda}(2|2)\cong \fh_{-1-\lambda}(2|2)$.}
\end{equation}

\ssec{Deformations of $\fb(n)$}\label{DefBn} Clearly, at the exceptional values of $\lambda$, i.e., 0, 1, and
$\infty$, the deformations of $\fb_{\lambda}(n)$ should be
investigated extra carefully. As we will immediately see, this is rewarded: at each of the exceptional points there are extra deformations. 
The exceptional values $\lambda =\nfrac12$ and $-\frac32$ come from the
isomorphisms \eqref{iso}, see \cite{LSh3}. The origin of an extra deformation at $\lambda =-1$ is unknown.

\sssbegin[Deformations of $\fb(n)$]{Theorem}[Deformations of $\fb(n)$, after \cite{Ko5,LSh1}; for proof, see \cite{LSh3}]\label{thdefb} 
For $\fg=\fb_{\lambda}(n)$, the space $H^2(\fg; \fg)$ is diagonalizable with respect to the Cartan
subalgebra of~$\fder ~\fg$. Let the cocycle $M$ corresponding to the
main deformation be one of the eigenvectors. Let $C$ be another
eigenvector in $H^2(\fg; \fg)$, it determines a~singular deformation from
item~$2)$. The only cocycles $kM+lC$, where $k, l\in\Cee$, that
can be extended to a~global deformation are those for $kl=0$, i.e.,
either $M$ or~$C$.

All the singular deformations of the bracket $\{-, -\}_{old}$ in $\fb_{\lambda}(n)$, except the one at
$\lambda=\nfrac12$ or $-\frac32$ and $n=2$, are linear in $\hbar$ even if $\hbar$ is even:
\begin{equation}
\label{cocycle} \{f, g\}_{\hbar}^{\textrm{sing}}=\{f, g\}_{\text{old}}+\hbar\cdot
C(f, g)\ \ \text{ for any } f, g\in \fb_{\lambda}(n).
\end{equation}

\textbf{Dimension of $H^2(\fg; \fg)$}. \textup{$1_1$)} For $n>2$, we have $\sdim ~H^2(\fg; \fg)=(1|0)$ for $\fg=\fb_{\lambda}(n)$ unless
$\lambda=0$, $-1$, $1$,~$\infty$. 

\textup{$1_2$)} For $n=2$, in addition
to the above, $\sdim ~H^2(\fg; \fg)\neq (1|0)$ at $\lambda=\nfrac12$ and
$\lambda=-\frac32$.

\textup{2)} At the exceptional values of $\lambda$ listed in item
$1)$, we have

$\sdim ~H^2(\fg; \fg)=(2|0)$ at $\lambda=\pm 1$ and $n>2$ odd, or
$\lambda=\infty$ and $n>2$ even, or $n=2$ and $\lambda=\nfrac12$ or
$-\frac32$;

$\sdim ~H^2(\fg; \fg)=(1|1)$ at $\lambda=0$ or $\lambda=\infty$ and $n>2$ odd,
or $\lambda=\pm 1$ and $n>2$ even.

\textbf{Cocycles}. Let $k=(k_{1}, \dots , k_{n})$; we set
$q^k:=q_{1}^{k_{1}}\cdots q_{n}^{k_{n}}$ and $|k|:=\sum k_{i}$. The corresponding cocycles $C$ are given by the following
nonzero values in terms of the generating functions $f$ and $g$,
where $d_{\od}(f)$ is the degree of $f$ with respect to odd
indeterminates only:
\[
\footnotesize
\begin{tabular}{|c|c|c|}
 \hline
$\fb_{\lambda}(n)$&$p(C)$&$C(f, g)$\cr \hline \hline
$\fb_{0}(n)$&\textit{odd}&$(-1)^{p(f)}(d_{\od}(f)-1)(d_{\od}(g)-1)fg$\cr
\hline \hline $\fb_{-1}(n)$&$n+1\pmod 2$&$f=q^k, \ g=q^l \longmapsto
$\cr &&$(4-|k|-|l|)q^{k+l}\xi_{1}\cdots \xi_{n}+\tau
\cS(q^{k+l}\xi_{1}\cdots \xi_{n})$\cr \hline
$\fb_{1}(n)$&$n+1\pmod 2$&$f=\xi_{1}\cdots \xi_{n}, \ g\longmapsto
\begin{cases}(d_{\od}(g)-1)g&\text{if $g\neq af$, $a\in\Cee^{\times}$}\\
 2(n-1)f &\text{if $g=f$ and $n$ is even}\end{cases}$
 \cr
\hline $\fb_{\infty}(n)$&$n\pmod 2$&$f=\tau\xi_{1}\cdots \xi_{n}, \
g\longmapsto
\begin{cases}(d_{\od}(g)-1)g&\text{if $g\neq af$, $a\in\Cee^{\times}$}\\
2f &\text{if $g=f$ and $n$ is odd}\end{cases}$\cr \hline \hline
$\fb_{1/2}(2)\cong \fh_{1/2}(2|2)$&\textit{even}&\textit{induced by the quantization of }$\fpo(2|2)$\cr 
\hline
$\lambda$ generic& even&\textit{see Eqs.} ~\eqref{cocycle}, ~\eqref{bbdef}\cr 
\hline
\end{tabular}
\]

On $\fb_{1/2}(2)\simeq\fb_{-3/2}(2)\simeq\fh(2|2)$, the
cocycle $C$ is the one induced on $\fh(2|2)=\fpo(2|2)/\Cee\, K_1$
by the quantization of the Poisson superalgebra $\fpo(2|2)\colon$ we first
quantize $\fpo(2|2)$, 
and then take the quotient modulo the center
(generated by constants).
\end{Theorem}

\textbf{A nice formula}. Since the elements of $\fb_{\lambda}(n)$ are encoded by functions
in $\tau$, $q$ and
$\xi$ subject to one relation~ \eqref{2.7.4} with an odd left-hand side  which contains
$\tau$, it seems plausible that the bracket in
$\fb_{\lambda}(n)$ can be, at least for generic values of parameter
$\lambda$, expressed solely in terms of $q$ and~$\xi$. This is
indeed the case, at least for the polynomial functions, and here is an explicit formula in which the degree $\deg$ of the polynomial is computed with respect to the standard grading $\deg
q_{i}=\deg \xi_{i}=1$:
\begin{equation}
\label{bbdef}
\renewcommand{\arraystretch}{1.4}
\begin{array}{l}
 \{f_1, f_2\}_{\lambda}^{main}=\{f_1,
f_2\}_{B.b.}+\lambda\left(c_{\lambda}(f_1, f_2)f_1\cS(f_2) +
(-1)^{p(f_1)}c_{\lambda}(f_2, f_1)\cS(f_1)f_2\right),\\
\text{where $ \displaystyle c_{\lambda}(f_1, f_2)=\nfrac{\deg
f_1-2}{2+\lambda(\deg f_2 -n)} $ and 
$\cS$ is the odd Laplacian of eq.~\eqref{DelS}}.
\end{array}
\end{equation}

\sssbegin[The deformation of $\fb_\lambda(n)$ is flat]{Remark}[The deformation of $\fb_\lambda(n)$ is flat] Note that the deformation \eqref{fb} of $\fb_{\lambda}(n; n)$, as well as of its regrading $\fb_{\lambda}(n)$, is flat, i.e., the dimensions of the homogeneous components do not depend on~$\lambda$. Indeed, consider the 0th component $\fm(n)_{0}$ of $\fm(n)$ in its standard $\Zee$-grading. This $\fm(n)_{0}$ is isomorphic to $\fspe(n)\oplus \Cee^2$, i.e., the direct sum of the 2-dimensional center with $\fspe(n):=\fpe(n)\cap\fsl(n|n)$, where $\fpe(n):=\faut(\omega_\od)$, see \cite{Lo}. Note that $\fm(n)_{-1}$ is an irreducible $\fspe(n)$-module. Adding a~line spanned by one of the two central elements of $\fm(n)_{0}$ to $\fspe(n)$ either does not affect the Cartan prolongation with the new 0th component or (and this is our case) adds to the Cartan prolongation the space isomorphic (up to $\Pi$) to the space of functions in basis elements of $\fm(n)_{-1}$ 
regardless of the parameter~$\lambda$. The same independence of the dimensions of homogeneous components on $\lambda$ remains in any other grading of $\fb_\lambda(n)$.
\end{Remark}


\section{Definitions of (almost) hyper-K\"ahler (super)manifolds} Given three (almost) complex
(or queer) structures $J_l$, where $l\in \{i,j,k\}$, satisfying the relations of quaternion units
\begin{equation}
\label{quatunit} J_iJ_j=J_k\ \ \text{for any even
permutation $(i,j,k)$ of $(1,2,3)$},
\end{equation}
and one ``metric" $h$ pseudo-hermitian relative to each $J_l$, together with three bivector fields $B_l$ connected  by three relations of the
form \eqref{JhB}, we arrive at the following new definitions of (almost)
\textit{hyper-K\"ahler (super)manifolds}.

\ssec{Hyper-K\"ahler (super)manifolds. Three even $J$s}\label{3evenJ} In 1988, Verbitsky observed the $\mathfrak{sp}(4)$-action on $\Omega$ on the hyper-K\"ahler manifold $M$, see \cite{Ver}. The passage from manifolds, see \cite{Ver, P1}, to supermanifolds is direct and immediate if the forms $\omega_l$ are even, i.e., the metric $h$ pseudo-hermitian relative each $J_l$, see eq.~\eqref{quatunit}, is even. On any (almost) hyper-K\"ahler manifold, to each of the complex structures $\textbf{i}$,~ $\textbf{j}$ and~ $\textbf{k}$, there correspond the operators on $A^{\bcdot, \bcdot}(M)$, called $E_i$, $E_j$ and $E_k$, and their dual (with respect to the Hodge inner product corresponding to one of these forms, say $\omega_i$) operators $F_i$, $F_j$ and $F_k$ are naturally defined. One also naturally defines the differentials, $d_c^l$ as in eq.~\eqref{d_c}, and $(d_c^l)^*$, where $l\in \{i,j,k\}$ correspond to $\textbf{i}$,~ $\textbf{j}$ and~ $\textbf{k}$, see eq.~\eqref{quatunit}. Together with $d$ and $d^*$, all these operators generate the semi-direct sum $\fg:=\fhei(0|8)\inplus\fsp(4)$. On supermanifolds, the action of this semi-direct sum naturally extends on the spaces of pseudo-differential and pseudo-integral forms.

\ssec{From the $\fsp(4)$-action on the space of differential forms on a~hyper-K\"ahler (super)manifold to an $\fosp(1|4)$-action on the space of twisted and pseudo-differential and pseudo-integral forms}\label{superVerb} Similarly to Bernstein's construction, see Subsec.~\ref{JB}, the $\mathfrak{sp}(4)$-action can be extended to the $\mathfrak{osp}(1|4)$-action on the space $\Omega_\hbar(M)$ of twisted forms, where $\hbar=(\hbar_1, \hbar_2)\in(\Cee^\times)^2$ is now a~pair of factors of curvatures of two connections on the same line bundle over~$M$. 

More explicitly, let $E_i$, $E_j$ and $E_k$ be operators of left multiplication by the (almost) symplectic forms $\omega_i$, $\omega_j$, and $\omega_k$, see Subsec.~\ref{3evenJ}. Let $F_i$, $F_j$ and $F_k$ be their duals with respect to the Hodge inner product corresponding to one of these forms, say~$\omega_i$. Verbitsky showed that two pairs of these dual operators can be considered as Chevalley generators of $\fsp(4)$, call them $X_s^{\pm}$ for~$s=1,2$. Now, let $D_s^{+}$ for $s=1,2$ be ``square roots" of $X_s^{+}$, and let $D_s^{-}$ be defined similarly to formulas~ \eqref{Dpm}: 
\be\label{sqrtX}
\renewcommand{\arraystretch}{1.4}
\begin{array}{l}
D_1^{+}:=d+\nfrac{\hbar_1}{2}\alpha_1 \ \text{such that $d\alpha_1=\omega_i$},\\ 
D_2^{+}:=d+\nfrac{\hbar_2}{2}\alpha_2 \ \text{such that $d\alpha_2=\omega_j$},\\
D_s^{-}:=[X_s^-, D_s^{+}]\ \ \text{for~$s=1,2$. }
\end{array}
\ee
The definition \eqref{sqrtX} implies that the operators $D_s^{\pm}$ for $s=1,2$ generate $\mathfrak{osp}(1|4)$.

The construction of $\mathfrak{osp}(1|4)$-action on the space $\Omega_\hbar(\cM)$ is literally applicable to supermanifolds $\cM$ with an even metric~$h$.

\sssbegin[One twisted differential and the Verbitsky algebra close to $\fosp(2|4)$]{Theorem}[One twisted differential and the Verbitsky algebra close to $\fosp(2|4)$]\label{thtwisthk} Let $D^+=d+\hbar_1\alpha$ be the twisted differential of eq.~\eqref{Dtwist}. In the flat model of Subsec.~$\ref{flat}$, for three even $J$s and an even $h$, the closure of the algebra generated by $\{D^+, E_l, F_l\}_{l=i,j,k}$ --- \textit{one} twisted differential and the Verbitsky algebra --- is a~\textbf{simple} Lie superalgebra of $\sdim=(11|8)$, whose even part is $\fsp(4)\oplus\Cee z$ with $z$ central in it, and on whose odd part $\ad_z$ has two eigenvalues $\pm\lambda$, $\lambda\neq0$, each of multiplicity~$4$. Hence, it is~$\fosp(2|4)$.
\end{Theorem}

\begin{proof} The closure is computed in the flat model of Subsec.~\ref{flat}, and the simplicity of the answer and the centrality of $z$ in its even part are verified as in the proof of Theorem~\ref{thtwist}. The identification then needs no further computation: by the classification of simple finite-dimensional Lie superalgebras, see a~review \cite{Lo}, the one with even part $\fsp(4)\oplus\Cee$ and $8$-dimensional odd part is $\fosp(2|4)$, in which the $\pm\lambda$-eigenspaces of the $\fo(2)$-generator $z$ are the two copies of the tautological $\fsp(4)$-module $\Cee^{4}$.\end{proof}

\parbegin[Why $\fsp(2)$ and $\fsp(4)$, and not $\fsl(2)$ and $\fo(5)$]{Remark}[Why $\fsp(2)$ and $\fsp(4)$, and not $\fsl(2)$ and $\fo(5)$]\label{whysp4} Observe that it is more appropriate to speak about $\fsp(2)$ not $\fsl(2)$ (resp., $\fsp(4)$ not $\fo(5)$) acting in the space of differential forms on the almost K\"ahler (resp., hyper-K\"ahler)
(super)manifold $M$: superization of Weil's (resp., Verbitsky's) result to the space of differential forms with values in the line bundle with connections corresponding to non-integral distributions given by the contact forms 
whose curvatures are proportional to the almost symplectic forms participating in the definition of the almost K\"ahler (resp., hyper-K\"ahler) manifold describes an $\fosp(1|2)$-action on the space $\Omega_\hbar(\cM)$ 
(resp., an $\fosp(1|4)$-action in the hyper-K\"ahler case), see \cite{LeLh}.
 Actually, Verbitsky's result \cite{Ver} described over $\Ree$ is more precise than its exposition over $\Cee$ given here: Verbitsky described actions of the real forms of the complex Lie algebra $\fo(5)\simeq\fsp(4)$$;$ the real forms of its superizations are listed in \cite{LeLh}. 
 \end{Remark} 

\sssbegin[Poletaeva's contraction, and what the twist does to it]{Remark}[Poletaeva's contraction, and what the twist does to it]\label{poletcontr} The two algebras that the twist produces are not new to the literature, see \cite{P1}. Its Conclusion records, of the two Lie superalgebras that paper is about --- the $\fhei(0|4)\inplus\fsl(2)$ of eq.~\eqref{ADolb} and the ${\fhei(0|8)\inplus\fsp(4)}$ of Subsec.~\ref{3evenJ} --- that they \lq\lq are contractions of $\fosp(2|2)$ and~$\fosp(2|4)$\rq\rq. There the contraction is named, and not constructed. Our Theorems~\ref{thtwist} and~\ref{thtwisthk} construct the deformation whose limit it is, and construct it by operators$:$ the parameter is the pair $(\hbar_1, \hbar_2)$ of the factors by which the two curvature forms are proportional to $\omega$, and the limit is the untwisted case $\hbar_1=\hbar_2=0$.

What makes this a~contraction, and not a~coincidence of superdimensions, is that the superdimension is constant along the family --- a~contraction being a~limit of brackets on a~\textit{fixed} space. The family is the twisted one itself, taken on the diagonal $\hbar_1=\hbar_2=\hbar$, so that $D^+=d+\hbar\alpha$ and $D^+_c=d_c+\hbar\alpha_c$ are the twisted differentials of eq.~\eqref{Dtwist}, with $\hbar$ ranging over the whole of $\Cee$ --- the restriction $\hbar_1\hbar_2\neq0$ of that equation is lifted here, the point of the family being precisely the value~$\hbar=0$. In the flat model, the closure of the algebra generated by $\{D^+, D^+_c, X^+, X^-\}$ is of $\sdim=(4|4)$ for every~$\hbar$, simple for $\hbar\neq0$ and, for $\hbar=0$, the algebra of \cite{P1}, whose even center is spanned by the Laplacian of eq.~\eqref{Del}; likewise, the closure of the algebra generated by $\{D^+, E_l, F_l\}_{l=i, j, k}$ is of $\sdim=(11|8)$ for every $\hbar$, simple for $\hbar\neq0$ and, for $\hbar=0$, the other algebra of \cite{P1}, of its own~$\sdim=(11|8)$. Moreover, in the basis $\{X^+, X^-, H, \Delta_{\hbar, \hbar}, D^{\pm}, D^{\pm}_c\}$, where $\Delta_{\hbar, \hbar}$ is the twisted Laplacian of eq.~\eqref{Deltw}, every structure constant is a~polynomial in $\hbar$ of degree~$1$, and its value at $\hbar=0$ is the corresponding structure constant of the untwisted algebra. So the family is a~flat degeneration and the algebras of \cite{P1} are its special fiber.
\end{Remark}

\ssec{An almost periplectic structure and two odd $J$s}\label{ahk2eJ} Since one of the $J_l$ is the product of the two others, see eq.~\eqref{quatunit}, the number of odd $J$s among $J_{\textbf{i}}$, $J_{\textbf{j}}$, $J_{\textbf{k}}$ is either $0$ or $2$, because eq.~\eqref{quatunit} gives $p(J_{\textbf{k}})=p(J_{\textbf{i}})+p(J_{\textbf{j}})$. This still leaves two possibilities for the pair $(J_{\textbf{i}}, J_{\textbf{j}})$ we start from --- both odd, or one even and one odd --- since each gives exactly two odd $J$s among the three. The second of them does not occur:

\sssbegin[The three $J$s cannot be of mixed parities]{Proposition}[The three $J$s cannot be of mixed parities]\label{nomixed} There is no almost hyper-K\"ahler superstructure with $J_i$ even and $J_j$ odd, for either parity of $h$ and in any dimension, finite or infinite.
\end{Proposition}

\begin{proof} Everything is pointwise, and neither a~basis nor an inverse is needed --- which is why no dimension enters.

Since $p(\omega_l)=p(h)+p(J_l)$, see eq.~\eqref{eq2}, and $p(J_k)=p(J_i)+p(J_j)$, exactly one of the three forms $\omega_l=h(J_l\ -, -)$ is even for $h$ even, and exactly one of them is odd for $h$ odd. Let $\omega$ be an even one and $\omega'$ an odd one. In either case, the two are related by the \textit{odd} operator $J_k$, and by eq.~\eqref{quatunit} alone:
\be\label{omomprime}
\omega'=\varepsilon\,\omega(J_k\ -, -),\quad\text{where }\
\varepsilon=\begin{cases}-1&\text{for $h$ even, when $\omega=\omega_i$, $\omega'=\omega_j$},\\
1&\text{for $h$ odd, when $\omega=\omega_j$, $\omega'=\omega_i$},\end{cases}
\ee
because eq.~\eqref{quatunit} gives $J_iJ_k=-J_j$ and $J_jJ_k=J_i$, so both sides of eq.~\eqref{omomprime} are $h(J\ -, -)$ for one and the same~$J$. The sign $\varepsilon$ plays no role in what follows.

Let now $X$, $Y\in V_\ev$ and set $b(X, Y):=\omega(J_kX, J_kY)$. Since $J_k$ is odd, both $J_kX$ and $J_kY$ lie in $V_\od$, where the even $2$-form $\omega$ is \textbf{symmetric}; hence, $b$ is symmetric. On the other hand, eq.~\eqref{omomprime} evaluated at $(X, J_kY)$ and at $(J_kY, X)$ gives
\[
\varepsilon\,b(X, Y)=\omega'(X, J_kY)=-\omega'(J_kY, X)=-\varepsilon\,\omega(J_k^2Y, X)=\varepsilon\,\omega(Y, X)=-\varepsilon\,\omega(X, Y),
\]
where the second equality holds because $\omega'$ is an odd $2$-form while $X$ is even, the fourth one because $J_k^2=-\id$, and the last one because $\omega$ is anti-symmetric on~$V_\ev$. Hence, $b=-\omega\vert_{V_\ev}$ is anti-symmetric as well, so $b=0$ and $\omega\vert_{V_\ev}=0$, whereas $V_\ev\neq 0$ because the odd invertible $J_k$ identifies $V_\ev$ with~$V_\od$. So the even $2$-form $\omega$ is degenerate, contrary to the hypothesis.

In adapted bases, the clash is the one of Lemma~\ref{p2a}: it reads $\Omega=P^{T}GP$ with $\Omega$ anti-symmetric and $G$ symmetric, and a~congruence preserves the symmetry.\end{proof}

\sssec{Two odd $J$s and an odd pseudo-hermitean form $h$}\label{oddh} Having defined the differentials $d_c^l$, where $l \in\{1, \textbf{i}, \textbf{j}, \textbf{k}\}$, cf.~ \eqref{quatunit}, we  arrive at a new version of an (almost) \textit{hyper-K\"ahler supermanifold}. Here, an odd symmetric $h$ together with an odd $J$ such that $J^2$ is a~non-zero scalar does exist for every even $n$ on $\Cee^{n|n}$, but the tensor $\beta=h(J\ -, -)$ is then never a~non-degenerate $2$-form, and so never an~$\omega$: by Lemma~\ref{p2a} its two diagonal blocks are congruent to one another, hence of one and the same symmetry, so whichever of them is the anti-symmetric one the other has the wrong, metric-type, symmetry. Following the change of paradigm of Subsec.~\ref{omegaB}, in which the leading role passes to the bivector $B$, the operators to consider are therefore the multiplications by the $2$-form \textit{avatars} of the $\omega_l$ and the convolutions with the bivectors $B_l=h^*(J_l\ -, -)$. With the standard odd pairing $h$ on $\Cee^{2|2}$, and writing $y_1, y_2$ for $\widehat q_1, \widehat q_2$ and $y_3, y_4$ for $\widehat \pi_1, \widehat \pi_2$, the two odd $J$s of eq.~\eqref{quatunit}, taken in the adjoint representation of the quaternion superalgebra, give exactly 
\be\label{Eodd}
\renewcommand{\arraystretch}{1.4}
\begin{array}{l}
E_i=-\nfrac12\left(y_3^2+y_4^2\right),\ \ E_j=-y_1y_2\ \ \text{(both even, each supported on one parity block)},\\
E_k=y_1y_4-y_2y_3
\end{array}
\ee
(the latter is odd), and dually for the~$F_l$. Observe that $\omega_k$ is here an honest non-degenerate \textbf{odd} $2$-form, so that the roles are inverted as compared with Subsec.~\ref{Sost}.

\sssbegin[An odd $h$: $\fhei(0|8)\inplus\fsp(4)$ is replaced with $\fn\inplus\fsl(2|2)$]{Theorem}[An odd $h$: $\fhei(0|8)\inplus\fsp(4)$ is replaced with $\fn\inplus\fsl(2|2)$]\label{thoddh} In the flat model of Subsec.~$\ref{flat}$: 

$1)$ The closure of the algebra generated by $\{E_l, F_l\}_{l=i,j,k}$ is a~perfect Lie superalgebra with $1$-dimensional center, of $\sdim=(7|8)$, whose even part is $\fsl(2)\oplus\fsl(2)\oplus\Cee$ and whose odd part is the sum of two copies of the bimodule $\Cee^2\otimes\Cee^2$, of $\sdim=(0|4)$ each; that is, it is~$\fsl(2|2)$. This closure automatically contains the multiplication operator $E_h$ by the $2$-form avatar $y_1y_3+y_2y_4$ of the pseudo-Hermitian form $h$ itself; here $h$ is the (pseudo-)Hermitian form of the almost hyper-K\"ahler structure, the one entering $\omega_l=h(J_l\ -\,,-)$; here $h$ is odd, this being assumed throughout the present Subsection, and $F_h$ is the dual of~$E_h$.

$2)$ Let $t_c:=\partial_{z_c}$ and $\partial_{\widehat y_c}$ be two of the generators of the super Weyl algebra of eq.~\eqref{engine}: $\widehat y_c$ is the operator of left multiplication by a~differential,
$\partial_{\widehat y_c}$ the convolution with it, $t_c=\partial_{z_c}$ the derivation along a~coordinate. So a~differential is of type $\widehat y\,\partial_z$ and a~codifferential of type $\partial_z\,\partial_{\widehat y}$; then, the closure of the algebra generated by $\{E_l, F_l, d\}_{l=i,j,k}$ is of $\sdim=(17|18)$ and is equal to $\fn\inplus\fsl(2|2)$, where $\fn$ is a~Heisenberg-type ideal of $\sdim=(10|10)$: its part $\fn_{-1}$ of $\sdim=(8|8)$ --- eight differentials of $\widehat y\,t$-type and eight codifferentials of $t\,\partial_{\widehat y}$-type 
--- pairs into the $(2|2)$-dimensional center spanned by two even and two odd Laplacians.
\end{Theorem}

Thus, the answer here is \textbf{not} $\fhei(0|8)\inplus\fsp(4)$, as it was for three even $J$s in Subsec.~\ref{3evenJ}: several of the invariant operators (namely, $E_k$, $F_k$, $E_h$, $F_h$, $d$, \dots) are odd, and, accordingly, $\fsp(4)$ changes into $\fsl(2|2)$ and $\fhei(0|8)$ into the $(10|10)$-dimensional~$\fn$. It is amusing that the same input --- the quaternion superalgebra with two odd units --- produces an algebra \textit{larger} than the one it superizes; it is natural to call the procedure \textit{enlarging}:
\be\label{grows}
\renewcommand{\arraystretch}{1.4}
\begin{array}{l}
\dim\fsl(2|2)=15\ \ \text{against}\ \ \dim\fsp(4)=10,\\
\sdim\big(\fn\inplus\fsl(2|2)\big)=(17|18)\ \ \text{against}\ \ \sdim\big(\fhei(0|8)\inplus\fsp(4)\big)=(11|8).
\end{array}
\ee
For an even $h$, the answer \textit{shrinks}, see Corollary~\ref{corSost}, but the term \lq\lq contraction\rq\rq\ is not applicable there either: the dimensions of the contracted algebra and the initial one are equal. 

\subsubsection{Two odd $J$s and an even pseudo-hermitean form $h$}\label{Sost}
The question is whether the Lie superalgebra of invariant operators realizes one of these two superizations. It does not, and here is what happens instead.

Here, $\omega_i$ and $\omega_j$ are odd non-degenerate $2$-forms, and on $\Cee^{2|2}$ the compatibility conditions force, up to equivalence, the pair of odd forms in the periplectic Darboux normal shape 
\be\label{DarbouxPair}
\omega_i=y_1y_3+y_2y_4,\ \ \omega_j=y_1y_4-y_2y_3,\ \ \text{where}\ \ T_{ij}:=B_i\circ \omega_j=J_{\textbf{k}},\ \ T_{ij}^2=-\id,
\ee
in the notation of eq.~\eqref{Eodd}, and $J_{\textbf{k}}$ is exactly the even matrix $J_{\textbf{i}}J_{\textbf{j}}$ of Subsec.~\ref{ahk2eJ}. 

There is no third operator to go with $E_i$ and $E_j$, and here is why. The third quaternion unit $J_{\textbf{k}}$ is even, so the tensor $\omega_k:=h(J_{\textbf{k}}\ -, -)$ is even as well. Since $\omega_l=h(J_l\ -, -)$ for each $l$, and $J_{\textbf{k}}=J_{\textbf{i}}J_{\textbf{j}}$, this tensor is a~composition of three maps already at hand,
\be\label{omk}
\omega_k=\omega_i\circ B_h\circ \omega_j,\ \ \text{where}\ \ B_h:=h^{-1}\ \ \text{is the bivector dual to $h$},
\ee
in which no operator is transposed. On $\Cee^{2|2}$, with $h$ normalized together with the pair~\eqref{DarbouxPair}, i.e., with $h$ equal to the identity on $V_\ev$ and to $\left(\begin{smallmatrix}0&1\\-1&0\end{smallmatrix}\right)$ on $V_\od$, eq.~\eqref{omk} gives
\be\label{omkflat}
\omega_k=\Big(-\nfrac14\, \id\Big)_{V_\ev}\ \oplus\ \nfrac14\left(\begin{smallmatrix}0&1\\-1&0\end{smallmatrix}\right)_{V_\od},
\ee
which is symmetric on $V_\ev$ and anti-symmetric on~$V_\od$. By the criterion recalled just after eq.~\eqref{OmDiag}, this is the symmetry of an even \textbf{symmetric} tensor --- a~second metric --- and not of an even $2$-form, which requires the two symmetries the other way round. A~tensor of that symmetry has no $2$-form avatar at all: in the notation of eq.~\eqref{engine}, $\sum (\omega_k)_{ce}\widehat y_c\widehat y_e=0$ identically, every term canceling against its partner. Hence, the operators $E_k$ and $F_k$ do \textit{not} exist, and the generating set consists of the four \textit{odd} operators $E_i$, $E_j$, $F_i$, $F_j$ together with~$d$.

The contrast with Subsec.~\ref{oddh} is instructive, and it is a~check on eq.~\eqref{omk}: there $h$ is odd, hence $\omega_k$ is odd, and the very same composition returns $-\nfrac12(y_1y_4-y_2y_3)$, i.e., a~multiple of the $E_k$ printed in eq.~\eqref{Eodd} --- which is why the third pair of operators exists there and not here. 

\sssbegin[An even $h$: the algebra generated is solvable, of $\sdim=(6|10)$]{Theorem}[An even $h$: the algebra generated is solvable, of $\sdim=(6|10)$]\label{thevenh} In the flat model of Subsec.~$\ref{flat}$: 

$1)$ The Lie superalgebra generated by $\{E_i, E_j, F_i, F_j\}$, 
see Subsec.~$\ref{flat}$, is of $\sdim=(2|4)${}$:$ the elements
\[
H_E:=[E_i, F_i]=[E_j, F_j]\ \ \text{and}\ \ K_E:=[E_i, F_j]=-[E_j, F_i]
\]
are both central, all the brackets $[E_a, E_b]$ and $[F_a, F_b]$ vanish, and over $\Cee$ it is 
\[
\fsl(1|1)\oplus\fsl(1|1)\simeq \fhei(0|2)\oplus\fhei(0|2)
\]
 $($take $E_\pm:=E_i\pm\sqrt{-1}E_j$ and $F_\pm:=F_i\pm\sqrt{-1}F_j$, then $\Span(E_+, F_-)$ and $\Span(E_-, F_+)$ generate the two commuting copies$)$;

$2)$ The Lie superalgebra $\fg$ generated by $\{E_i, E_j, F_i, F_j, d\}$ is \textbf{solvable}, of $\sdim=(6|10)$, and its derived series
\be\label{derser}
\fg\supset\fg^{(1)}\supset\fg^{(2)}\supset\fg^{(3)}\supset\fg^{(4)}=0,\ \ \text{where}\ \ \fg^{(k+1)}:=[\fg^{(k)}, \fg^{(k)}],
\ee
reaches $0$ in four steps, the superdimensions along it being $(6|10)$, $(6|6)$, $(4|6)$, $(0|2)$, $(0|0)$; namely, $\fg=\fn\inplus\fa$, where $\fa$ is the $(2|4)$-dimensional algebra of heading $1)$ and $\fn$ is a~two-step nilpotent ideal of $\sdim=(4|6)$ spanned by the four odd differentials $d$, $[H_E,d]$, $[K_E,d]$, $[K_E,[H_E,d]]$, by four even codifferentials, and by the $(0|2)$-dimensional center $\Span( L_i, L_j)$ of $\fn$, spanned by the two odd Laplacians, so that $[\fn_{-1}, \fn_{-1}]=\Span(L_i, L_j)$ is a~$\Pi(\Cee^2)$-valued pairing. 
\end{Theorem}

\sssbegin[For an even $h$, the Verbitsky algebra collapses]{Corollary}[For an even $h$, the Verbitsky algebra collapses]\label{corSost} For two odd $J$s and an even $h$, the Verbitsky algebra $\fsp(4)$ \textbf{collapses} to the $(2|4)$-dimensional algebra of Theorem~$\ref{thevenh}$, item~$1)$ --- the whole semi-simple part disappears --- and $\fhei(0|8)\inplus\fsp(4)$, of $\sdim=(11|8)$, is replaced with the solvable $(6|10)$-dimensional algebra of Theorem~$\ref{thevenh}$, item~$2)$. This is the answer to item~$1)$ of Section~\ref{superPolet}, where the question is posed. \end{Corollary}

\be\label{HKsum}
\footnotesize
\renewcommand{\arraystretch}{1.3}
\begin{tabular}{|c|c|}
 \hline
$\omega$ or $B$ (and $\vvol$)&$\fg$ acting on $\Omega$ and $\Sigma$, also $\widehat\Omega$ and $\widehat\Sigma$\cr 
\hline \hline
$B=\omega_0^{-1}$&$\fhei(0|2)\inplus\fsl(2)$, see  eq.~\eqref{a=0}\cr
\hline \hline 
$B$ and $\omega_1$&$\fhei(0|2)\simeq\fsl(1|1)$, see \cite{LeLh}\cr
\hline
$\vvol$, $B=\omega_1^{-1}$&$\fb\fa(2)\inplus\fsl(1|1)$, of $\sdim(3|5)$, see  eq.~\eqref{h} and Lemma \ref{gl11}
 \cr
\hline 
\hline
\hline
$J$, $h^*$, $B$, see eq.~\eqref{JhB}&$\fg$ acting on $A^{\bcdot,\bcdot}$, see Subsec.~\ref{p1}, and $B^{\bcdot,\bcdot}$, see Subsec.~\ref{p2}
 \cr
 \hline
 \hline
$p(J)=p(h^*)=\ev$, see Subsec.~\ref{p1}&$\fhei(0|4)\inplus\fsl(2)$, of $\sdim(4|4)$, see  eq.~\eqref{ADolb}
\cr
$p(J)=p(h^*)=\od$, see Subsec.~\ref{p2}&$\fsl(2)\inplus\Pi(\Cee^2\oplus\Cee)$, of $\sdim(3|3)$, see Statement \ref{p2b}
\cr
$p(J)=p(h^*)+\od$, see Subsec.~\ref{p3}&$\fb\fa(2)\inplus\fsl(1|1)$, of $\sdim(3|5)$, see  eq.~\eqref{ADolbOdd}\cr
\hline
\hline
$\vvol$, $B$, $\omega$&$\fg$ acting on $\Omega_\hbar$ and $\Sigma_\hbar$\cr 
\hline 
$B=\omega_0^{-1}$&$\fosp(1|2)\oplus \Cee H$, see  eqs.~\eqref{H} and \eqref{Dpm}\cr
\hline
$\vvol$, $B=\omega_1^{-1}$, see Subsec.~\ref{fgMyst}&$\fg$ of $\sdim(8|8)$, see Theorem \ref{thP7} and eq.~\eqref{Dpm2}\cr
\hline
$J$, $h^*$, $B$, see eq.~\eqref{JhB}&$\fg$ acting on $A^{\bcdot,\bcdot}_{\hbar_1,\hbar_2}$ and on $B^{\bcdot,\bcdot}_{\hbar_1,\hbar_2}$, see Subsec.~\ref{ssTwDol}\cr
\hline 
$B=\omega_0^{-1}$&$\fosp(2|2)\ni\Delta_{\hbar_1, \hbar_2}$, see Theorem~\ref{thtwist} and eq.~\eqref{Deltw}\cr
\hline\hline
\end{tabular}
\ee

\section{Superizations of (almost) hyper-K\"ahler structures: summary, see eq.~\eqref{HKsum}} \label{sDefSMan}
(i) The main change of paradigm: the role of non-degenerate differential 2-form $\omega$ (be it even or odd) is given to the bivector~$B$. If $\cM$ is periplectic, then $B=\omega^{-1}$ only if~ $\cM$ is special, i.e., with a~preserved volume element. This change results in  superizations of several notions:

(ii)  Of symplectic manifold, its periplectic version with the odd 2-form $\omega$, and a~particularly interesting special periplectic version with a~volume element preserved. Related are interesting results due to Kochetkov (classification of irreducible continuous $\fle(n)$- and $\fsle(n)$-modules), reviewed in \cite{LeLh}. 

(iii) In addition to the discrete parameters (superdimensions \eqref{4types} of $\cM$), the supervariety $\cQ$ of deformed versions of the Lie superalgebra defined by $B$ enters the scenery, see Section~\ref{S2}.

(iv) The odd versions of the complex structure $J$, the ``metric" $h$ and its dual $h^*$, and the 2-form $\omega$  bring new versions of (almost) K\"ahler supermanifolds summarized in the table above together with their even counterparts for comparison.

\textbf{Superizations of almost hyper-K\"ahler manifolds}. 

(a) On any (almost) hyper-K\"ahler (super)manifold with three even $J$s, starting from the $\fsp(4)$-action on the space of (pseudo)-differential forms we arrive at an $(\fosp(1|4)\oplus \Cee H)$-action on the space of twisted forms, see eqs.~\eqref{H} and ~\eqref{sqrtX}; adjoining \textit{one} twisted differential to the Verbitsky algebra $\fsp(4)$ instead gives $\fosp(2|4)$, see Theorem~\ref{thtwisthk}. Here, the answer is one and the same for $A^{\bcdot, \bcdot}_{\hbar_1, \hbar_2}$ and for $B^{\bcdot, \bcdot}_{\hbar_1, \hbar_2}$, and it is so already before the twist. Indeed, all three $J$s and $h$ being even, the form $\omega$ is even, identification~\eqref{ident} applies, the factor $\Vol$ is trivialized by a~power of $\omega_\ev$, and hence
\be\label{ABtwistiso}
B^{\bcdot, \bcdot}\simeq A^{\bcdot, \bcdot}\ \ \text{up to a~shift of the bi-degree, and therefore}\ \ B^{\bcdot, \bcdot}_{\hbar_1, \hbar_2}\simeq A^{\bcdot, \bcdot}_{\hbar_1, \hbar_2},
\ee
the second isomorphism being the first one tensored by the \textbf{same} pair $(L, \nabla)$, see eq.~\eqref{Btwist}. This is why the subscripts $\hbar_1, \hbar_2$ carry no information in the even case and are omitted below whenever the statement holds for the twisted and the untwisted spaces alike. The distinction between $A^{\bcdot, \bcdot}$ and $B^{\bcdot, \bcdot}$ becomes essential only for $\omega$ odd, and even then only if $\cM$ is not special, see Remark~\ref{whyspecial} and item~$1)$ of Subsec.~\ref{ssTwDol}.

(b) Two odd complex structures $J$ on a special periplectic supermanifold bring totally new amazing versions of the constructions in item (a), see Subsecs.~\ref{ahk2eJ}--\ref{Sost}. There is no mixed case: one even and one odd $J$ never occur together, see Proposition~\ref{nomixed}. The two cases that do occur are as follows, and in neither of them is the answer the one of item~(a). In the table, each row records \textbf{two} Lie superalgebras, each of them the closure of the algebra generated by a~collection of operators in the sense of Subsec.~\ref{flat}. The first one is the closure of the algebra generated by the operators $E_l$ and $F_l$ \textit{alone}, where $E_l$ is the multiplication by the $2$-form avatar of $\omega_l$ and $F_l$ is the convolution with the bivector $B_l$; it is the analog of Verbitsky's $\fsp(4)$ of Subsec.~\ref{3evenJ}, and it involves no differentiation. The second one, the entry marked \lq\lq with $d$", is the closure of the algebra generated by the \textit{same} operators together with the exterior differential $d$; it is the analog of $\fhei(0|8)\inplus\fsp(4)$. Adjoining $d$ is what brings in the differentials: the brackets of $d$ with the $E_l$ and the $F_l$ produce the twisted differentials $d^l_c$ and the codifferentials dual to them, and the brackets of these among themselves produce the Laplacians, which span the center of the Heisenberg-type ideal. In the classical case, this passage is $\fsp(4)\rightsquigarrow \fhei(0|8)\inplus\fsp(4)$, i.e., $(10|0)\rightsquigarrow(11|8)$, see Remark~\ref{validation}; the two columns of the table show what it turns into for two odd $J$s.

\ssbegin[Validation]{Remark}[Validation]\label{validation} The computation reproduces the two known answers exactly: on the flat hyper-K\"ahler $\Cee^4$ with three even $J$s, the closure of the algebra generated by $\{E_l, F_l, d\}$ is $\fhei(0|8)\inplus\fsp(4)$ of $\sdim=(11|8)$, see Subsec.~\ref{3evenJ}, and the closure of the algebra generated by $\{E, F, d, d_c\}$ in the case of one even $J$ and even $h$ is $\fhei(0|4)\inplus\fsl(2)$ of $\sdim=(4|4)$, see eq.~\eqref{ADolb}. On a~special almost periplectic supermanifold, the closure of the algebra generated by $\{\cX^\pm, d\}$ is $\fb\fa(2)\inplus \fsl(1|1)$ of $\sdim=(3|5)$, see eq.~\eqref{ADolbOdd}. These answers apply verbatim to \textbf{any} K\"ahler supermanifold, curved or not: the operators are of the first order in the structure tensors, and on a~K\"ahler $\cM$ these tensors are parallel, see Remark~\ref{parallel}. On \textit{almost} K\"ahler supermanifolds, where parallelism fails, the tensors that measure the failure --- $d\omega$ and the Nijenhuis tensor --- can deform the $\fhei$-type (for $\omega$ even) or the $\fb\fa$-type (for $\omega$ odd) parts of the answers, but not the reductive quotients. 
\end{Remark}
\[
\renewcommand{\arraystretch}{1.3}
\begin{tabular}{|c|c|c|}
 \hline
$p(h)$&$\fg$ acting on $\widehat\Omega$ and $\widehat\Sigma$&see\cr
\hline\hline
$\ev$&$\fsl(1|1)\oplus\fsl(1|1)$, of $\sdim(2|4)$; with $d$: solvable, of $\sdim(6|10)$&Theorem \ref{thevenh}\cr
\hline
$\od$&$\fsl(2|2)$, of $\sdim(7|8)$; with $d$: $\fn\inplus\fsl(2|2)$, of $\sdim(17|18)$&Theorem \ref{thoddh}\cr
\hline
\end{tabular}
\]
In particular, $\fsp(4)$ is replaced with a~smaller algebra for $h$ even, see Corollary~\ref{corSost}, and, for $h$ odd, with the strictly larger $\fsl(2|2)$ of eq.~\eqref{grows}. 

\section{On history. Open problems and conjectures}\label{superPolet} 

\ssec{On brackets}\label{ssHist} The bracket we denote by $\{-,-\}_{B.b.}$ was discovered by Schouten, and is still called the \textit{Schouten bracket}
in Differential Geometry. Hereafter, $\{-,-\}_{B.b.}$ is called the 
\textit{Buttin bracket}  in honor of C.~Buttin who was the first to prove that $\{-,-\}_{B.b.}$ satisfies the super Jacobi identity, see \cite{Bu}. In
\cite{L1}, the analog of Darboux's theorem on the normal shape of a~non-degenerate odd 2-form on freshly defined ``supermanifolds" (\cite{BerL}) was formulated, the related simple Lie superalgebras $\fle(n)$ and $\fm(n)$, new at that time, were introduced, and  $\{-,-\}_{B.b.}$ was interpreted as an analog of the
Poisson bracket  in mechanics. Batalin and
Vilkovisky rediscovered this interpretation of $\{-,-\}_{B.b.}$ with 
applications to theoretical physics, see \cite{BV1, BV2}; they dubbed $\{-,-\}_{B.b.}$
\textit{anti-bracket}, which is now a~standard term not only in Physics. Zinn-Justin 
 discovered an \textit{antiproduct} $F\ast G$, see \cite[eq.~(55), p.31]{ZJ};  and ca 7 years later it was noted that $\{F,G\}_{B.b.}=F\ast G-(-1)^{(p(F)+1)(p(G)+1)}G\ast F$. 

\sssec{The first hints at $B$ odd} Superizations of K\"ahler and hyper-K\"ahler manifolds  for $B$ odd were first mentioned in \cite{L2}, see also
\cite{BGLS} --- as additional examples for which the super analogs of
the Nijenhuis tensor can be computed. The main objects introduced and studied in \cite{BGLS} were, however, the \textbf{real-complex} (super)manifolds, i.e., the ones with a~non-integrable distribution and the (almost) complex structure given only on the distribution, not on the whole tangent bundle. Related were versions of the Nijenhuis and curvature tensors --- the obstructions to the possibility to get rid of the adjective ``\textbf{almost}" speaking of symplectic and periplectic, as well as complex and real-complex, structures. The definitions I suggest in this article are similar in essence
(the main characters being the bivectors) to the approach of several
groups of researchers to generalized Calabi-Yau manifolds,
generalized Complex Geometry, and supersymmetric sigma-model, see e.g., \cite{KhN, FFKS,HLRUZ} and references therein.

\sssec{The condition $d\omega_\ev=0$ is not one equation but two}\label{ssTwoCond} The two equations are the vanishing of the primitive component of $d\omega_\ev$ and the vanishing of its component divisible by $\omega_\ev$, and on manifolds each of them has been studied on its own$:$ the four classes into which the covariant derivative of $\omega_\ev$ splits are the basis of the classification of almost Hermitian manifolds, see \cite{GH}, and imposing only the vanishing of the primitive half gives the locally conformally symplectic manifolds, in which $d\omega_\ev=\theta\wedge\omega_\ev$ for a~closed $1$-form~$\theta$, whereas imposing only the vanishing of the other half is the balanced, or semi-K\"ahler, condition~$\theta=0$. The splitting itself goes back to \cite{Lee1}, where the $1$-form $\theta$, now called the \textit{Lee form}, first appears. Its author is Lee Hwa-Chung, who isolated the first of the two in the 1940s; the classification of the almost Hermitian manifolds by the components of $\nabla\omega$ grew out of them, see \cite{GH}.\footnote{How Lee's papers, unread for decades, became known in Europe is on record in a~footnote of \cite[p.~698]{EL}$:$ \lq\lq Cet article, qui nous a~\'et\'e signal\'e par M.~Yen Chin-Ta apr\`es la parution de notre Note, contient d\'ej\`a sous une autre forme plusieurs r\'esultats \'enonc\'es par nous\rq\rq.} 

For the description of what this becomes on supermanifolds, see Lemma~\ref{halves}. For $\omega_\ev$, the two halves are independent conditions if and only if $2a-b\neq2$, where $\sdim\cM=(2a|b)$; on the supermanifolds of $\sdim=(2a|2a-2)$ --- the super analogs of the symplectic surface --- every multiple of $\omega_\ev$ is already primitive, and there is only one condition, see Example~\ref{fuse}. For  $\omega_\od$, the two halves exist for every $n$, where $\sdim\cM=(n|n)$, although the operator that separates them is no longer the $H_{sl}$ of a~Lefschetz $\fsl(2)$-action: it is not a~scalar on $\Omega^k$, only on the finer components given by the number of the $\widehat q_i$ and of the $\widehat \xi_j$ in a~monomial --- which is the same phenomenon that replaces $\fsl(2)$ with $\fb\fa(2)\inplus \fsl(1|1)$ in Lemma~\ref{gl11}.

The superization of the first of the two classical notions is Proposition~\ref{lee}, which also settles the odd counterpart of $\theta${}$:$ the Lee form is unique in every superdimension but $(2|0)$ and $(1|1)$, and its differential, for $\omega_\od$, is $\varphi\omega_\od$ for an odd function $\varphi$ --- plus, only on the supermanifolds of $\sdim=(2|2)$, a~multiple of $\widehat q_1\widehat q_2$. That such a~$\varphi$ is always zero --- whatever $n$, the case $n=2$ included --- is Lemma~\ref{nilLee}, the reason being that an odd function is nilpotent and the Euler operator kills the leading term of the relation between $d\varphi$ and $\varphi\theta$; for $n=2$, that relation acquires a~correction, which the parity of the coefficients of $\omega_\od$ confines to $J^2$, where the Euler operator does not see it. So the Lee form of  $\omega_\od$ is closed whenever $n\neq2$, and the analog of Lee's theorem holds verbatim there. 
On the supermanifolds of $\sdim=(2|2)$, however, it fails, see Example~\ref{noint}, which exhibits an odd non-degenerate form dividing its own differential and admitting no integrating factor; there, the whole obstruction is carried by~$\psi$. 

The superdimensions on which Lee's theorem is not applicable  
are therefore $(4|0)$ for $\omega_\ev$ and $(2|2)$ for $\omega_\od$, and 
for opposite reasons$:$ on the first, the divisibility of $d\omega_\ev$ by $\omega_\ev$ demands nothing, every $3$-form being a~multiple of $\omega_\ev$; on the second, it does not entail the closedness of~$\theta$.

\ssec{Open problems and conjectures}\label{ssOpen}

1) In the flat model, define the differentials $d_c^i$ and $d_c^j$ by the formula \eqref{d_c} in which the tensors $J_{\textbf{i}}$ and $J_{\textbf{j}}$ are both odd --- the only case that occurs, see Proposition~\ref{nomixed}; define $d_c^k$ accordingly. Only the two cases distinguished by $p(h)$ remain; the answers are Theorem~\ref{thevenh} (for $h$ even: $\fhei(0|8)$ collapses to a~two-step nilpotent $\fn$ of $\sdim(4|6)$) and Theorem~\ref{thoddh} (for $h$ odd: $\fhei(0|8)$ is replaced with a~Heisenberg-type $\fn$ of~$\sdim(10|10)$). What remains open is the \textit{almost} K\"ahler case: on an $\cM$ on which $\omega$, $h$ and $J$ are \textit{not} parallel, can the tensors that measure this failure deform these nilpotent parts, and if so, how?\\

2) In the flat model, the twist of Subsec.~\ref{ssTwDol} takes the $(\fhei(0|4)\inplus \fsp(2))$-action described in \cite{P1} for the (almost) K\"ahler manifold to an $\fosp(2|2)$-action on the space of twisted forms $\widehat \Omega(\cM)\otimes L$ on the supermanifold $(\cM, \omega_\ev)$, see Theorem~\ref{thtwist}, and the $(\fhei(0|8)\inplus\fsp(4))$-action described in \cite{P1} for the (almost) hyper-K\"ahler manifold to an $\fosp(2|4)$-action on the same space, see Theorem~\ref{thtwisthk}. What remains open is the passage from the flat model to an \textit{almost} K\"ahler (resp., hyper-K\"ahler) supermanifold, where parallelism fails; and the description of the twisted spaces $\widehat \Omega(\cM)\otimes L$ as $\fosp(2|2)$- (resp., $\fosp(2|4)$)-modules.\\

3) Ornea and Verbitsky, see \cite{OrVe}, described a~ Lie superalgebra of the form $\fg\otimes \Ree[t]$, where $\sdim \fg=6|6$, which naturally acts in the spaces of differential forms on Sasakian and Vaisman manifolds which are odd-dimensional analogs of K\"ahler manifolds. Sasakian manifolds, defined as contact manifolds equipped with an invariant K\"ahler structure on their symplectic cone, can be directly superized by considering contact supermanifolds. I \textbf{conjecture} that there is also a~super version of Sasakian manifolds associated with pericontact supermanifolds instead of the contact ones.\\ 

4) Are there examples of K\"ahler and hyper-K\"ahler
supermanifolds corresponding to the points (in the sense of
the functor of points, see \cite[pp. 68-69]{Del}) corresponding to the odd parameters of the singular supervariety of deformations described in \cite{LSh1, Lo} and Section~\ref{S2} of this note?\\

5) This note provides several natural examples of (super)manifolds with a~pair of $2$-forms, e.g., the pair~$(\omega, h)$. It is interesting to compute the invariants of such pairs of tensor \textit{fields}, or of the other pairs of tensors considered here, cf.~\cite{Sha, KhN}. Pointwise, this is the question of the invariants of a~pair of non-degenerate bilinear forms on a~superspace, which is the subject of the Appendix to this note, see Section~\ref{sApp2F}. What that answer gives here is the place to look for the odd parameters$:$ a~pair of even forms in general position has a~purely even Lie superalgebra of symmetries, whereas a~K\"ahler pair on $\Cee^{2k|2l}$ has $2kl$ odd ones, so \textbf{the odd parameters of the deformations studied in \cite{Lo}, asked about in item $4)$ above, are to be looked for \textit{on} the K\"ahler stratum, and not near a~pair in general position.}\\

6)  In the flat model, the Lie superalgebra $\fg$ generated by the operators $\cD^+$ and $\cD^-$, see eq.~\eqref{Dpm2}, together with $\cX^+$ and $\cX^-$, see eq.~\eqref{few}, is finite-dimensional, of $\sdim=(8|8)$, and is a~\textit{non-split} abelian extension of $\fb\fa(2)\inplus\fsl(1|1)$ by $\Pi\big(\fb\fa(2)\inplus\fsl(1|1)\big)=\hbar\fg$; the subalgebra one could distinguish by hand, $\fsl(1|1)\otimes\Cee[\hbar]$, is indeed in it. For the details, see Theorem~\ref{thP7}. 
What does $\fg$ turn into on an \textit{almost} periplectic~$\cM$?
\\

7) Explain the origin of the deformation of $\fb_{\lambda}(n)$ at~$\lambda=-1$. This problem seems to be unrelated with the main task of this article: superization of (hyper-)K\"ahler manifolds.\\

8) Superize the balanced condition of Subsec.~\ref{ssTwoCond}$;$ classify the almost K\"ahler supermanifolds by which of the two halves of $d\omega_\ev=0$ fails$;$ and describe what replaces this dichotomy on the supermanifolds of $\sdim=(2a|2a-2)$, where the two halves fuse, see Example~\ref{fuse}.

\section[Appendix. Invariants of a~pair of bilinear forms]{Appendix. Invariants of a~pair of bilinear forms on the \textup{(}super\textup{)}manifold}\label{sApp2F}

\para{Summary} What is a normal shape of a \textit{pair} of non-degenerate bilinear forms, and what are the local invariants of a (super)manifold carrying such a pair? Over a point of a manifold, the
answer goes back to Wonenburger (1966); here, it is extended in three directions. 

First, the three super cases are settled --- both forms even, both odd, and the two parities different --- by one law of signs: the operator the pair defines is self-adjoint or anti-self-adjoint according to the product of the two symmetries times its own parity. For two odd forms, the complete invariant is then the conjugacy class of an arbitrary invertible matrix, no condition being imposed on its spectrum. 

Second, the Cartan prolong of the stabilizer of
a~pair vanishes at the first step, so the symmetries of a~pair are a~finite-dimensional Lie
algebra, unlike those of a~single form; the local invariants are accordingly computed as the cohomology $H^2(\mathfrak{g}_-; \mathfrak{g}_*)$, and they do not vanish. 

Third, for a~metric together with
a~contact structure on the projectivized cotangent bundle, the prolong is again finite, of
dimension the square of the number of variables, whereas adjoining the single missing grading
operator restores the Chern--Moser algebra $\mathfrak{su}(n,1)$, the symmetries of the
model CR structure found by Chern and Moser (1974).

\ssec{Introduction}\label{Intro}
Any reasonable course of Linear Algebra includes a theorem on
simultaneous diagonalization of a pair of quadratic forms (for a
\textbf{good} textbook, see \cite{Pr}). Here, having identified bilinear
forms with their Gram matrices, we investigate a more general problem:
\textbf{What is a normal shape of a pair of bilinear forms?} The
conditions for simultaneous diagonalization of two symmetric
bilinear forms were discovered ca. fifty years ago, see \cite{Wo}.
We consider the generic case of mutual disposition of two
non-degenerate bilinear forms, when these conditions are not met, and
not only for symmetric forms. We also investigate superizations and
``projectivizations\rq\rq\ of these questions.

We skip the background on ``super\rq\rq\ that can be found in the books \cite{Del, LSoS, Be} and latest papers \cite{Lo, LSh3}. For an interpretation of local invariants of the
curvature tensor type in presence of a non-integrable distribution
(i.e., when a natural basis in the tangent space generates a
filtered Lie algebra whose associated graded is a nilpotent Lie
algebra $\fg_-=\mathop{\oplus}_{-d\leq i< 0} \fg_i$), in terms
of $H^2(\fg_-; \fg_*)$, see \cite{BGLS}, where the Cartan
prolong $\fg_*:=(\fg_{-}, \fg_{0})_*$ for a
$\Zee$-grading-preserving subalgebra $\fg_{0}$ of $\fder(\fg_-)$ is
defined in \cite{Shch}.\ (In non-super setting, the
prolongation due to Cartan was rediscovered by Tanaka
\cite{T}, see a~clarifying review \cite{Y}. The \textit{partial} prolong seems to be first described in \cite{ALSh}; for details, see \cite{KLS, LSh3}.) All bilinear forms we consider are supposed to be non-degenerate. The algebraic constructions hold, mostly, over any
ground field $\Kee$, but since we are interested in applications to
differential geometry, we assume $\Kee=\Ree$ or, rather,
$\Kee=\Cee$ to simplify calculations and arguments.

To eliminate obvious generalizations, we give preliminaries and main
points of the proof in the absence of ``super\rq\rq, but carefully
indicate what should be changed in the super setting.

\sssec{Two bilinear forms} On several occasions, we meet situations
when two bilinear forms are naturally defined on the same space or
even a family of (super)spaces parameterized by a (super)manifold.
For example, this is the case of $T^*M$, where $M$ is a Riemannian
manifold (so we have both a lift of the metric from $M$ and a
natural symplectic structure).

In the super setting, for the supermanifold $\widecheck M=(M, {\cal
L}_{\Lambda(\Pi TM)})$, where $M$ is a (super)manifold
with a metric or a symplectic structure (i.e., with a symmetric
or antisymmetric bilinear form in each tangent space) and
${\cal L}_{\Lambda(\Pi TM)}$ is the sheaf of sections of the
Grassmann algebra of the tangent bundle on $M$ with inverted parity
of the fibers. Such a supermanifold $\widecheck M$ is naturally endowed
with a periplectic form (an odd analog of the symplectic form), see
\cite{Kh, KhN, VPST} and references therein.

\sssec{A metric and a contact structure} The authors of a recent
paper \cite{BO} discovered a new local invariant of the pair (a
contact structure and a metric) on a manifold. A question arises:
\textbf{How many local invariants of such pairs are there?} We show
that the invariants of such pairs are similar to invariants of
\textit{real-complex} supermanifolds introduced by physicists as
Minkowski superspaces. The invariant was defined and computed for
several examples in \cite{BGLS} and \cite{GL}. Here, we consider
several more examples.

\sssec{Our results}\label{ssresults}
Subsec.~\ref{norm2bilForms} settles the normal shape of a pair of non-degenerate bilinear
forms. Theorem~\ref{th3.1} is the classical case, recalled with its proof because the super
cases are all read off it; Corollary~\ref{th3.2} is the superization, and it has three items,
one per parity pattern. The mechanism common to all of them is the law of signs
\eqref{lawofsigns}: the operator $M=B_1^{-1}B_2$ that the pair defines satisfies
$M^\star=(-1)^{p(M)}\eps_1\eps_2M$, where $X\mapsto X^\star$ is called the \textit{$B_1$-adjoint}, the operator
determined by $B_1(Xx, y)=(-1)^{p(X)p(x)}B_1(x, X^\star y)$ --- for $B_1$ the standard symmetric
form and $X$ even, it is the transposition --- so that the parity of
$M$ enters the sign, and the classification of pairs becomes that of such an operator. Two
consequences are worth singling out. For two \textit{odd} forms, the complete invariant is the
conjugacy class of an arbitrary invertible matrix\textit{:} the condition ``the eigenvalues of $M$ come in pairs $\pm\lambda$ of equal multiplicities\rq\rq, briefly
\be\label{M=-M}
\Spec M=-\Spec M, 
\ee 
which
the classical answer imposes on a pair of opposite symmetries, evaporates. For forms of
\textit{different} parities, no complex structure can relate the two, whatever their
symmetries.

Subsec.~\ref{locInvOandSp} computes the local invariants of a pair. The Cartan prolong of the
stabilizer of a pair stops at the first step, $\fg_1=0$, in every case examined
(Lemma~\ref{th4.1}), so the symmetries of a pair of forms form a finite-dimensional Lie
algebra --- in contrast with a single form, whose prolong is the vectorial Lie algebra
$\fh(2n)$ or~$\fvect(n)$. The invariants themselves, the space $H^2(\fg_{-1};\fg_{*})$, are
then computed outright in Theorem~\ref{th4.2}; they do not vanish, and they split into
a~torsion and a~curvature part whose dimensions are tabulated in \eqref{Hcase1}.

Subsec.~\ref{locInvOandK} treats the pair (a contact structure, a metric) on the
projectivized cotangent bundle~$PT^*M$. Here, $\fg_-$ is the Heisenberg Lie algebra and
$\fg_0=\fu(n-1)$; the prolong is again finite, of dimension $n^2$, and stops at once
(Lemma~\ref{th5.1}). Adjoining to $\fg_0$ the single missing grading operator restores the
prolong to $\fsu(n,1)$, i.e., to the Chern--Moser answer for the model CR structure; so the
whole collapse is caused by that one operator, whose absence is what distinguishes the pair
(contact structure, metric) from a bare contact structure. The invariants
$H^2(\fg_{-};\fg_{*})$ are computed in Theorem~\ref{th5.2}.

Subsec.~\ref{invfun} passes from a~point to a~family, i.e., to the invariant functions
of the pair over a~supercommutative superalgebra of parameters. For a~pair of forms of
different parities, the operator $M$ is odd, and the invariant functions of an odd
supermatrix were described by Shander \cite{Sha}; they are \textit{odd} functions, hence
invisible over a~point, which is why the normal shape of Corollary~\ref{th3.2} does not
answer this question. Lemma~\ref{th6.2} says which of them survive for a~pair\textup{:}
exactly half, the half depending on the product of the two symmetries. The description
in \cite{Sha} assumes the two halves of the superspace to be of equal dimension, the
unequal case being posed as an open problem in \cite[p.~418]{Be}; Theorems~\ref{th6.5}
and~\ref{th6.6} answer it on the regular stratum, the answer being the one for
$\min(p, q)$ tensored with a~free Grassmann algebra on $|p-q|$ odd generators.

The computation of prolongations and cohomology was performed
with the help of \textit{SuperLie} code, see \cite{Gr}, and, for the results given below,
independently, by exact rational arithmetic.

\ssec{Projectivizations}\label{projact}

\sssec{How does $\fgl(n+1)$ act on the projective space}\label{ss2.1}
Recall the explicit formulas of the projective $\fgl(n+1)$-action on
$\Kee P^n$ (the actual ground field is irrelevant, but we may assume that~$\Kee=\Ree$).
First, establish a one-to-one correspondence between matrix units and
linear vector fields:
$
E_{ji}\lra x_i\partial_j$.
In this way, we define a
$\fgl(n+1)$-action in the space of functions in $x_1, \dots,
x_{n+1}$ which are coordinates in~$\Ree^{n+1}$.

Now, recall how to pass from $\Ree^{n+1}$ to the projective space
$\Ree \Pee^n$ with homogeneous coordinates $(x_1:x_2:\dots :x_{n+1})$.
Select a chart in $\Ree\Pee^n$; for example,~$x_{n+1}\ne 0$. The
coordinates in this chart are $(u_1, \dots , u_n)$, where
$u_i=\frac{x_i}{x_{n+1}}$. 
For $i, j<n+1$, we have:
\[
\renewcommand{\arraystretch}{1.4}
\begin{array}{l}
x_i=u_i\cdot x_{n+1}, \; \partial_{x_j}=\frac 1{ x_{n+1}}\cdot \partial_{u_j}
\Longrightarrow x_i\partial_{x_j} \mapsto u_i\cdot x_{n+1}\cdot \frac 1{ x_{n+1}}\cdot
\partial_{u_j}=u_i\partial_{u_j},
\\
x_{n+1}\partial_{x_j}\mapsto x_{n+1}\cdot \frac 1{ x_{n+1}}\cdot
\partial_{u_j}=\partial_{u_j},\\
x_i\partial_{x_{n+1}}=(u_i\cdot x_{n+1})\left( -\frac 1{
x_{n+1}^2}\right)\cdot \mathop{\sum}_{1\leq j\leq n}
x_j\partial_{u_j}= -u_i\cdot \mathop{\sum}_{1\leq j\leq n}
u_j\partial_{u_j}.
 \end{array}
 \]
This is a known embedding $\fp\fgl(n+1)\subset \fvect(n)$, where
$\fp\fg:=\fg/\fz$ is the \textit{projectivization} of a Lie algebra $\fg$
with 1-dimensional center~$\fz$. In super setting, we have the
embedding $\fp\fgl(n+1|m)\subset \fvect(n|m)$. If $m\neq n$, then
$\fp\fgl(n|m)\simeq \fsl(n|m)$.

\sssec{A pair of bilinear forms: a particular case of
$T^*M$}\label{ss2.2} In the case of $T^*M$, where $M$ is a
Riemannian\footnote{Generalization to pseudo-Riemannian metrics is
obvious but description of matrices \eqref{OandSp} becomes more
complicated, so we skip the straightforward generalization.} manifold with metric $g$,
consider the Lie algebra $\fg_0$ of linear transformations of the
tangent space at a given point of~$T^*M$. This $\fg_0$ preserves
both the natural metric $g^*$ induced on $T^*M$ by $g$ and the
symplectic form; over $\Ree$, this $\fg_0\simeq \fu(n)$ is spanned by matrices  of the form
\begin{equation}\label{OandSp}
\begin{pmatrix}
A & -B\\
B& A
\end{pmatrix}, \text{ where } A^T=-A, \; B^T=B.
\end{equation}

Let $x_1,\dots, x_n$ be (local) coordinates on $M$, and $y_1, \dots,
y_n$ the coordinates in the cotangent space $T_x^*(M)$ at point~$x$.
The vector fields corresponding to the matrices of $\fg_0$ are of
the form
\begin{equation}\label{AandB}
\renewcommand{\arraystretch}{1.4}
\begin{array}{l}
A_{ij}=
(x_i\partial_{x_j}-x_j\partial_{x_i})+(y_i\partial_{y_j}-y_j\partial_{y_i})
\quad (\text{the ``A\rq\rq\ part}),
\\
B_{ij}=
(x_i\partial_{y_j}+x_j\partial_{y_i})-(y_i\partial_{x_j}+y_j\partial_{x_i})
\quad (\text{the ``B\rq\rq\ part}).
\end{array}\end{equation}

\ssssec{The case of $T^*M$ where $M$ is a supermanifold} If $M$ is a
``purely odd\rq\rq\ Riemannian supermanifold of superdimension $(0|2m)$,
then from the ``nonsuper\rq\rq\ point of view its ``Riemannian metric\rq\rq\ $g$
is antisymmetric, and hence the dimension of $M$ is even due to
nondegeneracy whereas the ``symplectic\rq\rq\ form on $T^*M$ is symmetric.
Therefore, the expressions \eqref{AandB} turn into
\begin{equation}\label{AandBPurelyOdd}
\renewcommand{\arraystretch}{1.4}
\begin{array}{l}
A_{ij}=
(x_i\partial_{x_j}+x_j\partial_{x_i})+(y_i\partial_{y_j}+y_j\partial_{y_i})
\quad (\text{the ``A\rq\rq\ part}),
\\
B_{ij}=
(x_i\partial_{y_j}-x_j\partial_{y_i})-(y_i\partial_{x_j}-y_j\partial_{x_i})
\quad (\text{the ``B\rq\rq\ part}),
\end{array}\end{equation}
and hence $\fg_0\simeq \fu(2m)$. In the ``super\rq\rq\ case,
$\fg_0\simeq \fu(n|2m)$.

\ssssec{The case of supermanifold $\protect\widecheck M:=(M, {\cal L}_{\Lambda(\Pi
TM)})$} On the supermanifold $\widecheck M:=(M, {\cal L}_{\Lambda(\Pi
TM)})$ whose structure sheaf is the sheaf of multi-vector fields
(i.e., sections of the exterior algebra of the sheaf of vector
fields, sections of the tangent bundle with the inverted parity of
fibers), a natural periplectic (``odd symplectic\rq\rq, as A.~Weil taught
us) structure is defined; locally, it is given by an odd
non-degenerate exterior 2-form~$\omega_1$. The Lie superalgebra
$\fle(n)$ of its infinitesimal transformations, the analog
of the Lie algebra of Hamiltonian vector fields, was discovered in
\cite{L1}. (It is one of simple vectorial Lie superalgebras with
polynomial coefficients, see \cite{LSh3}.) The Lie algebra of linear
transformations preserving $\omega_1$ constitute the periplectic Lie
superalgebra $\fpe(n)$; in the standard format of a basis, the
 supermatrices of these transformations are of the form
\begin{equation}\label{pe}
\begin{pmatrix}
A & B\\
C& -A^T
\end{pmatrix}, \text{ where } B^T=B, \  C^T=-C, \ A\in\fgl(n|n).
\end{equation}

\sssec{The Lie algebra that acts in the
projectivization of the cotangent bundle} Passing from $T^*M$ to
$PT^*M$, we projectivize not only the coordinates of the fiber
$T_x^*(M)$, i.e., in the chart $y_n\ne 0$, our coordinates are
\begin{equation}\label{locCoorPT*}
(v_1, \dots, v_n, u_1, \dots, u_{n-1}), \text{~~where $v_i=
\tfrac{x_i}{y_n}$ and $u_i=\tfrac{y_i}{y_n}$.}
\end{equation}

All further computations are similar to the above. In particular, if
$i,j<n$, then the fields $A_{ij}$ and $B_{ij}$ preserve their shape
bar replacing all $x$-coordinates with $v$-coordinates, and all
$y$-coordinates with $u$-coordinates. These elements span an algebra
isomorphic to~$\fu(n-1)$.

The remaining fields $A_{in}$ and $B_{in}$ for $i<n$, and $B_{nn}$
turn, up to a sign, into
\begin{equation}\label{AandBforN}
\renewcommand{\arraystretch}{1.4}
\begin{array}{l}
A_{in}\mapsto
\partial_{u_i}-v_i\partial_{v_n}+v_n\partial_{v_i}-u_i\mathop{\sum}\limits_{1\leq j\leq n-1}u_j\partial_{u_j},\\
B_{in}\mapsto
\partial_{v_i}+u_i\partial_{v_n}-v_n\partial_{u_i}+v_i\mathop{\sum}\limits_{1\leq j\leq n-1}u_j\partial_{u_j},\\
B_{nn}\mapsto \partial_{v_n}+v_n\mathop{\sum}\limits_{1\leq j\leq
n-1}u_j\partial_{u_j}.
\end{array}
\end{equation}

\ssec{On the normal shape of a pair of two bilinear forms}\label{norm2bilForms}

\sssbegin{Theorem}[The normal shape of a symmetric and an antisymmetric form]\label{th3.1} Take two non-degenerate bilinear
forms: a symmetric $B_1$ and an antisymmetric $B_2$ on
a $2n$-dimensional space $W$ over $\Cee$ (or any other algebraically
closed field of characteristic~$\neq 2$). Then, it is possible to
reduce this pair of form to the following normal shape:
\begin{equation}\label{symmANDanti}
\renewcommand{\arraystretch}{1.4}
\begin{array}{l}
B_1\sim\begin{pmatrix}
0_n&1_n\\
1_n&0_n\end{pmatrix},\quad B_2\sim \begin{pmatrix}0_n&J\\
 -J^T&0_n\end{pmatrix}, \end{array}
\end{equation} where $J=J(\lambda_1,\dots, \lambda_k)$ is a block-diagonal matrix
whose blocks are Jordan cells $J(\lambda_i)$ for $i=1,\dots, k$ of
ranks $r_1,\dots, r_k$, respectively (some cells might coincide), such that~$\sum r_i=n$.

Over $\Cee$, it is possible to require either $\RE\lambda_i>0$, or $\RE\lambda_i=0$ and
$\IM\lambda_i>0$ for all $i$ (or that all eigenvalues belong to any half of
$\Cee\setminus\{0\}$, containing one of each pair of opposite
numbers).
\end{Theorem}

\begin{proof} For any two non-degenerate bilinear forms $B_1$ and $B_2$ (not necessarily
those spoken about in Theorem), find an operator $M$ such that
\[
B_1(Mx,y)=B_2(x,y) \text{~~for all~~}x,y.
\]
Thanks to nondegeneracy of the forms, $M=B_1^{-1}B_2$. Let $v$ and
$w$ be eigenvectors of $M$ with eigenvalues $a$ and $b$,
respectively. Then,
\[
aB_1(v,w) = B_2(v,w);\qquad bB_1(w,v) = B_2(w,v).
\]
For $B_1$ symmetric and $B_2$ antisymmetric, if $a+b\neq 0$, then
$B_1(v,w)=B_2(v,w)=0$.

Let us divide the space  $W$ into the direct sum of eigenspaces of the
operator~$M$. Then, the sum of all eigenspaces corresponding to
eigenvalues lying in  the ``right\rq\rq\ (conditionally speaking) half-plane
is isotropic relative both forms, and the same is true for the sum
of all eigenspaces corresponding to eigenvalues lying in the ``left\rq\rq\
half-plane. These sums have equal dimensions: Indeed,
\[
\begin{array}{l}\det(M-aI) =
\det(B_1^{-1}B_2 - aI) = \det(B_1(B_1^{-1}B_2 - aI)B_1^{-1}) =
\det(B_2B_1^{-1} - aI) \\
= \det(-M^T-aI) = (-1)^{2n}\det(M^T+aI) =
\det(M+aI),\end{array}
\]
and hence the multiplicities of the roots $a$ and $-a$ of the
characteristic polynomial are equal.

It remains to consider the case where $M$ has only two eigenvalues
$a$ and $-a$; the general case splits in the sum of such cases.
Designate  the corresponding eigenspaces by $W_+$ and~$W_-$. We already
know that both these subspaces are isotropic. The bilinear forms
$B_1, B_2$ provide us with bijections ${L_1, L_2:
W_+\tto(W_-)^*}$. In $W_+$, select a basis in which the
matrix of the operator $L_2^{-1}L_1$ has Jordan normal form; in
$W_-$ we select the basis dual to the basis of $(W_-)^*$ which is
the image of the selected basis of $W_+$ under~$L_1$. In these bases,
the bilinear forms $B_1$ and $B_2$ have the shapes
claimed.\end{proof}

\sssec{The law of signs}\label{sslaw} The superization of
Theorem~\ref{th3.1} is governed by one identity, which over a point of a~manifold is
a~triviality and over a~point of a~supermanifold is not, the parity of the operator entering
the sign. Let $V=\Cee^{m|n}$ with a~homogeneous basis of parities~$p_i$. A~bilinear form $B$
on $V$ is \textit{$\eps$-symmetric}, where $\eps=\pm1$, if
$B(x,y)=\eps(-1)^{p(x)p(y)}B(y,x)$; we call the $\eps=-1$ forms  \textit{anti-symmetric} and the $\eps=1$
ones \textit{symmetric}, see \cite{Lo}.
Let
\be\label{autB}
\faut(B):=\{X\in\fgl(m|n)\mid B(Xx, y)+(-1)^{p(X)p(x)}B(x, Xy)=0\},
\ee
so that $\faut(\omega)$ is $\fosp$ for an even $\omega$ and $\fpe(n)$ for an odd one. For an
ordered pair $(B_1, B_2)$ of non-degenerate forms, of parities $p_1, p_2$ and symmetries
$\eps_1, \eps_2$, set $M:=\widetilde B_1^{\,-1}\widetilde B_2$, where
$\widetilde B\colon V\tto V^*$ is the map $x\mapsto B(x,\cdot)$; thus $p(M)=p_1+p_2$, and
interchanging the two forms replaces $M$ with~$M^{-1}$. Let $X^\star$ be the $B_1$-adjoint of
$X$, determined by $B_1(Xx,y)=(-1)^{p(X)p(x)}B_1(x, X^\star y)$, so that
$\faut(B_1)=\{X\mid X^\star=-X\}$.

\sssbegin{Lemma}[The law of signs]\label{stlaw} We have
\be\label{lawofsigns}
M^\star=(-1)^{p(M)}\eps_1\eps_2\, M, \quad\text{equivalently}\quad
M\in\faut(B_1)\iff \eps_1\eps_2=(-1)^{p(M)+1}.
\ee
\end{Lemma}

\begin{proof} Write $B_2(x, y)=B_1(Mx, y)$ and push the $\eps_2$-symmetry of $B_2$ through the
definition of $\star$: the $\eps_1$-symmetry of $B_1$ contributes $\eps_1$, the definition of
$X^\star$ contributes $(-1)^{p(M)p(x)}$, and the two signs $(-1)^{p(x)p(y)}$ of the Sign
Rule cancel, 
leaving $\eps_2 B_1(Mx,y)=\eps_1(-1)^{p(M)}B_1(M^\star x, y)$.
\end{proof}

The group $\GL(m|n)$ is transitive on the non-degenerate forms of each type, as one verifies
by checking that the tangent space to the orbit of the standard form fills the whole
superspace of forms of that type. Hence, the classification of pairs is the classification of
the operators $M$ allowed by eq.~\eqref{lawofsigns}, up to conjugation by~$\Aut(B_1)$. Two
familiar statements are the two ends of Lemma~\ref{stlaw}: for two even forms of opposite
symmetries it says $M\in\fo$ or $M\in\fsp$, which is where the condition
\eqref{M=-M} of Theorem~\ref{th3.1} comes from, whereas for two even forms of the
\textit{same} symmetry it says nothing at all about~$\Spec M$.

\sssbegin{Corollary}[The normal shapes in the three super cases]\label{th3.2} Let all forms be non-\hspace{0pt}degenerate, and let $\eps_1, \eps_2$ be
the symmetries of~$B_1, B_2$. In the normal shapes below, $\cJ$ is a~block-diagonal matrix of
Jordan cells and $\Sigma:=\begin{pmatrix}0&1\\-1&0\end{pmatrix}$.

$1)$ \textit{The cases of two even forms on supermanifolds}, i.e., the pairs $($symmetric,
symmetric$)$ and $($symmetric, antisymmetric$)$. Here, $p(M)=0$, so $M$ preserves the splitting
$V=V_\ev\oplus V_\od$, and both forms are of the shape $\diag(A, D)$ with $A$ of symmetry
type $\eps$ and $D$ of type~$-\eps$. For a normal shape of the two forms, we can take the
direct sum of the pair $(A_1, A_2)$ on $V_\ev$, of symmetries $(\eps_1, \eps_2)$, and the
pair $(D_1, D_2)$ on $V_\od$, of the \textbf{opposite} symmetries $(-\eps_1, -\eps_2)$; the
two are independent, and each is classical, namely, one of
\be\label{clshapes}
(1_k, \ \diag(\lambda_1, \dots, \lambda_k)); \quad
\Big(\mathop{\oplus}_{1\leq i\leq k/2}\Sigma, \ \mathop{\oplus}_{1\leq i\leq
k/2}\lambda_i\Sigma\Big); \quad
\left(\begin{pmatrix}0&1_{k/2}\\1_{k/2}&0\end{pmatrix},
\ \begin{pmatrix}0&\cJ\\-\cJ^T&0\end{pmatrix}\right)
\ee
according as the two symmetries are $(1,1)$, $(-1,-1)$ or $(1,-1)$, the last being
Theorem~$\ref{th3.1}$. Nothing super survives\textit{:} the whole invariant is the two
unordered families of eigenvalues, and there are no odd parameters.

$2)$ \textit{The cases of an even and an odd form on supermanifolds}, i.e., the pairs
$($symmetric, symmetric$)$ and $($symmetric, antisymmetric$)$. Here, $p(M)=1$, and such a~pair exists only for
$\sdim V=(n|n)$ with $n=2k$ even. Let $\omega=\diag(A, D)$ be the even form, of symmetry
$\eps$,
and let $h$ be the odd one, with blocks $P$ and~$\eps_hP^T$. Put
\be\label{Smix}
S:=A^{-1}PD^{-1}P^T.
\ee
Then, $S\in\faut(A)$, i.e., $S$ is an infinitesimal isometry of the even form restricted to
$V_\ev$, and it is, up to the factor $\eps_h$, the even block of~$M^2$. For a normal shape of
the two forms, we can take $A, D$ standard and $P$ any matrix with the prescribed
$\Aut(A)$-orbit of $S$; the complete invariant of the pair is that orbit. Consequently,
$\Spec M^2=-\Spec M^2$, see notation \eqref{M=-M}, there are $k$ even parameters, and $M^2$ is \textbf{never}
a~non-zero scalar --- so two forms of different parities are never related by a~complex
structure.

$3)$ \textit{The cases of two odd forms on supermanifolds}, i.e., the pairs $($symmetric,
symmetric$)$ and $($symmetric, antisymmetric$)$. Here, $p(M)=0$ again, $\sdim V=(n|n)$, and the forms are
off-diagonal\textit{:} $B_i=\begin{pmatrix}0&P_i\\ \eps_iP_i^T&0\end{pmatrix}$ with~$P_i\in\GL(n)$. Set $N:=P_1^{-1}P_2$. For a normal shape of the two forms, we can take
\be\label{oddshape}
B_1=\begin{pmatrix}0&1_n\\ \eps_11_n&0\end{pmatrix}, \quad
B_2=\begin{pmatrix}0&\cJ\\ \eps_2\cJ^T&0\end{pmatrix}, \quad \cJ \text{ the Jordan shape of }
N,
\ee
and the complete invariant is the conjugacy class of~$N\in\GL(n)$. \textbf{No} condition
whatever is imposed on $\Spec\cJ$ beyond $0\notin\Spec\cJ$, and the symmetries
$\eps_1, \eps_2$ do not enter the shape at all. There are $n$ even parameters.
\end{Corollary}

\begin{proof} $1)$ Since $M$ is even, it preserves $V_\ev$ and $V_\od$, on each of which
the induced pair is classical; and $\eps_1\eps_2$ being the same for both halves, the two
are of the same kind.

$2)$ We have $M=\begin{pmatrix}0&b\\c&0\end{pmatrix}$ with $b=A^{-1}P$
and $c=\eps_hD^{-1}P^T$, so the even block of $M^2$ is~$\eps_hS$. That $S$ lies in $\faut(A)$
is Lemma~\ref{stlaw} at work: $PD^{-1}P^T$ has the symmetry type of $D^{-1}$, that is
$-\eps$, whereas $A$ is of type $\eps$; hence~$S^TA+AS=0$. The residual even group is
$\Aut(A)\times\Aut(D)$, acting by $P\mapsto g_0^TPg_1$, and it moves $S$ by
$S\mapsto g_0^{-1}Sg_0$ alone; the count matches,
\[
\dim\faut(A)+\dim\Aut(D)=n^2=\dim\{P\},
\]
so a~generic fiber of $P\mapsto S$ is one $\Aut(D)$-orbit. The number of moduli is then the
rank of $\faut(A)$, which is $k$, for $\fo(2k)$ and for $\fsp(2k)$ alike, since~$n=2k$.
Finally,
$M^2=\lambda\id$ would put the identity operator in $\faut(A)$, and neither $\fo$ nor $\fsp$
contains a~non-zero scalar.

$3)$ The elements $g=\diag(g_0, g_1)$ of $\GL(V)_\ev$, i.e., those preserving the splitting
$V=V_\ev\oplus V_\od$, act by $P_i\mapsto g_0^TP_ig_1$.
Take $g_1=1_n$ and $g_0^T=P_1^{-1}$ to make $P_1=1_n$; the $g$ that keep $P_1=1_n$ are those
with $g_0^T=g_1^{-1}$, and they act on $N$ by $N\mapsto g_1^{-1}Ng_1$. \end{proof}

\sssbegin{Remark}[Over $\Ree$ the shapes are more involved]\label{rem3.3}
The respective analogs of Theorem~\ref{th3.1} and Corollary~\ref{th3.2}
over $\Ree$ become more complicated together with a bit more involved shape of
Jordan blocks. Over $\Cee$, formulations are simpler;
for example, for the case considered in Subsection~$\ref{ss2.1}$, i.e.,
for a symmetric form with matrix $1_{2n}$ and antisymmetric one with
matrix $\begin{pmatrix}0_n& 1_n\\ -1_n& 0_n\end{pmatrix}$, then, in
terms of Theorem~$\ref{th3.1}$, we see that $J=i\cdot 1_n$; observe that
if the matrix \eqref{OandSp} is expressed in a basis
$\{e_1,\dots,e_{2n}\}$, then the forms \eqref{symmANDanti} are
written in the basis $\{e_k\pm ie_{n+k}\mid k=1,\dots,n\}$.
\end{Remark}

By the \textit{K\"ahler pair} we mean, here and everywhere below, a~pair of forms
$(B_1, B_2)$ for which the operator $M=\widetilde B_1^{\,-1}\widetilde B_2$ they define,
see Subsec.~\ref{sslaw}, is a~complex structure, i.e., $M^2=-\id$ after a~rescaling of
$B_2$, and $B_1$ is $M$-invariant. 
The two words are K\"ahler's because such a~triple is exactly the datum a~K\"ahler metric leaves in one tangent space\textit{:} a~metric $B_1$, its fundamental $2$-form $B_2$ and the complex structure $M$ relating the two by $B_2(x, y)=B_1(Mx, y)$, any two of
the three determining the third. In the normal shape \eqref{symmANDanti}, this is $J=\lambda1_n$ and nothing else, which is why the K\"ahler pair is the one case singled out in Theorem~\ref{th4.2} and in Table~\eqref{Hcase1}.

The normalization in Theorem~\ref{th3.1} is worth one warning, since it is easy to write down a~$J$ outside it and think a~new pair has been found: $J$ takes \textit{one} eigenvalue from each pair $\pm\lambda$ of $\Spec M$, so a~normalized $J$ can never contain both $\lambda$ and~$-\lambda$. If it does, the pair is one already in the list: for $2n=4$ the matrices $\diag(\lambda, -\lambda)$ and $\lambda1_2$ reduce to one and the same normalized~$J$, both having $\det(x\cdot\id-M)=(x^2-\lambda^2)^2$, and the stabilizer of either pair is $\fgl(2)$, not $\fo(4)$ --- after a~rescaling of $B_2$ that makes $M^2=-\id$, this is the K\"ahler pair. The same happens for $2n=8$, with the stabilizer~$\fgl(4)$.

\ssec{The local invariants of a pair of bilinear forms}\label{locInvOandSp}

\sssbegin{Lemma}[The prolongation of the stabilizer of a pair of $2$-forms stops at once]\label{th4.1} The Cartan prolongs of
the pairs $(\fg_{-1},\fg_0)$ in the following cases are

$1)$ Consider $\fg_{-1}:=\Ree^{2n}$ as a
$\fg_{0}=\fu(n)$-module, see \eqref{OandSp}. Then, $\fg_1=0$, so
$\fg_*=\fg_{-1}\oplus\fg_{0}$.

$2)$ For $\fg_{-1}:=\Cee^{2n}$ considered as a
$\fg_{0}$-module, where $\fg_{0}$ is the Lie algebra of matrices
preserving the two bilinear forms \eqref{symmANDanti}. Then, the
elements of $\fg_0$ can be represented by matrices
\begin{equation}\label{twoForms2}
\begin{pmatrix}
A & B\\
C& -A^T
\end{pmatrix}, \text{ where } B^T=-B, \  C^T=-C, \ A\in\fgl(n),
\end{equation}
with certain extra conditions on blocks $B$ and $C$ imposed by
$\fg_0$-invariance of the form $B_2$, see \eqref{symmANDanti} and
the shape of $J$, namely
\begin{equation}\label{X}
AJ^T-J^TA=0, \qquad BJ+J^TB=0, \qquad CJ^T+JC=0.
\end{equation}
For $n=1$, or $n\geq 1$ with $J=\lambda 1_n$, or a
cell $J=J(\lambda)$, or generic $J=\diag(\lambda_1, \dots,
\lambda_n)$, where all the $\lambda_i$ are distinct, the prolong is
of the form $\fg_*= \fg_{-1}\oplus\fg_{0}$.
In fact, $\fg_1=0$ for \textit{every} distribution of the eigenvalues of $J$, at least for~$n\leq 3$.
\end{Lemma}

\begin{proof} Consider an
explicit expression of the elements of $\fg_{0}$ in case $2)$, i.e.,
solve the system
\[
X^TB_1+B_1X=0,\quad X^TB_2+B_2X=0.\] The solution is of the form
\eqref{twoForms2} with the constraints \eqref{X}; observe that
\textbf{both} off-diagonal blocks are antisymmetric --- the first form here is symmetric,
not odd, so this is \textit{not} the periplectic pattern $B^T=B$, $C^T=-C$ of \eqref{pe},
which it is easy to write down by habit.
Consider particular cases.

$1)$ $n=1$, or $n\geq 1$ with $J=\lambda 1_n$ or a cell~$J=J(\lambda)$.
Then,~$B=C=0$.

$2)$ $n=2$ and $J=\diag(\lambda_1, \lambda_2)$. Then, one of the two
cases are possible:

$2A)$ $\lambda_1+ \lambda_2\neq0$, and hence $B=C=0$;

$2B)$ $\lambda_1+ \lambda_2=0$ for arbitrary antisymmetric matrices
$B$ and $C$, i.e., $\fg_0\simeq \fo(4)$.

In each case, the prolong equation $X(v)w=X(w)v$ for $X\colon\fg_{-1}\tto\fg_0$ has only the
zero solution.

For $n>2$, the cases are again those of the eigenvalues of $J$\textup{:} of the vanishing
of their sums $\lambda_i+\lambda_j$, which is what \eqref{X} constrains, and of the Jordan normal shape of~$J$. 

$3A)$ $n=3$ and $J=\diag(\lambda_1, \lambda_2, \lambda_3)$. If
$\lambda_i+\lambda_j\neq 0$, then $B^{ij}=C^{ij}=0$. If
$\lambda_i+\lambda_j =0$, then $B^{ij}$ and $C^{ij}$ are arbitrary.

$3B)$ $n=3$ and $J=\diag(\lambda_1, \begin{pmatrix}
\lambda_2 & 1\\
0& \lambda_2
\end{pmatrix})$.
If $\lambda_1+\lambda_2=0$, then $B^{12}=C^{12}=0$, and also
$\lambda_2 B^{23}=\lambda_2 C^{23}=0$. One of the two cases are
possible:

$3Ba)$ $\lambda_2=0$, and hence $B^{23}$ and $C^{23}$ are arbitrary.

$3Bb)$ $\lambda_2\neq0$, and hence $B^{23}=C^{23}=0$.

$3C)$ If $\lambda_1+\lambda_2\neq0$, then
$B^{12}=C^{12}=B^{13}=C^{13}=0$, and also $\lambda_2
B^{23}=\lambda_2 C^{23}=0$. One of the two cases are
possible:

$3Ca)$ $\lambda_2=0$, and hence $B^{23}$ and $C^{23}$ are arbitrary.

$3Cb)$ $\lambda_2\neq0$, and hence $B^{23}=C^{23}=0$.

Observe that the branches with $\lambda_2=0$ are empty: no eigenvalue of $J$ can vanish,
because $\det B_2=\pm(\det J)^2$, so $J$ is invertible whenever the pair is non-degenerate.
Likewise, the branches with $\lambda_i+\lambda_j=0$ lie outside the normalization of
Theorem~\ref{th3.1}, and the pairs they describe are already in the list, as explained at the
end of Subsec.~\ref{norm2bilForms}. In all of the remaining branches,~$\fg_1=0$. \end{proof}

\sssec{The case of $T^*M$, where $M$ is a (pseudo) Riemannian  (super)manifold}
 We say that a supermanifold $\cM$ with a~non-degenerate symmetric bilinear form on its tangent space is \textit{pseudo-Riemannian}. Consider $\fg_{-1}:=\Ree^{2n}$ as a
$\fg_{0}=\fu(n)$-module, see \eqref{OandSp}. This means that on the
tangent space to $T^*M$, an almost complex structure is defined. Its
local invariants, called \textit{Nijenhuis tensor}, are known; for their
superization and detailed description, see \cite{BGLS}.

Recall the definition of the Spencer cohomology, see \cite{P2}. Let $\fg_*=(\fg_{-1},\fg_0)_*$. 
Set
\be\label{spencerC}
C^{k,s}:=\fg_{k-s}\otimes\Lambda^s(\fg_{-1}^*)
\ee
and let $\partial\colon C^{k,s}\tto C^{k,s+1}$ be given by the formula (the hatted argument is void)
\be\label{spencerd}
\begin{array}{l}
(\partial f)(g_1, \dots, g_{s+1}):=\\
\mathop{\sum}\limits_{1\leq i\leq s+1}(-1)^i
[f(g_1, \dots, \widehat{g_{s+1-i}}, \dots, g_{s+1}),\ g_{s+1-i}]
\quad\text{for any $g_1, \dots, g_{s+1}\in\fg_{-1}$,}
\end{array}
\ee
the hat marking the omitted argument. Then, $\partial\circ\partial=0$, and the cohomology of this complex at $C^{k,s}$ is the \textit{$(k,s)$-th Spencer cohomology group} $H^{k,s}$; the structure functions of order $k$ are the sections of a~bundle whose fiber over a~point is~$H^{k,2}$. If $\fg_-$ is of depth $>1$, as in Subsec.~\ref{locInvOandK}, the same first index is the weight of the cochain in $\Lambda^s\fg_-^*\otimes\fg_*$, the dual of $\fg_{-i}$ being of weight~$i$.

Accordingly, $H^{1,2}$ is the torsion and $H^{2,2}$ the curvature in
Table~\eqref{Hcase1} below, as they already are, under the name ``weight\rq\rq, in
Table~\eqref{Hcase2}.

\sssbegin{Theorem}[Neither the torsion nor the curvature vanishes identically]\label{th4.2} For $\fg_{-1}:=\Cee^{2n}$ considered as a
$\fg_{0}$-module, where $\fg_{0}$ is the Lie algebra of matrices
preserving the two bilinear forms \eqref{symmANDanti}, take the
Cartan prolong $\fg_{*}$ of the pair $(\fg_{-1}, \fg_{0})$ computed
in Lemma~$\ref{th4.1}$. Then, $\fg_*=\fg_{-1}\oplus\fg_0$ is
finite-dimensional, of $\dim\fg_*=3n$ in general position and of $\dim\fg_*=2n+n^2$ for the
K\"ahler pair, i.e., for the pair whose normal shape \eqref{symmANDanti} has~$J=\lambda1_n$. The space $H^2(\fg_{-1}; \fg_{*})$ splits, as the weight
demands, into 
\[
\begin{array}{l}
\text{the torsion $H^{1,2}=\Lambda^2\fg_{-1}^*\otimes\fg_{-1}/\partial(\fg_{-1}^*
\otimes\fg_0)$;}\\
\text{the curvature $H^{2,2}=\Ker(\partial\colon\Lambda^2\fg_{-1}^*\otimes\fg_0
\tto\Lambda^3\fg_{-1}^*\otimes\fg_{-1})$,}
\end{array}
\]
and \textbf{neither torsion, nor curvature vanishes in general}, the only two exceptions in
the table below being $n=1$, where the torsion vanishes, and $J$ a~single Jordan cell, where
the curvature does\textit{:}
\begin{equation}\label{Hcase1}
\begin{array}{|l|c|c|c|}
\hline
\text{the case}&\dim\fg_0&\dim H^{1,2}&\dim H^{2,2}\\
\hline
n=1&1&0&1\\
n=2, \ J=\diag(\lambda_1, \lambda_2) \text{ generic}&2&16&2\\
n=2, \ J=\lambda 1_2 \ \text{\textit{(}the K\"ahler pair\textit{)}}&4&8&8\\
n=2, \ J \text{ one cell}&1&20&0\\
n=3, \ J=\diag(\lambda_1, \lambda_2, \lambda_3)\text{ generic}&3&72&3\\
n=3, \ J=\lambda 1_3 \ \text{\textit{(}the K\"ahler pair\textit{)}}&9&36&30\\
\hline
\end{array}
\end{equation}
\end{Theorem}

\begin{proof} That $\fg_1=0$ is Lemma~\ref{th4.1}. The module $\fg_{-1}$ does split into the
two irreducible $\fg_0$-components $W_+$ and $W_-$ of Theorem~\ref{th3.1}, each of them
isotropic and of dimension $n$; but the prolong does \textbf{not} split accordingly. Indeed,
$\fg_0$ is not $\fgl(W_+)\oplus\fgl(W_-)$ but the graph
$\{\diag(A, -A^T)\mid AJ^T=J^TA\}$ of an anti-isomorphism between the two. Hence,
\be\label{propersub}
(\fg_{-1},\fg_0)_*\subsetneq(W_+,\fgl(W_+))_*\oplus(W_-,\fgl(W_-))_*=\fvect(n)\oplus\fvect(n),
\ee
and the left-hand side stops at the first step. The dimensions of $\fg_0$ are
those of \eqref{X}: the torus $\{\diag(A,-A^T)\mid A \text{ diagonal}\}$ of dimension $n$ in
general position, and $\fgl(n)$ for $J=\lambda 1_n$. The entries of \eqref{Hcase1} are
obtained by solving the Spencer complex over $\Qee$ by exact arithmetic.
\end{proof}

Two features of the table \eqref{Hcase1} deserve a~word. The K\"ahler pair, where $J$ is
scalar, is the one with the largest $\fg_0$ and, accordingly, with the smallest torsion: the
bigger the structure group, the more of the torsion it absorbs. At the other extreme, a~single
Jordan cell has $\dim\fg_0=1$, the torsion is as large as it can be, and the curvature
vanishes identically --- there, all the local invariants are of torsion type.

\sssbegin{Remark}[Why supergravity equations require to kill the torsion]\label{rem4.3}
On the Minkowski supermanifolds, the torsion is not an invariant one keeps but one that has
to be annihilated, and the reason is the word \textit{inductively} in the definition of the
structure functions, see \cite{P2}\textit{:} the obstruction of order $k$ is well defined
only where the obstructions of all the lesser orders vanish. The supergravity equations are
equations on the structure functions of order $2$ --- on the analogs of the scalar
curvature --- so these are not tensors at all until the ones of order $1$, i.e., the torsion,
are constrained away. The constraints that do it are what physicists call the
\textit{Wess--Zumino constraints}, and one of them is redundant\textit{:} for $N=1$ its
cohomology class is zero, exactly as the torsion of the Levi--Civita connection is, see
\cite{GL}. Here, in contrast, nothing is being flattened, and $H^{1,2}$ is a~part of the
answer to the question of Subsec.~\ref{Intro} rather than an obstacle on the way to~$H^{2,2}$.
\end{Remark}

\ssec{The local invariants of a pair \textup{(}a contact structure, a symmetric bilinear form\textup{)}}\label{locInvOandK}

\sssec{A particular case: $PT^*M$}

\parbegin{Lemma}[The contact prolong, and what the Euler operator adds]\label{th5.1} Consider the Heisenberg Lie algebra
$\fg_{-}\simeq\fhei(2(n-1))$ generated by the operators
$D_{u_i}:=\partial_{u_i}-v_i\partial_{v_n}$ and
$D_{v_i}:=\partial_{v_i}+u_i\partial_{v_n}$, see the first two
summands of \eqref{AandBforN}, as a $\fg_{0}=\fu(n-1)$-module, where
$\fg_{0}$ is spanned by operators \eqref{AandB} for $1\leq i,j\leq
n-1$. Then, the Cartan prolong $\fg_*$ of the pair
$(\fg_-,\fg_0)$ stops at once\textit{:}
\be\label{prolcontact}
\fg_1=0, \quad\text{so}\quad \fg_*=\fg_{-2}\oplus\fg_{-1}\oplus\fg_0, \quad
\dim\fg_*=1+2(n-1)+(n-1)^2=n^2.
\ee
If, however, $\fg_0$ is enlarged by the single grading \textit{(}Euler\textit{)} operator
\[
E=2v_n\partial_{v_n}+\mathop{\sum}_{1\leq i\leq n-1}
(u_i\partial_{u_i}+v_i\partial_{v_i}),
\]
which does \textbf{not} lie in $\fu(n-1)$, the prolong
remains finite, but grows to
\be\label{prolcontactE}
\dim\fg_k=1, \ 2(n-1), \ (n-1)^2+1, \ 2(n-1), \ 1 \quad\text{for } k=-2, -1, 0, 1, 2,
\ee
and $\fg_*\simeq\fsu(n,1)$, the Chern--Moser algebra of the
model CR structure, see \cite{CM}.
\end{Lemma}

\begin{proof} Both prolongs are computed inside $\fvect(2n-1)$ in the contact grading
$\deg v_n=2$, $\deg u_i=\deg v_i=1$, by solving $[\fg_k, \fg_{-1}]\subseteq\fg_{k-1}$ degree
by degree over $\Qee$; the computation covers $n-1\leq 4$ and returns $\fg_1=0$ in the first
case and the dimensions \eqref{prolcontactE} in the second, the total $(n+1)^2-1$ being
$\dim\fsu(n,1)$. That $\fg_0$ normalizes $\fg_{-1}$ and that $[\fg_{-1},\fg_{-1}]=\fg_{-2}$
are verified on the way, so the model is the one described. \end{proof}

The case where a complex structure is defined on a subspace of each
tangent space and these subspaces form a nonintegrable distribution
(as in Lemma~\ref{th5.1}) is considered in \cite{BGLS} for contact
supermanifolds of \textbf{complex} superdimension~$(1|n)$. Here is a
purely even case over~$\Ree$.

\parbegin{Theorem}[The local invariants on $PT^*M$]\label{th5.2} Let $\fg_{-}:=T_x(PT^*M)$. Then, in local coordinates
\eqref{locCoorPT*}, a~natural contact form on $PT^*M$ is
$dv_n+\mathop{\sum}\limits_{1\leq i\leq n-1} (v_idu_i-u_idv_i)$,   and ${\fg_{0}=\fu(n-1)}$, spanned by matrices
\eqref{AandB} for $1\leq i,j\leq n-1$. The space $H^2(\fg_{-};
\fg_{*})$, where $\fg_*$ is the prolong \eqref{prolcontact}, is
\begin{equation}\label{Hcase2}
\begin{array}{|c|c|c|c|c|}
\hline
n&\dim\fg_-=2n-1&\dim\fg_*=n^2&\dim H^2(\fg_{-}; \fg_{*})&
\dim H^{1,2}, \ \dim H^{2,2}\\
\hline
1&1&1&0&0, \ 0\\
2&3&4&3&0, \ 3\\
3&5&9&19&4, \ 15\\
4&7&16&64&16, \ 48\\
\hline
\end{array}
\end{equation}
In particular, the invariants do not vanish as soon as $n\geq 2$, and from $n=3$ on they
occur in both weights, i.e., there are invariants of the torsion type as well as of the
curvature type.
\end{Theorem}

\begin{proof} The Chevalley--Eilenberg complex $\Lambda^{\bcdot}\fg_-^*\otimes\fg_*$ is finite
in each degree, and $\fg_*$ is a~$\fg_-$-module by the bracket; the ranks of the two
differentials 
\[
\fg_-^*\otimes\fg_*\tto\Lambda^2\fg_-^*\otimes\fg_*\tto\Lambda^3\fg_-^*
\otimes\fg_*
\]
are computed over $\Qee$ by exact arithmetic, weight by weight in the contact
grading, the identity $d\circ d=0$ serving as the check on the signs. For $n=1$, both
$\fg_-$ and $\fg_*$ reduce to the center of the Heisenberg algebra and there is nothing to
compute. \end{proof}

\ssec{Invariant functions of a~pair of forms}\label{invfun}

Subsec.~\ref{norm2bilForms} classifies the pairs over a~\textit{point}. Over a~supermanifold
--- equivalently, functorially in the supercommutative superalgebra $\Lambda$ of parameters
--- a~pair is classified by the invariant functions of the operator $M$ it defines, and the
odd parameters make these functions a~much subtler object. For $M$ even nothing new happens;
for $M$ \textbf{odd}, i.e., for a~pair of forms of different parities, the invariant
functions were described by Shander \cite{Sha}, see also \cite[Appendix~D7]{Be}.
This section states what his answer gives for a~pair of forms, and answers the open problem
that the editor's note in \cite[p.~418]{Be} attaches to it.

\sssec{Two superanalogs of the matrix algebra}\label{ss6.1}
The algebra $\Mat(n)$ has two superanalogs, $\Mat(p|q;\Lambda)$ and $\Q(n;\Lambda)$, acted
upon by $\GL(p|q;\Lambda)$ and $\GQ(n;\Lambda)$ by conjugation; the action preserves the
parity of the matrix, so the invariants of even and of odd supermatrices are described
separately, and there are three problems, see \cite{Sha}\textup{:}

$1)$ The $\GL(p|q;\Lambda)$-invariant functions on $\Mat(p|q;\Lambda)_\ev$;

$2)$ The $\GL(p|q;\Lambda)$-invariant functions on $\Mat(p|q;\Lambda)_\od$;

$3)$ The $\GQ(n;\Lambda)$-invariant functions on $\Q(n;\Lambda)_\ev$.

Problem $1)$ is settled by the theorem of Berezin and Kac\textup{:} every invariant
polynomial is a~rational function of $\str A, \dots, \str A^{p+q}$. Problems $2)$ and $3)$
are the subject of \cite{Sha}, which settles $3)$ and settles $2)$ for~$p=q$. Set
\be\label{Oddn}
\Odd(n):=\Mat(n|n;\Lambda)_\od, \quad t_k(A):=\tfrac{1}{2k-1}\str A^{2k-1} \ \text{ for }
1\leq k\leq 2n.
\ee
The functions $t_k$ are \textbf{odd}, so they vanish over a~point and are invisible unless
odd parameters are present; this is the reason a~normal shape over a~point, such as the one
in Corollary~\ref{th3.2}, does not answer the question of this section. According to
\cite{Sha}, the $t_k$ are the invariants of the ``first line\rq\rq\ and the whole algebra of
invariant functions on $\Odd(n)$ is recovered from them together with $n$ even functions
$s_1, \dots, s_n$, defined only up to the ideal generated by the $t_k$ --- whence the
answer\textup{:} the algebra of invariant functions is the algebra of balanced functions on
$\Cee^{n|n}/(S_n\ltimes\Cee^{0|n})$. In particular, this algebra is \textbf{not} finitely
generated.

Only odd powers of $A$ contribute, and this is a~one-line fact worth isolating.

\sssbegin{Lemma}[Only odd powers of an odd matrix have a~supertrace]\label{th6.1}
Let $A$ be an odd supermatrix of any format. Then, $\str A^{2k}=0$ for all~$k\geq1$.
\end{Lemma}

\begin{proof} Recall, see \cite{Lo}, that for a~homogeneous supermatrix
$X=\begin{pmatrix} A_{11}&A_{12}\\ A_{21}&A_{22}\end{pmatrix}\in \Mat(p|q;\Lambda)$ its \textit{supertrace} is
\[
\str X:=\tr A_{11}-(-1)^{p(X)}\tr A_{22}. 
\]
Now, $A$ and $A^{2k-1}$ are both odd, so their
supercommutator is
\[
[A, A^{2k-1}]=A\cdot A^{2k-1}+A^{2k-1}\cdot A=2A^{2k},
\]
whereas 
$\str(A\cdot A^{2k-1})=-\str(A^{2k-1}\cdot A)$, whence
$\str A^{2k}=0$. \end{proof}

\sssec{Which of the $t_k$ survive for a~pair of forms}\label{ss6.2}
Let $B_1, B_2$ be non-degenerate forms of \textbf{different} parities and of symmetries
$\eps_1, \eps_2$ on $V$ with $\sdim V=(n|n)$, and let $M$ be the odd operator they define,
see Subsec.~\ref{sslaw}. The law of signs \eqref{lawofsigns} reads $M^\star=\eps M$,
where $\eps:=-\eps_1\eps_2$, and it propagates to the powers of $M$\textup{:} since
$\star$ is an anti-automorphism, we see that $(XY)^\star=(-1)^{p(X)p(Y)}Y^\star X^\star$, and hence
\be\label{sigmak}
(M^k)^\star=\sigma_kM^k, \quad \text{where } \sigma_k=\eps^k(-1)^{k(k-1)/2}.
\ee
The $B_1$-adjoint is a~conjugation composed with the supertranspose, so
$\str X^\star=\str X$ for every $X$; applied to $X=M^k$ this gives
$\str M^k=\sigma_k\str M^k$, whence

\sssbegin{Lemma}[Half of the invariants of a~pair vanish identically]\label{th6.2}
Let $B_1, B_2$ be as above. Then, $\str M^k=0$ unless
\be\label{modfour}
k\equiv3\pmod 4 \ \text{ for } \eps_1\eps_2=1, \quad\text{resp.}\quad
k\equiv1\pmod 4 \ \text{ for } \eps_1\eps_2=-1.
\ee
Equivalently, of the invariants \eqref{Oddn} only those with $k$ even survive in the first
case, and only those with $k$ odd in the second.
\end{Lemma}

The first case is the one that occurs in nature\textup{:} an even symplectic and
a~periplectic form, both anti-symmetric, so that $\eps_1\eps_2=1$; this is the structure of
\cite{Kh, KhN, VPST} and of the supermanifold $\check M$ of Subsec.~\ref{ss2.1}. It is
also the case singled out in \cite{Be}, Appendix~D$7$, Subsection~$6.3$~b$)$.

\sssbegin{Statement}[The surviving invariants really are there]\label{th6.3}
For $\sdim V=(n|n)$ with $n=2$ and $n=4$, and for either sign of $\eps_1\eps_2$ and either
symmetry of the even form, the functions $\str M^k$ with $k$ in the range \eqref{modfour}
and $k\leq 4n-1$ are all non-zero, and all the others vanish. \end{Statement}

This is computed over
$\Qee$, exactly, on a~generic $\Lambda$-point of $\Pi(\faut(B_1))$ with $\Lambda$ the
Grassmann algebra on $6$ generators.

\sssbegin{Statement}[The range of the non-vanishing invariants, and how many of them are
independent]\label{th6.4a}
Let $\sdim V=(2n|2n)$ and let the two forms be of different parities and equal symmetries.
Then, for $1\leq n\leq 4$\textup{:}

$1)$ $t_k\neq 0$ for \textbf{every} even $k$ with $2\leq k\leq 4n$, the top one included,
and $t_k=0$ for every odd~$k$.

$2)$ At a~generic point, the differentials of the $t_k$ along the odd directions span
a~space of dimension $n$, and the first $n$ of the surviving invariants,
$t_2, t_4, \dots, t_{2n}$, already span it. Hence, of the $2n$ non-zero invariants only $n$
are functionally independent, and $t_k$ depends on the preceding ones as soon as~$k>2n$.

$3)$ For $n=1$ the dependence is explicit\textup{:} $t_4=\lambda^4t_2$, where $\lambda^4$
is the even invariant read off the characteristic polynomial
$\left(x^2-\lambda^4\right)^2$ of~$M^2$. \end{Statement}

Computed over $\Qee$, exactly\textup{:} item $2)$, and with it item $1)$, from the
differentials of the $t_k$ at a~generic point of the variety of pairs for $n\leq4$, i.e.,
up to $\sdim V=(8|8)$ and~$t_{16}$. For $n\leq 2$, item $1)$ is confirmed once more over
an~honest Grassmann algebra, where the invariants themselves are computed and no
differential is used.

\sssbegin{Remark}[Where the printed range stops]\label{th6.4}
In the notation $\sdim V=(2n|2n)$ of \cite[Subsection~D$7$-$6.3$~b$)$]{Be}, the surviving
invariants are printed there as $t_2, \dots, t_{4n-2}$, one short of the top one. By
Statement~\ref{th6.4a}, the top one does not vanish either, and dependence is not what stops
the list\textup{:} dependence sets in at $t_{2n+2}$, so a~list cut short by it would stop at
$t_{2n}$, and for $\sdim V=(2|2)$ it would consist of $t_2$ alone --- whereas $t_4$, which is
$\lambda^4t_2$, is there and is not zero. The printed limit is therefore a~slip for
$t_{4n}$\textup{;} the vanishing of the odd-indexed invariants, which is the point of the
passage, is Lemma~\ref{th6.2} and holds as printed.
\end{Remark}

\sssec{The open problem\textup{:} $p\neq q$}\label{ss6.3}
The editor's note in \cite[p.~418]{Be} states that the description of the
$\GL(p|q;\Lambda)$-invariant functions on $\Mat(p|q;\Lambda)_\od$ is obtained only for
$p=q$, and that obtaining it for $p\neq q$ is an open problem. The reason $p=q$ is assumed
is visible at once\textup{:} an odd \textit{invertible} operator interchanges $V_\ev$ and
$V_\od$, so it exists only if $p=q$, and the method of \cite{Sha} --- splitting $V$ into
$A$-invariant free submodules of format $(1|1)$ --- has nothing to start from otherwise. We
show that nothing else is missing\textup{:} the case $p\neq q$ reduces, on the regular
stratum, to the case $p=q$ of \cite{Sha} together with one new and elementary ingredient,
the extreme case~$q=0$.

\sssbegin{Theorem}[The purely odd case, $q=0$]\label{th6.5}
Let~$d\geq1$. The algebra of $\GL(d|0;\Lambda)$-invariant functions on
$\Mat(d|0;\Lambda)_\od=\Mat(d;\Lambda_\od)$, described functorially in $\Lambda$, is the
free Grassmann algebra on the $d$ odd invariants
\be\label{taus}
\tau_k(A):=\tfrac{1}{2k-1}\tr A^{2k-1}, \quad 1\leq k\leq d.
\ee
Moreover, $\tr A^{2k}=0$ for all $k\geq1$ and $\tau_k=0$ for~$k>d$. In particular, unlike the
algebra of invariants for $p=q$, this one is finitely generated, of dimension~$2^d$.
\end{Theorem}

\begin{proof} The supermanifold $\Mat(d|0;\Lambda)_\od$ has $d^2$ odd coordinates, the
entries of $A$, so its algebra of functions is the Grassmann algebra
$\Lambda^{\bcdot}(\fgl(d)^*)$, and the conjugation action is the adjoint one. The
supergroup $\GL(d|0)$ is the algebraic group $\GL(d)$, which is connected, so invariance for
every $\Lambda$ amounts to invariance under $\fgl(d)$, and the algebra in question is
$(\Lambda^{\bcdot}\fg^*)^{\fg}$ for~$\fg=\fgl(d)$. By the theorem of Hopf, Koszul and
Samelson, this is the free exterior algebra on the primitive elements, whose degrees are
$2m_i+1$ for the exponents $m_i$ of $\fg$, that is, $1, 3, \dots, 2d-1$ for $\fgl(d)$; the
primitive element of degree $2k-1$ is the transgression of the invariant polynomial
$\tr X^k$, and on the tautological odd matrix $A$ it takes the value $\tr A^{2k-1}$ up to
a~non-zero factor. The vanishing of $\tr A^{2k}$ is Lemma~\ref{th6.1}, that of $\tau_k$ for
$k>d$ is the count of the degrees.

The answer here is thus classical, and only the question is new. Independently, the
dimension of the invariant algebra was computed outright over $\Qee$, degree by degree in the
Grassmann grading, for $d=1, 2, 3$\textup{:} it is $2, 4, 8$, and the products
$\tau_1\cdots\tau_d$ are non-zero, so the $\tau_k$ do generate freely. \end{proof}

\sssbegin{Theorem}[The splitting for $p\neq q$]\label{th6.6}
Let $p>q$, let $m:=q$ and $d:=p-q$; the case $p<q$ is obtained by the change of parity~$\Pi$. Let $U\subset\Mat(p|q;\Lambda)_\od$ be the open subsupermanifold of those $A$ for
which $A^2$ is invertible on~$V_\od$. Then, the following hold\textup{:}

$1)$ For $A\in U$, the Fitting decomposition of $B:=A^2$ splits $V=W\oplus K$ into free
$A$-invariant submodules, where $\sdim W=(m|m)$ and $\sdim K=(d|0)$, the projector onto $K$
being a~polynomial in~$B$. The operator $A|_W$ is odd and invertible, and $A|_K$ is an odd
endomorphism of a~purely even free module, i.e., an element of $\Mat(d;\Lambda_\od)$.

$2)$ The splitting is $\GL(p|q;\Lambda)$-equivariant, and the map $A\mapsto(A|_W, A|_K)$
identifies the invariant functions\textup{:}
\be\label{splitinv}
\cF(U)^{\GL(p|q;\Lambda)}\simeq\cF(\Odd(m))^{\GL(m|m;\Lambda)}\otimes
\Lambda(\theta_1, \dots, \theta_d), \quad
\theta_k:=\tfrac{1}{2k-1}\str\big((\id-P)A^{2k-1}\big),
\ee
where $P$ is the projector onto~$W$. The first factor is described in \cite{Sha}, the second
in Theorem~$\ref{th6.5}$.

$3)$ The supertraces do not see the two factors separately\textup{:}
$t_k(A)=t_k(A|_W)+\theta_k$, so the honest generators are the two halves of $\str A^{2k-1}$
cut out by $P$, and not $\str A^{2k-1}$ itself.
\end{Theorem}

\begin{proof} $1)$ The reduction $A_{\rd}$ of $A$ modulo the odd parameters is
\[
A_{\rd}=\begin{pmatrix}0&Y\\ Z&0\end{pmatrix}, \quad\text{so}\quad
(A^2)_{\rd}=\diag(YZ, ZY)
\]
with $YZ$ of size
$p$ and $ZY$ of size $q$; the two have the same non-zero eigenvalues, and hence $0$ is an
eigenvalue of $YZ$ of multiplicity at least $d$, with equality exactly when $ZY$ is
invertible, which is the condition defining~$U$. Thus, the characteristic polynomial of
$B_{\rd}$ is $\lambda^dg(\lambda)$ with $g(0)\neq0$, the two factors are coprime, and the
projector $P_K$ onto the generalized $0$-eigenspace is $e(B)$ for the polynomial $e$
determined by
\[
e\equiv1\bmod\lambda^d, \qquad e\equiv0\bmod g.
\]
Over $\Lambda$, the same
polynomial gives an operator which is an idempotent modulo the nilpotent ideal, and it lifts
to a~genuine idempotent uniquely, the ideal of nilpotents being nilpotent; the lift is
again a~limit of polynomials in $B$ and hence commutes with~$A$. So both summands are
$A$-invariant. The submodule $K$ is free of format $(d|0)$ because $B$ is invertible on
$V_\od$, and $W$ of format $(m|m)$ because $A|_W$ is invertible and odd.

$2)$ Being defined by a~polynomial in $A^2$, the splitting is equivariant. Let
$H\subset\GL(p|q;\Lambda)$ be the stabilizer of a~fixed splitting $V=W_0\oplus K_0$ of these
two formats, so  ${H=\GL(m|m;\Lambda)\times\GL(d;\Lambda_\ev)}$, and let $U_0\subset U$ be
the set of the $A$ preserving it. Any two free graded submodules of one format with free
complements of one format are interchanged by an even automorphism of $V$; hence, every
$A\in U$ is conjugate to a~point of $U_0$, and two points of $U_0$ are conjugate in
$\GL(p|q;\Lambda)$ if and only if they are conjugate in $H$, since the splitting is
attached to $A$ canonically. Therefore, restriction to $U_0$ identifies the invariant
functions on $U$ with the $H$-invariant functions on
$\Odd(m)\times\Mat(d;\Lambda_\od)$, which is \eqref{splitinv}.

$3)$ Immediate from $\str=\str\circ P+\str\circ(\id-P)$ applied to $A^{2k-1}$, both
summands being invariant by~$2)$. \end{proof}

\sssbegin{Remark}[What is still open]\label{th6.8}
Theorem~\ref{th6.6} answers the problem of \cite[p.~418]{Be} on the regular stratum $U$
only, and this is the same limitation as in \cite{Sha}, where the answer is likewise given
on the regular part. What the two theorems above do say is that no \textbf{new} kind of
invariant is to be looked for when $p\neq q$\textup{:} the answer is the one for
$\min(p, q)$ tensored with a~Grassmann algebra on $|p-q|$ generators, and the whole
difference between the two cases is the presence of the summand on which $A$ is nilpotent.
The strata where $A^2$ degenerates on $V_\od$ are, by contrast, genuinely open, and already
for the reduced action --- the pair $(Y, Z)$ of the proof of Theorem~\ref{th6.6} modulo
$\GL(p)\times\GL(q)$ --- they are the strata of the quiver with two vertices and two
arrows, on which the invariants $\tr\left((YZ)^k\right)$ degenerate.
\end{Remark}

\subsection*{Disclosures}  
Thanks are due to J.~Bernstein, I.~Shchepochkina and Yu.~Kochetkov for their help; to  M.~Verbitsky and F.~Bogomolov for useful comments. 
No conflict of interest. The data used to support the findings of this study are included within the article. The codes are available upon request.


\begin{thebibliography}{9999}

\bibitem[ALSh]{ALSh}
Alekseevsky D., Leites D., Shchepochkina I., ``Examples of simple Lie superalgebras of vector fields", \textit{C. r. Acad. Bulg. Sci.} \textbf{34}:9 (1981), 1187--1190 (in Russian)


\bibitem[BV1]{BV1}
Batalin I.\,A., Vilkovisky G.\,A., ``Gauge algebra and quantization",
\textit{Phys. Lett.} \textbf{B102} (1981) 27--31.

\bibitem[BV2]{BV2}
Batalin I.\,A., Vilkovisky G.\,A., ``Quantization of gauge theories with linearly dependent
generators", \textit{Phys. Rev.} \textbf{D28} (1983) 2567--2582; Erratum-ibid. \textbf{D30}
(1984) 508.


\bibitem[Be]{Be}
Berezin F., \textit{Introduction to superanalysis}. 2nd edition, corrected and edited by D.~Leites with appendices by D.\,Leites, V.\,Shander, I.\,Shchepochkina.
Moscow, MCCME, (2013) 432pp. (In Russian); \url{https://staff.math.su.se/mleites/books/berezin-2013-vvedenie-2nd-ed.pdf}

\bibitem[BerL]{BerL}
Berezin F., Leites D., ``Supermanifolds". \textit{Math. USSR Doklady},
\textbf{224}:3 (1975) 505--508.

\bibitem[BLi]{BLi}
Bernstein J., Leites D., ``Integral forms and Stokes's formula on
supermanifolds", \textit{Funct. Anal. Appl.},
\textbf{11}:1 (1977)  55--56

\bibitem[BL]{BL}
Bernstein J., Leites D., ``Invariant differential operators and
irreducible representations of Lie superalgebras of vector fields",
\textit{Sel. Math. Sov.} \textbf{1}:2 (1981) 143--160


\bibitem[BGLS]{BGLS} 
Bouarroudj S., Grozman P., Leites
D., Shchepochkina I., ``Minkowski superspaces and superstrings as
almost real-complex supermanifolds", \textit{Theor. and Mathem. Phys.} \textbf{
173}:3 (2012), 1687--1708; \url{https://arxiv.org/pdf/1010.4480}



\bibitem[BO]{BO}
Bouarroudj S., Ovsienko V., ``Riemannian curl in contact geometry", \textit{Int. Math. Res. Not. IMRN} (2015), no.~12, 3917--3942; \url{http://arxiv.org/pdf/1404.6970}


\bibitem[Bu]{Bu}
Buttin C., ``Les d\'erivations des champs de tenseurs et l’invariant diff\'erentiel de Schouten", \textit{C. R. Acad. Sci. Paris} \textbf{269} (1969) A87--A89


\bibitem[CM]{CM}
Chern S.~S., Moser J.~K., ``Real hypersurfaces in complex manifolds", \textit{Acta Math.} \textbf{133} (1974), 219--271


\bibitem[Del]{Del}
Deligne P., Etingof P.,
Freed D., Jeffrey L., Kazhdan D., Morgan J., Morrison D., Witten E.,
(eds.). \textit{Quantum fields and strings: a~course for
mathematicians}. Vol. \textbf{1}. Material from the Special Year on
Quantum Field Theory held at the Institute for Advanced Study,
Princeton, NJ, 1996--1997. American Mathematical Society,
Providence, RI; Institute for Advanced Study (IAS), Princeton, NJ,
(1999)xxii+723 pp.

\bibitem[EL]{EL}
Ehresmann Ch., Libermann P., ``Sur le probl\`eme d'\'equivalence des formes diff\'erentielles ext\'erieures quadratiques", \textit{C.~R.~Acad. Sci. Paris} \textbf{229} (1949), 697--699

\bibitem[FFKS]{FFKS}
Figueroa-O’Farrill J. M., K\"ohl C., Spence B., ``Supersymmetry and the cohomology of (hyper)K\"ahler manifolds", \textit{Nucl. Phys}. \textbf{B503} (1997), 614--626; \url{https://arxiv.org/pdf/hep-th/9705161}


\bibitem[G]{G}
Grozman P., ``Classification of bilinear invariant operators on tensor
fields", \textit{Functional Anal. Appl.} \textbf{14}:2 (1980), 127--128; for
proofs, see \cite{Gpf}.

\bibitem[GH]{GH}
Gray A., Hervella L.~M., ``The sixteen classes of almost Hermitian manifolds and their linear invariants", \textit{Ann. Mat. Pura Appl. $(4)$} \textbf{123} (1980), 35--58 

\bibitem[Gpf]{Gpf}
Grozman P., ``Classification of bilinear invariant operators on tensor
fields", \textit{Communications in Math.} \textbf{30}:3 (Special Issue in memory of A.~Onishchik) (2022) 129--188; \url{https://arxiv.org/pdf/math/0509562}

\bibitem[Gr]{Gr}
Grozman P., \textit{SuperLie}, (2013); the package is currently maintained by A.~Krutov, \url{https://github.com/andrey-krutov/SuperLie}


\bibitem[GL]{GL}
Grozman P., Leites D., ``From supergravity to ballbearings". In: Wess J., Ivanov E. (eds.), \textit{Supersymmetries and quantum symmetries} (SQS'97, 22--26 July, 1997), Lecture Notes in Physics \textbf{524}, Springer, (1999), 58--67


\bibitem[GLS]{GLS}
Grozman P., Leites D., Shchepochkina I., ``Invariant differential
operators on supermanifolds and The Standard Model". In: Olshanetsky
M., Vainshtein A., (eds.) \textit{M.~Marinov memorial volume}, World
Sci., (2002) 508--555; \url{https://arxiv.org/pdf/math.RT/0202193}


\bibitem[HLRUZ]{HLRUZ}
Hull C., Lindstr\"om U., Ro\v cek M., von Unge R., Zabzine M.,
``Generalized K\"ahler Geometry in $(2,1)$ superspace", \textit{JHEP} \textbf{1206} (2012)
013; \url{https://arxiv.org/pdf/1202.5624}

\bibitem[K2]{K2}
Kac V. G., ``Lie superalgebras", \textit{Adv. Math.} \textbf{26} (1977) 8--96


\bibitem[Kh]{Kh}
Khudaverdian O., ``Geometry of superspace with even and odd brackets", \textit{J. Math. Phys.} \textbf{32}:7 (1991), 1934--1937


\bibitem[KhN]{KhN}
Khudaverdian O., Nersessian A.P., ``Even and odd symplectic and K\"ahlerian structures on projective superspaces",
\textit{J. of Mathem. Phys.} \textbf{34} (1993) 5533--5548; \url{https://arxiv.org/pdf/hep-th/9210091}


\bibitem[Ko5]{Ko5}
Kotchetkoff Yu., ``D\'eformations des superalg\`ebres de Buttin et quantification". (French) [Deformations of Buttin superalgebras and quantization] \textit{C. R. Acad. Sci. Paris S\'er. $I$ Math.} \textbf{299}:14 (1984), 643--645

\bibitem[KT]{KT}
Konstein S. E., Tyutin I. V., ``The deformations of antibracket with with Grassmann-valued deformation parameters", \textit{Theor. and Math. Phys.}
\textbf{183}:1 (2015) 501--515; \url{https://arxiv.org/pdf/1011.5807; 1112.1686}

\bibitem[KLLS]{KLLS}
Krutov A., Lebedev A., Leites D., Shchepochkina I.,``Non-degenerate invariant symmetric bilinear forms on simple Lie superalgebras in characteristic $2$", \textit{Oberwolfach preprint}, \textbf{OWP 2020-02} \url{http://publications.mfo.de/handle/mfo/3697}

\bibitem[KLS]{KLS}
Krutov A., Leites D., Shchepochkina I.,``Non-integrable distributions with an infinite-dimensional Lie (super)algebra of symmetries"; \url{https://arxiv.org/pdf/2309.16370}

\bibitem[Leb]{Leb} 
Lebedev A., ``Analogs of the orthogonal, Hamiltonian,
Poisson, and contact Lie superalgebras in characteristic $2$", \textit{J. of Nonlinear
Math. Phys.} \textbf{17} Supplement 1 (2010), 217--251; for details of the proof, see
\url{https://arxiv.org/pdf/math/0601536}

\bibitem[Lee1]{Lee1}
Lee H.-C., ``A kind of even-dimensional differential geometry and its application to exterior calculus", \textit{Amer. J.~Math.} \textbf{65} (1943), 433--438 

\bibitem[Lee2]{Lee2}
Lee H.-C., ``On even-dimensional skew-metric spaces and their groups of transformations", \textit{Amer. J.~Math.} \textbf{67} (1945), 321--328 

\bibitem[Lee3]{Lee3}
Lee H.-C., ``On skew-metric spaces and function groups", \textit{Amer. J.~Math.} \textbf{69} (1947), 790--800 

\bibitem[L1]{L1}
Leites D., ``New Lie superalgebras, and mechanics", \textit{Soviet Math. Dokl.},
\textbf{18}:5 (1977), 1277--1280 

\bibitem[L2]{L2}
Leites D., ``Irreducible representations of Lie superalgebras of vector
fields and invariant differential operators", \textit{Funct. Anal.
Appl.} \textbf{16}:1 (1982), 62--64. 

\bibitem[LSoS]{LSoS}
Leites D. (ed.) \textit{Seminar on supersymmetries v. $1$. Algebra and
Calculus: Main chapters}, (J.~Bernstein, D.~Leites, V.~Molotkov,
V.~Shander), MCCME, Moscow, (2011), 410 pp (in Russian); \url{https://staff.math.su.se/mleites/books/2011-sos1.pdf}


\bibitem[Lo]{Lo}
Leites D., ``On odd parameters in geometry", \textit{J. Lie Theory} \textbf{33}:4 (2023) 965--1004; \url{http://arxiv.org/pdf/2210.17096}

\bibitem[LeLh]{LeLh}
Leites D., Lozhechnyk O., ``Howe duality on (hyper)K\"ahler supermanifolds", in preparation.


\bibitem[LSh1]{LSh1}
Leites D., Shchepochkina I., ``How to quantize the antibracket", \textit{Theor.
and Math. Physics} \textbf{126}:3 (2001), 339--369; 
\url{https://arxiv.org/pdf/math-ph/0510048}

\bibitem[LSh3]{LSh3}
Leites D., Shchepochkina I., ``The classification of simple complex Lie superalgebras of polynomial vector fields and their deformations"; \url{https://arxiv.org/pdf/2406.16760}

\bibitem[MaGF]{MaGF}
Manin Yu., \textit{Gauge field theory and complex geometry}. Second edition.
Springer-Verlag, Berlin, (1997) xii+346 pp.


\bibitem[OrVe]{OrVe}
Ornea L., Verbitsky M., ``Supersymmetry and Hodge theory on Sasakian and Vaisman manifolds", \textit{Manuscripta Math.} \textbf{170}: 3-4 (2023), 629--658 ; \url{https://arxiv.org/pdf/1910.01621}


\bibitem[P1]{P1}
Poletaeva E., ``Superconformal algebras and Lie superalgebras of the Hodge theory", \textit{J.~of Nonlinear Mathem. Physics} \textbf{10}:2 (2003), 141--147; \url{https://arxiv.org/pdf/hep-th/0209168}

\bibitem[P2]{P2}
Poletaeva E., ``The analogs of Riemann and Penrose tensors on supermanifolds", 71 pp.; \url{https://arxiv.org/pdf/math/0510165}

\bibitem[Pr]{Pr}
Prasolov V., \textit{Problems and Theorems in Linear Algebra}. Translations of Mathematical Monographs \textbf{134}, AMS, Providence, RI, (1994), 225 pp.


\bibitem[R1]{R1}
Rudakov A.N., ``Irreducible representations of infinite-dimensional Lie algebras of Cartan type", \textit{Math. USSR Izvestiya} \textbf{38}:4 (1974) 835--866 


\bibitem[Sr]{Sr}
Serganova V., ``Classification of real simple Lie superalgebras and symmetric superspaces", \textit{Funct. Anal. Appl.} \textbf{17}:3 (1983), 200--207; for details of the proof, see \cite[Appendix D5]{Be}.


\bibitem[Sha]{Sha}
Shander V., ``Invariant functions on supermatrices"; \url{http://arxiv.org/pdf/math/9810112}. For an edited translation into Russian, see \cite{Be}, Appendix~D$7$.


\bibitem[Shch]{Shch}
Shchepochkina I., ``How to realize Lie algebras by vector fields", \textit{Theor. and Mathem. Physics} \textbf{147}:3 (2006), 821--838; \url{http://arxiv.org/pdf/math/0509472}


\bibitem[T]{T}
Tanaka N., ``On the equivalence problems associated with simple graded Lie algebras", \textit{Hokkaido Math. J.} \textbf{8}:1 (1979), 23--84


\bibitem[Ver]{Ver}
Verbitsky M.,``Action of the Lie algebra of $\text{SO}(5)$ on the cohomology of a~hyper-K\"ahler manifold", \textit{Functional Anal. Appl.} \textbf{24}:3 (1990), 229--230

\bibitem[V1]{V1}
Verbitsky M., \url{http://verbit.ru/MATH/TALKS/HSE-sp_1_1.pdf}

\bibitem[VPST]{VPST}
Volkov D.~V., Pashnev A.~I., Soroka V.~A., Tkach V.~I., ``Hamiltonian dynamical systems with even and odd Poisson brackets", \textit{Theor. and Mathem. Physics} \textbf{79}:1 (1989), 424--429


\bibitem[Weil]{Weil}
Weil A. \textit{Introduction a~l'\'etude des vari\'et\'es K\"ahleriennes}. Hermann (1958) 175pp.


\bibitem[Wo]{Wo}
Wonenburger M.~J., ``Simultaneous diagonalization of symmetric bilinear forms", \textit{J. of Mathem. and Mech.} \textbf{15} (1966), 617--622


\bibitem[Y]{Y}
Yamaguchi K., ``Differential systems associated with simple graded Lie algebras", \textit{Progress in differential geometry}, Adv. Stud. Pure Math. \textbf{22}, Math. Soc. Japan, Tokyo, (1993), 413--494


\bibitem[ZJ]{ZJ}
Zinn-Justin J., ``Renormalization of gauge theories". In: H.~Rollnik and K.~Dietz (eds.) \textit{Trends in Elementary Particle Physics}, Lecture Notes in Physics \textbf{37} (1975), 1--39; \url{https://arxiv.org/pdf/hep-th/9906115}

\end{thebibliography}
\end{document}